     %%%%%%%%%%%%%%%%%%%%%%%%%%%%%%%%%%%%%%%%%%%%%%%%%%%%%
     %%                                                 %%
     %%         QUANTIZATION OF POISSON GROUPS          %%
     %%                                                 %%
     %%                       by                        %%
     %%                                                 %%
     %%                 Fabio GAVARINI                  %%
     %%                                                 %%
     %%                                                 %%
     %%%%%%%%%%%%%%%%%%%%%%%%%%%%%%%%%%%%%%%%%%%%%%%%%%%%%

         %%%%%%%%%%%%%%%%%%%%%%%%%%%%%%%%%%%%%%%%%%%%%
         %%                                         %%
         %%      final version published on the     %%
         %%                                         %%
         %%     "Pacific Journal of Mathematics"    %%
         %%                                         %%
         %%%%%%%%%%%%%%%%%%%%%%%%%%%%%%%%%%%%%%%%%%%%%

%&amstex
\input amstex
\documentstyle{amsppt}

\magnification=\magstep1
\hsize=6.5truein
\vsize=9truein

\font \smallrm=cmr10 at 10truept
\font \smallbf=cmbx10 at 10truept
\font \smallit=cmti10 at 10truept
\font \smallsl=cmsl10 at 10truept
 at 10truept
 at 10truept

\baselineskip=.15truein

\def \loongrightarrow {\relbar\joinrel\relbar\joinrel\rightarrow}
\def \llongrightarrow {\relbar\joinrel\relbar\joinrel\relbar\joinrel
\rightarrow}

\def \longhookrightarrow {\lhook\joinrel\relbar\joinrel\rightarrow}
\def \llonghookrightarrow
{\lhook\joinrel\relbar\joinrel\relbar\joinrel\relbar\joinrel\rightarrow}
\def \lllonghookrightarrow
{\lhook\joinrel\relbar\joinrel\relbar\joinrel\relbar\joinrel\relbar\joinrel
\relbar\joinrel\relbar\joinrel\rightarrow}
\def \longtwoheadrightarrow {\relbar\joinrel\twoheadrightarrow}
\def \llongtwoheadrightarrow
{\relbar\joinrel\relbar\joinrel\relbar\joinrel\twoheadrightarrow}

\def \smallcirc {\, {\scriptstyle \circ} \,}
\def \aij {a_{ij}}
\def \unon {1,\dots,n}

\def \N {\Bbb N}
\def \Z {\Bbb Z}
\def \Q {\Bbb Q}

\def \qm {q^{-1}}
\def \kq {k(q)}
\def \kqqm {k\!\left[q,\qm\right]}

\def \ebar {\overline E}
\def \fbar {\overline F}
\def \rbar {\overline r}
\def \otimeshat {\,\widehat\otimes\,}

\def \gerg {{\frak g}}
\def \gerh {{\frak h}}
\def \gerk {{\frak k}}
\def \gerK {{\frak K}}
\def \gert {{\frak t}}
\def \gern {{\frak n}}
\def \gerb {{\frak b}}
\def \gerE {{\frak E}}
\def \gerU {{\frak U}}
\def \gerF {{\frak F}}

\def \calU {\hbox{$ {\Cal U} $}}
\def \calF {\hbox{$ {\Cal F} $}}
\def \calE {\hbox{$ {\Cal E} $}}

\def \ug {U(\gerg)}
\def \ugtau {U(\gerg^\tau)}
\def \uh {U(\gerh)}
\def \uhtau {U(\gerh^\tau)}

\def \uqMg#1 {U_q^{\scriptscriptstyle#1}(\gerg)}
\def \gerUMg#1 {\gerU^{\scriptscriptstyle#1}(\gerg)}
\def \gerUunoMg#1 {\gerU_1^{\scriptscriptstyle #1}(\gerg)}
\def \gerUepsilonMg#1 {\gerU_\varepsilon^{\scriptscriptstyle #1}(\gerg)}
\def \calUMg#1 {{\Cal U}^{\scriptscriptstyle #1}(\gerg)}

\def \uzM#1 {U_0^{\scriptscriptstyle#1}}
\def \calUzM#1 {{\Cal U}_0^{\scriptscriptstyle#1}}
\def \gerUzM#1 {\gerU_0^{\scriptscriptstyle#1}}

\def \calUp {{\Cal U}_{\scriptscriptstyle +}}
\def \gerUp {\gerU_{\scriptscriptstyle +}}

\def \uMbm#1 {U^{\scriptscriptstyle#1}_{\scriptscriptstyle  \leq}}
\def \calUMbm#1 {{\Cal U}^{\scriptscriptstyle#1}_{\scriptscriptstyle \leq}}
\def \gerUMbm#1 {\gerU^{\scriptscriptstyle#1}_{\scriptscriptstyle \leq}}

\def \uMbp#1 {U^{\scriptscriptstyle#1}_{\scriptscriptstyle \geq}}
\def \calUMbp#1 {{\Cal U}^{\scriptscriptstyle#1}_{\scriptscriptstyle \geq}}
\def \gerUMbp#1 {\gerU^{\scriptscriptstyle#1}_{\scriptscriptstyle \geq}}

\def \uqMh#1 {U_q^{\scriptscriptstyle#1}(\gerh)}
\def \calUMh#1 {{\Cal U}^{\scriptscriptstyle #1}(\gerh)}

\def \gerUMh#1 {\gerU^{\scriptscriptstyle#1}(\gerh)}
\def \gerUunoMh#1 {\gerU_1^{\scriptscriptstyle #1}(\gerh)}
\def \gerUepsilonMh#1 {\gerU_\varepsilon^{\scriptscriptstyle #1}(\gerh)}

\def \fg {F[G]}
\def \fh {F[H]}
\def \finfg {F^\infty[G]}

\def \fqMg#1 {F_q^{\scriptscriptstyle#1} [G]}
\def \fqMinfg#1 {F_q^{\scriptscriptstyle#1, \infty} [G]}
\def \gerFMg#1 {{\frak F}^{\scriptscriptstyle#1} [G]}
\def \gerFMinfg#1 {{\frak F}^{\scriptscriptstyle #1, \infty} [G]}
\def \gerFunoMg#1 {{\frak F}_1^{\scriptscriptstyle #1} [G]}

\def \calFMg#1 {{\Cal F}^{\scriptscriptstyle #1} [G]}
\def \calFMinfg #1{{\Cal F}^{\scriptscriptstyle #1, \infty} [G]}
\def \calFunoMg#1 {{\Cal F}_1^{\scriptscriptstyle #1} [G]}
\def \calFepsilonMg#1 {{\Cal F}_\varepsilon^{\scriptscriptstyle #1} [G]}
\def \calFunoMinfg#1 {{\Cal F}_1^{\scriptscriptstyle #1, \infty} [G]}

\def \fqMbs #1#2 {{F_q^{{\scriptscriptstyle#1}} [B_{#2}]}}
\def \gerFMbs #1#2 {{{\frak F}^{{\scriptscriptstyle#1}} [B_{#2}]}}

\def \uqphiMg#1 {U_{q,\varphi}^{\scriptscriptstyle#1}(\gerg)}
\def \gerUphiMg#1 {\gerU_{\varphi}^{\scriptscriptstyle#1}(\gerg)}
\def \gerUunophiMg#1 {\gerU_{1,\varphi}^{\scriptscriptstyle #1} (\gerg)}
\def \gerUepsilonphiMg#1 {\gerU_{\varepsilon, \varphi}^{\scriptscriptstyle
#1} (\gerg)}\def \calUphiMg#1 {{\Cal U}_\varphi^{\scriptscriptstyle #1}(\gerg)}
\def \calUunophiMg#1 {{\Cal U}_{1,\varphi}^{\scriptscriptstyle #1}(\gerg)}
\def \calUepsilonphiMg#1 {{\Cal U}_{\varepsilon, \varphi}^{\scriptscriptstyle
#1}(\gerg)}

\def \uphizM#1 {U_{\varphi,0}^{\scriptscriptstyle#1}}
\def \calUphizM#1 {{\Cal U}_{\varphi,0}^{\scriptscriptstyle#1}}
\def \gerUphizM#1 {\gerU_{\varphi,0}^{\scriptscriptstyle#1}}

\def \uphim {U_{\varphi,{\scriptscriptstyle -}}}
\def \uphip {U_{\varphi,{\scriptscriptstyle +}}}
\def \calUphim {{\Cal U}_{\varphi,{\scriptscriptstyle -}}}
\def \gerUphim {\gerU_{\varphi,{\scriptscriptstyle -}}}
\def \calUphip {{\Cal U}_{\varphi,{\scriptscriptstyle +}}}
\def \gerUphip {\gerU_{\varphi,{\scriptscriptstyle +}}}

\def \uphiMbs #1#2 {U^{\scriptscriptstyle#1}_{\varphi,{\scriptscriptstyle
#2}}}
\def \gerUphiMbs #1#2 {\gerU^{\scriptscriptstyle#1}_{\varphi,
{\scriptscriptstyle #2}}}
\def \calUphiMbs #1#2 {{\Cal U}^{\scriptscriptstyle#1}_{\varphi,
{\scriptscriptstyle #2}}}

\def \uphiMbm#1 {U^{\scriptscriptstyle#1}_{\varphi,{\scriptscriptstyle
\leq}}}
\def \calUphiMbm#1 {{\Cal U}^{\scriptscriptstyle#1}_{\varphi,
{\scriptscriptstyle \leq}}}
\def \gerUphiMbm#1 {\gerU^{\scriptscriptstyle#1}_{\varphi,
{\scriptscriptstyle \leq}}}

\def \uphiMbp#1 {U^{\scriptscriptstyle#1}_{\varphi, {\scriptscriptstyle
\geq}}}
\def \calUphiMbp#1 {{\Cal U}^{\scriptscriptstyle#1}_{\varphi,
{\scriptscriptstyle \geq}}}
\def \gerUphiMbp#1 {\gerU^{\scriptscriptstyle#1}_{\varphi,
{\scriptscriptstyle \geq}}}

\def \uqphiMh#1 {U_{q,\varphi}^{\scriptscriptstyle#1} (\gerh)}
\def \calUphiMh#1 {{\Cal U}_\varphi^{\scriptscriptstyle#1}(\gerh)}
\def \calUunophiMh#1 {{\Cal U}_{1,\varphi}^{\scriptscriptstyle #1} (\gerh)}
\def \calUepsilonphiMh#1 {{\Cal U}_{\varepsilon, \varphi}^{\scriptscriptstyle
#1} (\gerh)}
\def \gerUphiMh#1 {\gerU_\varphi^{\scriptscriptstyle#1}(\gerh)}
\def \gerUunophiMh#1 {\gerU_{1,\varphi}^{\scriptscriptstyle #1}(\gerh)}
\def \gerUepsilonphiMh#1 {\gerU_{\varepsilon, \varphi}^{\scriptscriptstyle
#1}(\gerh)}

\def \uqphiMinfh#1 {U_{q,\varphi}^{\scriptscriptstyle #1, \infty} (\gerh)}
\def \calUphiMinfh#1 {{\Cal U}_\varphi^{\scriptscriptstyle#1, \infty}
(\gerh)}
\def \calUunophiMinfh#1 {{\Cal U}_{1,\varphi}^{\scriptscriptstyle #1, \infty}
(\gerh)}
\def \calUepsilonphiMinfh#1 {{\Cal U}_{\varepsilon,
\varphi}^{\scriptscriptstyle #1, \infty} (\gerh)}
\def \gerUphiMinfh#1 {\gerU_\varphi^{\scriptscriptstyle#1, \infty} (\gerh)}
\def \gerUunophiMinfh#1 {\gerU_{1,\varphi}^{\scriptscriptstyle #1, \infty}
(\gerh)}
\def \gerUepsilonphiMinfh#1 {\gerU_{\varepsilon, \varphi}^{\scriptscriptstyle
#1, \infty} (\gerh)}

\def \fgtau {F\left[G^\tau\right]}
\def \fhtau {F\left[H^\tau\right]}
\def \finfgtau {F^\infty \left[G^\tau\right]}

\def \fqphiMg#1 {F_{q,\varphi}^{\scriptscriptstyle#1} [G]}
\def \fqphiMinfg#1 {F_{q, \varphi}^{\scriptscriptstyle#1, \infty} [G]}
\def \gerFphiMg#1 {{\frak F}_{\varphi}^{\scriptscriptstyle#1} [G]}
\def \gerFphiMinfg#1 {{\frak F}_{\varphi}^{\scriptscriptstyle#1,
\scriptscriptstyle \infty} [G]}
\def \gerFunophiMg#1 {{\frak F}_{1,\varphi}^{\scriptscriptstyle #1} [G]}
\def \gerFepsilonphiMg#1 {{\frak F}_{\varepsilon,
\varphi}^{\scriptscriptstyle #1} [G]}
\def \gerFunophiMinfg#1 {{\frak F}_{1,\varphi}^{\scriptscriptstyle M, \infty}
[G]}
\def \gerFepsilonphiMinfg#1 {{\frak F}_{\varepsilon,
\varphi}^{\scriptscriptstyle #1,\infty} [G]}
\def \calFphiMg#1 {{\Cal F}_{\! \varphi}^{\scriptscriptstyle#1} [G]}
\def \calFphiMinfg#1{{\Cal F}_{\! \varphi}^{\scriptscriptstyle#1,
\scriptscriptstyle \infty} [G]}
\def \calFunophiMg#1 {{\Cal F}_{1, \varphi}^{\scriptscriptstyle #1} [G]}
\def \calFepsilonphiMg#1 {{\Cal F}_{\varepsilon, \varphi}^{\scriptscriptstyle
#1}
[G]} \def \calFunophiMinfg#1 {{\Cal F}_{1,\varphi}^{\scriptscriptstyle #1,
\infty} [G]}

\def \fqphiMbs#1#2 {F_{q,\varphi}^{\scriptscriptstyle#1} [B_{#2}]}
\def \gerFphiMbs#1#2 {{\frak F}_\varphi^{\scriptscriptstyle#1} [B_{#2}]}
\def \calFphiMbs#1#2 {{\Cal F}_{\! \varphi}^{\scriptscriptstyle#1} [B_{#2}]}

\def \fqphiMbm #1 {F_{q,\varphi}^{\scriptscriptstyle#1} [B_-]}
\def \gerFphiMbm #1 {{\frak F}_\varphi^{\scriptscriptstyle#1} [B_-]}
\def \calFphiMbm #1 {{\Cal F}_{\! \varphi}^{\scriptscriptstyle #1} [B_-]}
\def \fqphiMbp #1 {F_{q,\varphi}^{\scriptscriptstyle#1} [B_+]}
\def \gerFphiMbp #1 {{\frak F}_\varphi^{\scriptscriptstyle#1} [B_+]}
\def \calFphiMbp #1 {{\Cal F}_{\! \varphi}^{\scriptscriptstyle #1} [B_+]}

\def \gerFrgtau {{\frak F}{\frak r}_{\gerg^\tau}}
\def \gerFrhtau {{\frak F}{\frak r}_{\gerh^\tau}}
\def \calFrgtau {{{\Cal F}r}_{\gerg^\tau}}
\def \calFrhtau {{{\Cal F}r}_{\gerh^\tau}}
\def \gerFrGtau {{\frak F}{\frak r}_{G^\tau}}
\def \calFrHtau {{{\Cal F}r}_{H^\tau}}

%\NoBlackBoxes

\document

\topmatter

{\ }

\vskip-51pt

 \centerline{ \smallrm {\smallsl Pacific Journal of Mathematics\/}  {\smallbf 186},  no.~2 (1998), 217--266 }
 \vskip1pt
 \centerline{ \smallrm  http://pjm.berkeley.edu/pjm/1998/186-2/pjm-v186-n2-p02-p.pdf }
% 
% \hfill   {{\sl Pacific Journal of Mathematics\/}  {\bf 186},
% no.~2 (1998), 217--266}
% 

\vskip39pt

 \title
  Quantization of Poisson groups
 \endtitle

\author
  Fabio Gavarini
\endauthor

\affil
   Dipartimento di Matematica, Istituto G. Castelnuovo,\\ Universit\`a
degli studi di Roma "La Sapienza" \\
\endaffil

\address\hskip-\parindent
            Dipartimento di Matematica  \newline
            Universit\`a di Roma "La Sapienza"  \newline
            Piazzale Aldo Moro 5  \newline
            I-00185 Roma --- ITALY  \newline
            e-mail:  \  gavarini\@mat.uniroma1.it  \newline
  \phantom{e-mail:}  \  gavarini\@mat.uniroma3.it
\endaddress

 \abstract
   Let  $ G^\tau $  be a connected simply connected semisimple algebraic
group, endowed with generalized Sklyanin-Drinfel'd structure of Poisson
group; let  $ H^\tau $  be its dual Poisson group.  By means of quantum
double construction and dualization via formal Hopf algebras, we construct
new quantum groups  $ {\uqphiMh M } $   --- dual of  $ {\uqphiMg {M'} } $
---   which yield infinitesimal quantization of  $ H^\tau $  and  $ G^\tau
\, $;  we study their specializations at roots of 1 (in particular, their
classical limits), thus discovering new quantum Frobenius morphisms. The
whole description dualize for  $ H^\tau $  what was known for  $ G^\tau $,
completing the quantization of the pair  $ (G^\tau,H^\tau) $.
 \endabstract

\endtopmatter

\footnote""{ 1991 {\it Mathematics Subject Classification,}
Primary 17B37, 81R50 }

\vskip0,5truecm

 \centerline{ \bf  Introduction }

\vskip10pt

\hfill  \hbox{\vbox{ \hbox{\it  \   "Dualitas dualitatum }
                     \hbox{\it \ \ \;\, et omnia dualitas" }
                     \vskip4pt
                     \hbox{\sl    N.~Barbecue, "Scholia" } } \hskip1truecm }

\vskip8pt

   Let  $ G $  be a semisimple, connected, and simply
connected affine algebraic group over an algebraically closed field  $ k $
of characteristic zero; we consider a family of structures of Poisson group on
$ G $,  indexed by a multiparameter  $ \tau $,  which generalize the
Sklyanin-Drinfel'd one.  Then every such
Poisson group  $ G^\tau $  has a dual Poisson
group  $ H^\tau $,  and  $ \, \gerg^\tau := Lie(G^\tau) \, $  and  $ \,
\gerh^\tau := Lie(H^\tau) \, $  are Lie bialgebras dual of each other.
                                                      \par
  In 1985 Drinfel'd and Jimbo provided a quantization
of  $ \, \ug = U(\gerg^0) \, $,  namely a Hopf algebra  $ {\uqMg Q } $  over
$ \kq $,  presented by generators and relations, with a
$ \kqqm $--form  $ {\gerUMg Q } $  which for  $ \, q \rightarrow 1 \, $
specializes to  $ \ug $  as a Poisson Hopf coalgebra.  This has been
extended to general parameter  $ \tau $
introducing multiparameter quantum groups  $ {\uqphiMg Q } $  (cf.~[R],
[CV-1], [CV-2]).  Dually, one
constructs a Hopf algebra  $ {\fqMg P } $  of matrix coefficients of
$ {\uqMg Q } $  with a  $ \kqqm $--form  $ {\gerFMg P } $  which specializes
to  $ \fg $,  as a Poisson Hopf algebra, for  $ \, q \rightarrow 1 \, $;  in
particular  $ {\gerFMg P } $  is nothing but the Hopf
subalgebra of "functions" in  $ {\fqMg P } $  which take values in  $ \kqqm $
when "evaluated" on  $ {\gerUMg Q } $  (in a word, the  $ \kqqm $--integer
valued functions on  $ {\gerUMg Q } $).  This again extends to general
$ \tau $  (cf.~[CV-2]).
                                                 \par
  So far the quantization only dealt with the Poisson group  $ G $
(or  $ G^\tau \, $);  the dual group  $ H $  is involved defining a different
$ \kqqm $--form  $ {\calUMg P } $  (of a quantum group  $ {\uqMg P } $)
which specializes to  $ \fh $  (as a Poisson Hopf algebra) for  $ \, q
\rightarrow 1 \, $  (cf.~[DP]), with generalization to the multiparameter
case possible again.  Here sort of a "mixing dualities" (Hopf duality   ---
among enveloping and function algebra ---   and Poisson duality   --- among
dual Poisson groups) occurs, which was described (in a formal setting) by
Drinfel'd (cf.~[Dr], \S 7), and by  Etingof and Kazhdan (cf.~[EK-1], [EK-2]).
This leads to consider the following: let  $ {\fqMg Q } $  be the
quantum function algebra dual of  $ {\uqMg P } $,  and look at the "dual" of
$ {\calUMg P } $  within  $ {\fqMg Q } $,  call it  $ {\calFMg Q } $,
namely the Hopf algebra of  $ \kqqm $--integer  valued functions on
$ {\calUMg P } $;  then this should specialize to  $ \uh $  (as a Poisson
Hopf coalgebra) for  $ \, q \rightarrow 1 \, $;  the same conjecture can be
formulated in the multiparameter case too.
                                                   \par
   Our starting aim was to achieve this
goal, i.~e.~to construct  $ {\fqMg Q } $  and its
$ \kqqm $--form  $ {\calFMg Q } $,  and to prove that
$ {\calFMg Q } $  is a deformation of the Poisson Hopf coalgebra  $ \uh $.
This goal is succesfully attained by performing a suitable dualization of the
quantum double construction; but by the way, this leads to discover a  {\it new
quantum group},  which we call  $ {\uqMh M } $,  which is for  $ \uh $  what
$ {\uqMg M } $  is for  $ \ug $;  in particular it has an integer form
$ {\gerUMh Q } $  which is a quantization of  $ \uh $,  and an integer form
$ {\calUMh P } $  which is a quantization of  $ \finfg $  (the function algebra
of the formal Poisson group associated to  $ G \, $).  Furthermore, we exhibit a
Hopf pairing between  $ {\uqMg {M'} } $  and  $ {\uqMh M } $  which gives a
quantization of the various pairings occurring among the algebras attached to
the pair  $ (G,H) $.  Once again all this extends to the multiparameter case.
  Thus in particular we provide a (infinitesimal) quantization for a wide
class of Poisson groups (the  $ H^\tau $'s);  now, in the
summer of 1995 (when the present work was already
accomplished) a quantization of any Poisson group was
presented in [EK-1] and [EK-2]; but greatest generality
implies lack of concreteness: in contrast, our construction
is extremely concrete; moreover, it allows specialization
at roots of 1, construction of quantum Frobenius morphisms, and so on (like
for  $ \gerUMg{Q} $  and  $ \calUMg{P} $),  which is not possible in the
approach of [EK-1], [EK-2].
                                                      \par
  Finally, a brief sketch of the main ideas of the paper.  Our aim is to
study the "dual" of a quantum group  $ {\uqphiMg M } $  ($ M $  being a
lattice of weights).  First, we select as operation of
"dualization" the most na\"\i{}ve one, namely taking the  {\it full linear
dual\/}  (rather than the usual   --- restricted ---   Hopf dual), the latter
being a  {\it formal\/}  Hopf algebra (rather than a common Hopf algebra).
Second, as  $ {\uqphiMg M } $  is a quotient of a quantum double  $ \,
D_{q,\varphi}^{\scriptscriptstyle M} (\gerg) := D \left( U_{q,
\varphi}^{\scriptscriptstyle M} (\gerb_-), U_{q,\varphi}^{\scriptscriptstyle
M} (\gerb_+), \pi_\varphi \right) \, $,  its linear dual  $ {\uqphiMg M }^* $
embeds into  $ {D_{q,\varphi}^{\scriptscriptstyle M}
(\gerg)}^* $.  Third, since  $ \, D_{q,\varphi}^{\scriptscriptstyle M}
(\gerg) \cong U_{q,\varphi}^{\scriptscriptstyle M} (\gerb_+) \otimes
U_{q,\varphi}^{\scriptscriptstyle M} (\gerb_-) \, $  (as coalgebras) we have
$ \, {D_{q,\varphi}^{\scriptscriptstyle M} (\gerg)}^* \cong {U_{q,
\varphi}^{\scriptscriptstyle M} (\gerb_+)}^* \otimeshat {U_{q,
\varphi}^{\scriptscriptstyle M} (\gerb_-)}^* \, $  (as algebras), where
$ \, \otimeshat \, $  denotes topological tensor product.  Fourth, since
quantum Borel algebras of opposite sign are perfectly paired their linear
duals are suitable completions of quantum Borel algebras again: thus we find
a presentation of  $ {\uqphiMg M }^* $  by generators and relations which
leads us to  {\sl define}  $ \, {\uqphiMh M } := {\uqphiMg{M'\!} }^* \, $
(where  $ M' $  depends on  $ M \, $)  and gives us all claimed results;
because of their construction we call the new objects  $ {\uqphiMh M } $
{\it (multiparameter) quantum formal groups}.
                                                      \par
   In contrast, we also present an alternative
approach, yielding other  {\sl new}  objects   --- denoted by
$ {\fqphiMinfg M } $ ---   which we call  {\it (multiparameter) formal
quantum groups};  the  {\sl similar}  but  {\sl different}  terminology
reveals the fact that  $ {\uqphiMh M } $  and  $ {\fqphiMinfg M } $  provide
two  {\sl different}  quantizations of the  {\sl same}  classical objects
$ \uhtau $  and  $ \finfgtau $,  arising from two different ways of realizing
$ \finfgtau $.

\vskip13pt

\centerline { ACKNOWLEDGEMENTS }

\vskip4pt

  The author wishes to thank C.~Procesi for several fruitful
conversations, and C.~De Concini for some useful talk; he is also endebted
with M.~Costantini, M.~Varagnolo, and I.~Damiani for many helpful
discussions.

\vskip1,7truecm

\centerline{ \bf  \S \; 1 \,  The classical objects }

\vskip10pt

   {\bf 1.1  Cartan data.}  \  Let  $ A:= {(\aij)}_{i,j=\unon} $  be a
$ n \times n $   symmetrizable Cartan matrix; thus  $ \, \aij \in \Z
\, $,  $ \, a_{ii} = 2 $,  $ a_{ij} \leq 0 \, $  if  $ i \neq
j $,  and there exists a vector  $ (d_1, \ldots, d_n) $
with relatively prime positive integral entries  $ d_i $
such that  $ \, (d_i \aij)_{i,j=\unon;} \, $  is a symmetric
positive definite matrix.  Define the  {\it weight lattice\/}  $ P $  to be
the lattice with basis  $ \{ \omega_1, \dots, \omega_n \} $;  let  $ \, P_+
:= \sum^n_{i=1} \N \omega_i \, $  be the subset of  {\it dominant integral
weights\/},  $ \, \alpha_j := \sum^n_{i=1} \aij \omega_i \; (j=\unon) \, $
the  {\it simple roots\/},  $ \, Q := \sum^n_{j=1} \Z \alpha_j \; (\,\subset P
\,) \, $  the  {\it root lattice\/},  and  $ \, Q_+ := \sum^n_{j=1} \N \alpha_j
\, $  the  {\it positive root lattice\/}.  Let  $ W $  be the Weyl group
associated to  $ A $,  and let  $ \, \Pi := \{ \alpha_1, \dots,
\alpha_n \} \, $:  then  $ \, R := W \big( \Pi \big) \, $  is the set of
{\it roots\/},  $ \, R^+ = R \cap Q_+ \, $  the set of {\it positive
roots\/};  finally, we set $ N:= \#(R^+) \; (\, = \vert W \vert) $.  Define
bilinear pairings  $ \, \langle \ \vert \ \rangle \colon Q \times P \rightarrow
\Z \, $  and  $ \, (\ \vert \ ) \colon Q \times P \rightarrow \Z \, $  by
$ \langle \alpha_i \vert \omega_j \rangle = \delta_{ij} $  and
$ (\alpha_i \vert \omega_j) = \delta_{ij} d_i $.  Then
$ (\alpha_i \vert \alpha_j) = d_i \aij $,  giving a
symmetric  $ \Z $--valued  $ W $--invariant  bilinear form on  $ Q $  such
that  $ (\alpha \vert \alpha) \in 2 \Z $.  For all  $ \, \alpha \in R^+ $,  let
$ \, d_\alpha := {\, (\alpha \vert \alpha) \, \over \, 2 \,} \, $;  then  $ \,
d_{\alpha_i} = d_i \, $  for all  $ \, i= 1, \dots, n \, $.  We also extend the
$ \, (\ \vert \ ) \colon Q \times P \rightarrow \Z \, $  to a
(non-degenerate) pairing  $ \, (\ \vert \ ) \colon \Q Q \times
\Q P \rightarrow \Q \, $  of  $ \Q $--vector  spaces by scalar extension,
where  $ \, \Q T := \Q \otimes_\Z T (T = Q, P) \, $:  then restriction gives
a pairing  $ \, (\ \vert \ ) \colon P \times P \rightarrow \Q \, $  (looking
at  $ P $  as a sublattice of  $ \, \Q Q \, $),  which takes values in
$ \Z \left[ D^{-1} \right] $,  where  $ \, D:= det \left( {(\aij)}_{i,j=1}^n
\right) \, $.  Given any pair of lattices  $ (M,M') $,  with  $ \, Q \leq
M, M' \leq \Q P \, $,  we say that they are  {\it dual of each  other\/}  if
$ \, M' = \big\{\, y \in \Q P \,\big\vert\, \langle M, y \rangle \subseteq
\Z \,\big\} \, $,  $ \, M = \big\{\, x \in \Q P \,\big\vert\, \langle x,
M' \rangle \subseteq \Z \,\big\} \, $,  the two conditions being equivalent;
then for any lattice  $ M $  with  $ \, Q
\leq M \leq \Q P \, $  there exists a unique dual lattice  $ M' $  such that
$ \, Q \leq M' \leq \Q P \, $  and  $ \, (\ \vert \ ) \colon \Q P
\times \Q P \rightarrow \Q \, $  restricts to a perfect pairing
$ \, (\ \vert \ ) \colon M \times M' \rightarrow \Z \, $;  in particular
$ \, P' = Q \, $  and  $ \, Q' = P \, $.  In the sequel we denote
by  $ \, \{\, \mu_1, \dots, \mu_n \,\} \, $  and  $ \, \{\, \nu_1, \dots,
\nu_n \,\} \, $  fixed  $ \Z $--bases  of  $ M $  and  $ M' $  dual of each
other, i.~e.~such that  $ \, (\mu_i \vert \nu_j) = \delta_{i j} \, $  for all
$ \, i, j= \unon \, $,  and we set  $ \, M_+ := M \cap P_+ \, $.
                                                \par
  In the following our constructions will work in general for the pairs of
dual lattices  $ (P,Q) $  and  $ (Q,P) $;  but in the simply laced case
(in which  $ \, \langle \ ,\ \rangle = (\ ,\ ) \, $)  $ \left( M, M' \right) $
will be  {\sl any}  pair of dual lattices.

\vskip7pt

   {\bf 1.2  The Poisson groups  $ G $  and  $ H $.}  \    Let  $ G $  be a
connected simply-connected semisimple affine algebraic group over an
algebraically closed field   $ k $  of characteristic 0.  Fix a maximal
torus  $ T \leq G $  and opposite Borel subgroups  $ B_\pm $,  with unipotent
subgroups  $ \, U_\pm \, $,  such that  $ \, B_+ \cap B_- = T \, $,  and let
$ \, \gerg := Lie(G) \, $,  $ \, \gert := Lie(T) \, $,  $ \, \gerb_\pm :=
Lie(B_\pm) \, $,  $ \, \gern_\pm := Lie(U_\pm) \, $;  fix also
$ \, \tau := (\tau_1, \ldots, \tau_n) \in Q^n \, $  such that  $ \, (\tau_i,
\alpha_j) = - (\tau_j, \alpha_i) \, $  for all  $ \, i, j= \unon \, $:  when
$ \, \tau = (0, \ldots, 0) \, $  we shall simply skip it throughout.  Set
$ \, K = G \times G \, $,  define  $ \, G^\tau:= G \, $  embedded in  $ K $
as  the diagonal subgroup, and define a second subgroup
  $$  H^\tau := \Big\{\, \big( u_- t_-, t_+ u_+ \big) \,\Big\vert\, u_\pm \in
U_\pm , t_\pm \in T , t_- t_+ \in \exp \big( \gert^\tau \big) \,\Big\}  \quad
\big(\, \leq B_- \times B_+ \leq K \big)  $$
where  $ \, \gert^\tau := \sum_{i=1}^n k \cdot h_{-\alpha_i
+ 2 \tau_i} \oplus h_{\alpha_i + 2 \tau_i} \leq \gert \oplus
\gert \leq \gerg \oplus \gerg = \gerk := Lie(K) \, $;  hence we have
$ \, \gerh^\tau := Lie(H^\tau) = (\gern_-,0) \oplus
\gert^\tau \oplus (0,\gern_+) \, $.  The triple  $ (K,G^\tau,H^\tau) $  is an
algebraic Manin triple (cf.~[DP], \S 11), whose invariant form is defined as
follows: first normalize the Killing form  $ \, (\ ,\ )\, $  on  $ \gerg $
so that short roots have square length 2; then define the form on  $ \, \gerk =
\gerg \oplus
          \gerg \, $  by\footnote{{\it Warning:\/}  beware of the
normalization of the invariant form of  $ \gerk $,  which is
different from [DP]. }
  $$  \big\langle x_1 \oplus y_1 , x_2 \oplus y_2 \big\rangle := {\,1\, \over
\,2\,} \, (y_1,y_2) - {\,1\, \over \,2\,} \, (x_1,x_2) \, .  $$
   \indent   In general, if  $ \big(\gerk', \gerg', \gerh'\big) $  is any
Manin triple, the bilinear form on  $ \gerk' $  gives by restriction a
non-degenerate pairing  $ \, \langle \ , \ \rangle \colon \gerh' \otimes \gerg'
\rightarrow k \, $  which is a pairing of Lie bialgebras, that is
  $$  \big\langle x, [y_1,y_2] \big\rangle = \big\langle \delta (x),
y_1 \otimes y_2 \big\rangle \, ,  \quad  \big\langle [x_1,x_2], y \big\rangle
= \big\langle x_1 \otimes x_2, \delta (y) \big\rangle  $$
where  $ \delta $  is the Lie cobracket;  we shall call it
{\it Poisson pairing\/}.  In the present case we denote it by
$ \, \pi^\tau_{\Cal P}(h,g) := \langle h, g \rangle  \, $;  it is described
by
  $$  \matrix
  \langle \text{f}^{\,\tau}_i, f_j \rangle = 0  &  \langle
\text{f}^{\,\tau}_i, h_j \rangle = 0
&  \langle \text{f}^{\,\tau}_i, e_j \rangle = - {\,1\, \over \,2\,} \,
\delta_{ij} d_i^{-1}  \\
  \langle \text{h}^\tau_i, f_j \rangle = 0  &  \langle
\text{h}^\tau_i, h_j \rangle = a_{ij} d_j^{-1} = a_{ji} d_i^{-1}  &  \langle
\text{h}^\tau_i, e_j \rangle = 0 \\
  \langle\text{e}^\tau_i, f_j \rangle = {\,1\, \over \,2\,} \, \delta_{ij}
d_i^{-1}  &  \langle \text{e}^\tau_i, h_j \rangle = 0  &  \langle
\text{e}^\tau_i, e_j \rangle = 0  \\
\endmatrix   \eqno (1.1)  $$
where the  $ \, \text{f}^{\,\tau}_s $,  $ \text{h}^\tau_s $,
$ \text{e}^\tau_s $,  resp.~$ \, f_s $,  $ h_s $,  $ e_s $,  are
Chevalley-type generators of  $ \gerh^\tau $,  resp.~$ \gerg^\tau $,
embedded inside  $ \, \gerk = \gerg \oplus \gerg \, $,  namely
$ \, \text{f}^{\,\tau}_s = f_s \oplus 0 \, $,  $ \, \text{h}^\tau_s =
h_{-\alpha_s + 2 \tau_s} \oplus h_{\alpha_s + 2 \tau_s} \, $,
$ \, \text{e}^\tau_s = 0 \oplus e_s \, $,  and  $ \, f_s = f_s \oplus f_s
\, $,  $ \, h_s = h_s \oplus h_s \, $,  $ \, e_s  = e_s \oplus e_s \, $
(see \S\S 1.3--4 below).

\vskip7pt

   {\bf 1.3  The Poisson Hopf coalgebra  $ \ugtau $.}  \  The universal
enveloping algebra  $ \, \ugtau = \ug \, $  can be presented as the
associative  \hbox{$k$--alge}bra  with 1 generated
by elements  $ , f_i $,  $ h_i $,  $ e_i $  ($ i=\unon $)  (the
{\it Chevalley generators\/})  satisfying Serre's relations;
it has a canonical structure of Hopf algebra,
given by  $ \, \Delta(x) = x \otimes 1 + 1 \otimes x $,  $ S(x) = -x $,
$ \epsilon(x) = 0 \, $  for  $ \, x = f_i, h_i, e_i \, $;  finally,  the Lie
cobracket  $ \, \delta = \delta_{\gerg^\tau} \colon \gerg^\tau
\longrightarrow \gerg^\tau \otimes \gerg^\tau \, $  extends to a Poisson
cobracket  $ \, \delta \colon \ugtau \longrightarrow \ugtau \otimes
\ugtau \, $  (compatible with the Hopf structure) given by  $ \; \delta (f_i) =
{\,(\alpha_i + 2 \tau_i \vert\alpha_i + 2 \tau_i)\, \over \,2\,} \, h_{\alpha_i
+ 2 \tau} \otimes f_i - {\,(\alpha_i + 2 \tau_i \vert\alpha_i + 2 \tau_i)\,
\over \,2\,} \, f_i \otimes h_{\alpha_i + 2 \tau} \, $,  $ \; \delta (h_i) = 0
\, $,  $ \; \delta (e_i) = {\,(\alpha_i - 2 \tau_i \vert\alpha_i - 2 \tau_i)\,
\over \,2\,} \, h_{\alpha_i - 2 \tau} \otimes e_i - {\,(\alpha_i - 2 \tau_i
\vert\alpha_i - 2 \tau_i)\, \over \,2\,} \, e_i \otimes h_{\alpha_i - 2 \tau}
\, $.

\vskip7pt

   {\bf 1.4  The Poisson Hopf coalgebra  $ \uhtau $.}   \  From the very
definition and the previous presentation of  $ \ugtau $  we get for
$ \uhtau $  the following presentation.  $ \uhtau $  is the
associative  \hbox{$k$--alge}bra  with 1 generated by
$ \text{f}^{\,\tau}_i $,  $ \text{h}^\tau_i $,
$ \text{e}^\tau_i $  ($ i=\unon $)  with relations
  $$  \eqalign {
       \text{h}^\tau_i \text{h}^\tau_j - \text{h}^\tau_j
\text{h}^\tau_i = 0 \, ,  &  \qquad  \text{e}^\tau_i
\text{f}^{\,\tau}_j - \text{f}^{\,\tau}_j \text{e}^\tau_i = 0  \cr
       \text{h}^\tau_i \text{f}^{\,\tau}_j - \text{f}^{\,\tau}_j
\text{h}^\tau_i = \langle \alpha_i - 2 \tau_i ,
\alpha_j \rangle \, \text{f}^{\,\tau}_j \, ,  &  \qquad
\text{h}^\tau_i \text{e}^\tau_j - \text{e}^\tau_j
\text{h}^\tau_i = \langle \alpha_i + 2 \tau_i , \alpha_j \rangle \,
\text{e}^\tau_j  \cr
       \sum_{k=0}^{1-\aij} {(-1)}^k {1 - \aij \choose k}  &
{(\text{f}^{\,\tau}_i)}^{1-\aij-k} \text{f}^{\,\tau}_j
{(\text{f}^{\,\tau}_i)}^k = 0  \qquad  (i \neq j)  \cr
       \sum_{k=0}^{1-\aij} {(-1)}^k {1 - \aij \choose k}  &
{(\text{e}^\tau_i)}^{1-\aij-k} \text{e}^\tau_j
{(\text{e}^\tau_i)}^k = 0  \qquad  (i \neq j)  \cr }
   \eqno (1.2)  $$
for all  $ \, i, j= \unon \, $;  its natural Hopf structure is given by
  $$  \Delta(x) = x \otimes 1 + 1 \otimes x \, ,  \qquad  S(x) = -x \, ,
\qquad  \epsilon(x) = 0   \eqno (1.3)  $$
for  $ \, x= \text{f}^{\,\tau}_i, \text{h}^\tau_i, \text{e}^\tau_i $,  and the
co-Poisson structure  $ \, \delta = \delta_{\gerh^\tau}
\colon \uhtau \longrightarrow \uhtau \otimes \uhtau \, $  by
  $$  \displaylines{
   \delta(\text{f}^{\,\tau}_i) = d_i \cdot \big( \text{h}^{\,\tau}_i
\otimes \text{f}^{\,\tau}_i - \text{f}^{\,\tau}_i \otimes \text{h}^{\,\tau}_i
\big) + 2 \, d_i^{\,-1} \cdot \!\! \sum_{\alpha, \beta \in R^+} \!\!
c^{i,+}_{\alpha,\beta} \, d_\alpha d_\beta \cdot \big( \text{e}^\tau_\alpha
\otimes \text{f}^{\,\tau}_\beta - \text{f}^{\,\tau}_\beta \otimes
\text{e}^\tau_\alpha \big)  \cr
   \hfill   \delta(\text{h}^{\,\tau}_i) = 4 \, d_i^{-1} \cdot
\sum_{\gamma \in R^+} d_\gamma \, (\gamma \vert \alpha_i)
\cdot \big( \text{e}^\tau_\gamma \otimes \text{f}^{\,\tau}_\gamma -
\text{f}^{\,\tau}_\gamma \otimes \text{e}^\tau_\gamma \big)   \hfill  (1.4)  \cr
   \delta(\text{e}^\tau_i) = d_i \cdot \big( \text{e}^\tau_i
\otimes \text{h}^{\,\tau}_i - \text{h}^{\,\tau}_i \otimes \text{e}^\tau_i
\big) + 2 \, d_i^{\,-1} \cdot \!\! \sum_{\alpha, \beta \in R^+} \!\!
c^{i,-}_{\alpha,\beta} \, d_\alpha d_\beta \cdot \big( \text{f}^{\,\tau}_\beta
\otimes \text{e}^\tau_\alpha - \text{e}^\tau_\alpha \otimes
\text{f}^{\,\tau}_\beta \big)  \cr }  $$
with the  $ \text{e}^\tau_\gamma $'s  and the  $ \text{f}^{\,\tau}_\gamma $'s
given by  $ \, \big\langle \text{e}^\tau_\gamma, f_\eta \big\rangle = +
\delta_{\gamma,\eta} d_\gamma \big/ 2 \, $,  $ \, \big\langle
\text{e}^\tau_\gamma, h_i \big\rangle = 0 \, $,  $ \, \big\langle
\text{e}^\tau_\gamma, e_\eta \big\rangle = 0 \, $,  and  $ \; \big\langle
\text{f}^{\,\tau}_\gamma, f_\eta \big\rangle = 0 \, $,  $ \, \big\langle
\text{f}^{\,\tau}_\gamma, h_i \big\rangle = 0 \, $,  $ \, \big\langle
\text{f}^{\,\tau}_\gamma, e_\eta \big\rangle = - \delta_{\gamma,\eta} d_\gamma
\big/ 2 \, $  ($ f_\eta $  and  $ e_\eta $  being root vectors in
$ \gerg^\tau $),  and the  $ c^{i,\pm}_{\alpha,\beta} $'s  given by  $ \, \big[
f_\alpha, e_\beta \big] = c^{i,-}_{\alpha,\beta} \cdot f_i \, $,  $ \, \big[
f_\alpha, e_\beta \big] = c^{i,+}_{\alpha,\beta} \cdot e_i \, $.

\vskip1,7truecm

\centerline{ \bf  \S \; 2 \,  Quantum Borel algebras and DRT pairings }

\vskip10pt

   {\bf 2.1  Notations.}  \  For all
$ \, s, n \in \N \, $,  let  $ \, {(n)}_q := {q^n - 1 \over q - 1}
\; (\in k[q]) \, $,  $ \, {(n)}_q! := \prod_{r=1}^n {(r)}_q $,
$ {({n \atop s})}_q := {{(n)}_q! \over {(s)}_q! {(n-s)}_q! } \;
(\in k[q]) \, $,  and  $ \, {[n]}_q := {q^n - q^{-n} \over q - \qm}
\; (\in \kqqm) \, $,  $ \, {[n]}_q! := \prod_{r=1}^n {[r]}_q $,
                            $ {[{n \atop s}]}_q := $\break
$ := {{[n]}_q! \over {[s]}_q! {[n-s]}_q! } \; (\in \kqqm)
\, $;  let  $ \, q_\alpha := q^{d_\alpha} \, $  for all  $ \, \alpha
\in R^+ \, $,  and  $ \, q_i := q_{\alpha_i} \, $.  Let  $ Q $,  $ P $
be as in \S 1; we fix an endomorphism  $ \varphi $  of the  $ \Q $--vector
space  $ \, \Q P := \Q \otimes_\Z P \, $  which is antisymmetric  --- with
respect to  $ (\ \mid \ ) $ ---  and satisfies the conditions
  $$  \varphi(Q) \subseteq Q \, ,  \quad  {\,1\, \over \,2\,} (\varphi(P) \mid
P) \subseteq \Z \, ,  \quad  2 A Y A^{-1} \in {Mat}_n(\Z)  $$
where, letting  $ \, \tau_i := {\,1\, \over \,2\,} \, \varphi(\alpha_i) =
\sum_{j=1}^n y_{ji} \alpha_j \, $,  we set  $ \, Y := {(y_{ij})}_{i,j=
\unon;} \, $.  We also define  $ \, \tau_\alpha := {\,1\, \over \,2\,} \,
\varphi(\alpha) \, $  for all  $ \, \alpha \in R \, $  (so
$ \, \tau_{\alpha_i} = \tau_i \, $).  It is proved in [CV-1] that
$ \, (\text{id}_{\Q P} + \varphi) \, $  and  $ \, (\text{id}_{\Q P} -
\varphi) \, $  are isomorphisms:  then we set
$ \, r := {(\text{id}_{\Q P} + \varphi)}^{-1} \, $,  $ \, \rbar :=
{(\text{id}_{\Q P} - \varphi)}^{-1} \, $.

\vskip7pt

   {\bf 2.2  Quantum Borel algebras.}  \  From now on  $ M $  will be any
lattice such that  $ \, Q \leq M \leq P \, $;  then  $ M' $  will be the dual
lattice defined in \S 1.1, according to the conditions therein.  As in [CV-1],
$ \, U_{q,\varphi}^{\scriptscriptstyle M} (\gerb_-) \, $,  resp.~$ \,
U_{q,\varphi}^{\scriptscriptstyle M} (\gerb_+) \, $,  is
the associative  $ \kq $--algebra  with 1 generated by  $ \, L_\mu $
($ \mu \in M \, $),  $ F_1, \dots, F_n $,  resp.~$ \, L_\mu $
($ \mu \in M \, $),  $ E_1, \dots, E_n $,  with  relations
  $$  \eqalign{
      L_0 = 1 \, ,   &  \qquad  L_\mu L_\nu = L_{\mu+\nu} \, ,  \cr
      L_\mu F_j = q^{-( \alpha_j | \mu )} F_j L_\mu \, ,  &
\qquad \sum_{p+s=1-\aij} {(-1)}^s {\left[ {1 - \aij \atop s} \right]}_{\! q_i}
F_i^p F_j F_i^s = 0  \cr
      \hbox{resp. \ \ \ }  L_\mu E_j = q^{+( \alpha_j | \mu )}
E_j L_\mu \, ,  &  \qquad  \sum_{p+s=1-\aij} {(-1)}^s {\left[ {1 - \aij \atop
s} \right]}_{\! q_i} E_i^p E_j E_i^s = 0  \,  \cr }   \eqno (2.1)  $$
for all  $ \, i, j= \unon $  and  $ \, \mu, \nu \in M  \, $;
both are Hopf algebras, with
  $$  \matrix
   \Delta_\varphi (F_i) = F_i \otimes L_{-\alpha_i - \tau_i}
+ L_{\tau_i} \otimes F_i \, ,  &  \quad  \epsilon_\varphi (F_i) = 0 \, ,  &
\quad  S_\varphi (F_i) = - F_i L_{\alpha_i}  \\
   \Delta_\varphi (L_\mu) = L_\mu \otimes L_\mu \, ,  &  \quad  \epsilon_\varphi
(L_\mu) = 1 \, ,  &  \quad  S_\varphi (L_\mu) = L_{-\mu}  \\
   \Delta_\varphi (E_i) = E_i \otimes L_{\tau_i} +
L_{\alpha_i - \tau_i} \otimes E_i \, ,  &  \quad  \epsilon_\varphi (E_i) = 0
\, ,  &  \quad  S_\varphi (E_i) = - L_{-\alpha_i} E_i  \\
      \endmatrix  $$
for all  $ \, i=\unon \, $,  $ \, \mu \in M \, $.  We also consider the
subalgebras  $ U_{q,\varphi}^{\scriptscriptstyle M} (\gert) $  (generated by
the  $ L_\mu $'s),  $ U_{q,\varphi} (\gern_-) $  (generated by the
$ F_i $'s),  $ U_{q,\varphi} (\gern_+) $  (generated by the  $ E_i $'s).  In
the sequel we shall use the notation  $ \, K_\alpha := L_\alpha \, $,  $ \,
M_\mu := L_\mu \, $,  $ \, \Lambda_\nu := L_\nu \, (\forall \, \alpha \in Q,
\mu \in M, \nu \in M' \,) \, $  (and in particular  $ \, K_i := K_{\alpha_i}
\, $,  $ \, M_i := M_{\mu_i} \, $  $ \, \Lambda_i := \Lambda_{\nu_i} \, $),
and  $ \, \uMbm{M} := U_{q,\varphi}^{\scriptscriptstyle M}(\gerb_-) \, $,
$ \, \uMbp{M} := U_{q,\varphi}^{\scriptscriptstyle M}(\gerb_+) \, $,
$ \, \uphizM{M} := U_{q,\varphi}^{\scriptscriptstyle M} (\gert) \, $,
$ \, U_{\scriptscriptstyle -} := U_{q,\varphi} (\gern_-) \, $,
$ \, U_{\scriptscriptstyle +} := U_{q,\varphi} (\gern_+) \, $.  If
$ \varphi = 0 $  we just skip it througout.
                                                     \par
   Finally, multiplication yields isomorphisms
  $$  {\uphiMbm M } \cong \uphim \otimes {\uphizM M } \cong {\uphizM M } \otimes
\uphim \, ,  \quad  {\uphiMbp M } \cong \uphip \otimes {\uphizM M } \cong
{\uphizM M } \otimes \uphip  $$

\vskip7pt

{\bf 2.3  DRT pairings.}  \;  If  $ H $ is any Hopf algebra, we let
$ H^{op} $  be the same coalgebra with opposite multiplication, and
$ H_{op} $  the same algebra with opposite comultiplication.
                                       \par
   From [CV-1], \S 3, there existe perfect (i.~e.~non-degenerate) pairings
of Hopf algebras
  $$  \eqalign{
   \pi_\varphi \colon \, {\left( {\uphiMbm M } \right)}_{op} \otimes
{\uphiMbp {M'} } \longrightarrow \kq  \, ,  &  \quad \qquad  \pi_\varphi \colon
\, {\uphiMbm M } \otimes {\left( \uphiMbp{M'} \right)}^{op} \longrightarrow
\kq   \cr
   \overline{\pi_\varphi} \colon \, {\left( {\uphiMbp M } \right)}_{op}
\otimes {\uphiMbm {M'} } \longrightarrow \kq  \, ,  &  \quad \qquad
\overline{\pi_\varphi} \colon \, {\uphiMbp M } \otimes {\left(
{\uphiMbm {M'} } \right)}^{op} \longrightarrow \kq  \cr }  $$
  $$  \eqalign{
   \pi_\varphi (L_\mu,L_\nu) = q^{-(r(\mu) \vert \nu)} \, ,  \ \pi_\varphi
(L_\mu,E_j) = 0 \, ,  &  \ \pi_\varphi (F_i,L_\nu) = 0 \, ,  \  \pi_\varphi
(F_i,E_j) = \delta_{ij} {\,q^{-(r(\tau_i) \vert \tau_i)}\, \over
\,\left( q_i^{-1} - q_i \right)\,}  \cr
   \overline{\pi_\varphi} (L_\mu,L_\nu) = q^{+(r(\mu) \vert \nu)} \, , \
\overline{\pi_\varphi} (E_i,L_\nu) = 0 \, ,  &  \  \overline{\pi_\varphi}
(L_\mu,F_j) = 0 \, ,  \  \overline{\pi_\varphi} (E_i,F_j) = \delta_{ij}
{\,q^{+(r(\tau_i) \vert \tau_i)}\, \over \,\left( q_i - q_i^{-1}
\right)\,}  \cr }  $$
   \indent   These pairings were introduced by Drinfel'd, Rosso, Tanisaki,
and others, whence we shall call them   {\it DRT pairings\/}.  If  $ \pi $
is any DRT pairing we shall also set  $ \, {\langle x, y \rangle}_\pi \, $
for  $ \, \pi(x,y) \, $.

\vskip7pt

   {\bf 2.4 PBW bases.}  \  Let  $ \, N:= \#(R^+) \, $,  and fix any total
convex ordering (cf.~[Pa] and [DP], \S 8.2)  $ \, \alpha^1,
\alpha^2, \dots, \alpha^N \, $  of  $ R^+ \, $:  following Lusztig we can
construct root vectors  $ E_{\alpha^r} $,
($ r=1, \dots, N \, $)  as in [DP] or [CV-1] and get PBW bases of  {\it
increasing\/}  ordered monomials \  $ \left\{\, L_\mu \cdot \prod_{r=1}^N
F_{\alpha^r}^{f_r} \, \Big\vert \, \mu \in M ; f_1, \dots, f_N \in \N
\,\right\} $  \  for  $ \, {\uphiMbm {M} } \, $  and \  $ \Big\{\, L_\mu
\cdot \prod_{r=1}^N E_{\alpha^r}^{e_r} \, \Big\vert \, \mu \in M; e_1, \dots,
e_N \in \N \,\Big\} $  \  for  $ \, \uphiMbp{M} \, $  or similar PBW bases
of  {\it decreasing\/}  ordered monomials; the same construction also provide
PBW bases for  $ U_{\scriptscriptstyle -} $,
$ \uzM{M} $,  and  $ U_{\scriptscriptstyle +} $.
                                                     \par
   Now, for every monomial  $ {\Cal E} $  in the  $ E_i $'s,  let
$ \, s({\Cal E}) := {\,1\, \over \,2\,} \, \varphi(wt({\Cal E})) $,
$ r({\Cal E}) := {\,1\, \over \,2\,} \, r(\varphi(wt({\Cal E}))) $,
$ \, \rbar({\Cal E}) := {\,1\, \over \,2\,} \, \rbar(\varphi(wt({\Cal E})))
\, $,  where  $ \, wt({\Cal E}) \, $  denotes the weight of  $ \, {\Cal E}
\, $  ($ E_i $  having weight  $ \alpha_i \, $),  and similarly for every
monomial  $ {\Cal F} $  in the  $ F_i $'s,  ($ F_i $  having weight
$ -\alpha_i \, $).  Then
\hfill\break
  $$  \eqalign {
   \pi_\varphi \Bigg( \prod_{r=N}^1 F_{\alpha^r}^{f_r} \,\cdot\,  &  L_\mu ,
\prod_{r=N}^1 E_{\alpha^r}^{e_r} \cdot L_\nu \Bigg) =  \cr
   &  = q^{-\left( r(\mu) - r \left( \prod_{r=N}^1 F_{\alpha^r}^{f_r} \right)
\big\vert \nu - s \left( \prod_{r=N}^1 E_{\alpha^r}^{e_r} \right) \right)}
\prod_{r=1}^N \delta_{e_r,f_r} {{[e_r]}_{q_{\alpha^r}}\! ! \, \,
q_{\alpha^r}^{+{e_r \choose 2}} \over {\left( q_{\alpha^r}^{-1} -
q_{\alpha^r} \right)}^{e_r}}  \cr
   \overline{\pi_\varphi} \Bigg( L_\mu \cdot\,  &  \prod_{r=N}^1
E_{\alpha^r}^{e_r} , L_\nu \cdot
\prod_{r=N}^1 F_{\alpha^r}^{f_r} \Bigg) =  \cr
   &  = q^{\left( \rbar(\mu) - \rbar \left( \prod_{r=N}^1 E_{\alpha^r}^{f_r}
\right) \big\vert \nu - s \left( \prod_{r=N}^1 F_{\alpha^r}^{e_r} \right)\right)} \prod_{r=1}^N \delta_{e_r,f_r} {{[e_r]}_{q_{\alpha^r}}\! ! \, \,
q_{\alpha^r}^{-{e_r \choose 2}} \over {\left( q_{\alpha^r} -
q_{\alpha^r}^{-1} \right)}^{e_r}}  \cr }   \eqno (2.2)  $$
gives the values of DRT pairings on PBW monomials (cf.~[CV-1], Lemma 3.5, and
[CV-2], \S 1, up to normalizations).  Now define  {\sl modified root vectors}
$ \, F^\varphi_\alpha := L_{\tau_\alpha} F_\alpha =
F_\alpha L_{\tau_\alpha} \, $,  $ \, E^\varphi_\alpha :=
L_{\tau_\alpha} E_\alpha = E_\alpha L_{\tau_\alpha} \, $  for all  $ \,
\alpha
\in R^+ \, $  (and set  $ \, F^\varphi_i := F^\varphi_{\alpha_i}
\, $,  $ \, F^\varphi_i := F^\varphi_{\alpha_i} \, $).  Then
  $$  \eqalign {
   \pi_\varphi  &  \left( \prod_{r=N}^1 {\left( F^\varphi_{\alpha^r}
\right)}^{f_r}  \,\cdot\, L_{(1+\varphi) (\mu)} , \prod_{r=N}^1
E_{\alpha^r}^{e_r} \cdot L_\nu \right) =  \cr
   &  = q^{-\left( \mu \vert \nu - s \left( \prod_{r=N}^1 E_{\alpha^r}^{e_r}
\right) \right) - \sum_{h<k} (f_h \tau_{\alpha^h} \vert f_k \alpha^k)} \cdot
\prod_{r=1}^N \delta_{e_r,f_r} {{[e_r]}_{q_{\alpha^r}}\! ! \, \,
q_{\alpha^r}^{+{e_r \choose 2}} \over {\left( q_{\alpha^r}^{-1} -
q_{\alpha^r} \right)}^{e_r}}  \cr
   \overline{\pi_\varphi}  &  \left( L_{(1-\varphi)(\mu)} \cdot \prod_{r=N}^1
{\left( E^\varphi_{\alpha^r} \right)}^{e_r} , L_\nu \cdot \prod_{r=N}^1
F_{\alpha^r}^{f_r} \right) =  \cr
   &  = q^{+ \left( \mu \vert \nu - s \left( \prod_{r=N}^1 F_{\alpha^r}^{e_r}
\right) \right) + \sum_{h<k} (e_h \tau_{\alpha^h} \vert e_k \alpha^k)}
\prod_{r=1}^N \delta_{e_r,f_r} {{[e_r]}_{q_{\alpha^r}}\! ! \, \,
q_{\alpha^r}^{-{e_r \choose 2}} \over {\left( q_{\alpha^r} -
q_{\alpha^r}^{-1} \right)}^{e_r}}  \cr }   \eqno (2.3)  $$
(cf.~[C-V1], Lemma 3.5, and [C-V2], Proposition 1.9).  In the sequel
$ \uphim $,  resp.~$ \uphip $,  will be the  $ \kq $--subalgebra  of
$ \uphiMbm{M} $,  resp.~$ \uphiMbp{M} $,  generated by the
$ F^\varphi_i $'s, resp.~by the  $ E^\varphi_i $'s  ($ i = \unon $);  these
too have PBW bases of ordered monomials of modified root vectors.

\vskip7pt

   {\bf 2.5 Integer forms.}   \  Let  $ \, X^{(m)} := X^m \big/
{{[m]}_{q_i}\! !} \, $  and  $ \, \left(Y; c \atop t \right) := \prod_{s=1}^t
{{q_i^{c-s+1} Y - 1} \over  {q_i^s - 1}} \, $  be the so-called "divided
powers";  let  $ \gerUphiMbm{M} $  be the  $ \kqqm $--subalgebra  of
$ \uphiMbm{M} $  generated by  $ \, \left\{\, F_i^{(m)}, \left( M_i; c \atop
t \right), M_i^{-1} \, \bigg\vert \, m,c,t \in \N; i=\unon \,\right\} \, $.
Then  $ \gerUphiMbm{M} $  is a Hopf subalgebra of  $ \uphiMbm{M} $,
(cf.~[CV-2]) having a PBW basis (as a
$ \kqqm $--module)  of increasing ordered monomials
  $$  \left\{\, \prod_{i=1}^n \left( M_i; 0 \atop t_i
\right) M_i^{-Ent({t_i/2})} \cdot \prod_{r=1}^N
F_{\alpha^r}^{(n_r)} \, \bigg\vert \, t_1, \dots, t_n, n_1, \dots, n_N \in
\N \,\right\}  $$
and a similar PBW basis of decreasing ordered monomials; in particular
$ \gerUphiMbm{M} $  is a  $ \kqqm $--form  of
$ \uphiMbm{M} $.  \! Similarly we define the Hopf subalgebra
$ \gerUphiMbp{M} $  and locate PBW bases for it.
                                               \par
  Let  $ \, \ebar_{\alpha^r} := \left( q_{\alpha^r} - q_{\alpha^r}^{-1}
\right) E_{\alpha^r} \, $,  $ \forall \, r=1, \dots, N \, $,  and let
$ \, \calUphiMbp{M} \, $  be the  $ \kqqm $--subalgebra  of  $ \uphiMbp{P} $
generated by  $ \, \{ \ebar_{\alpha^1}, \dots, \ebar_{\alpha^N} \} \cup
\{\, M_1^\pm, \dots , M_n^\pm \,\} \, $;  then (cf.~[DKP], [DP])
$ \calUphiMbp{M} $  is a Hopf subalgebra of  $ \uphiMbp{M} $,  having a PBW
basis (as a  $ \kqqm $--module)
  $$  \left\{\, \prod_{i=1}^n M_i^{t_i} \cdot \prod_{r=1}^N
\ebar_{\alpha^r}^{n_r} \, \bigg\vert \, t_1, \dots, t_n \in \Z ; n_1, \dots,
n_N \in \N  \,\right\}  $$
of increasing ordered monomials and a similar PBW basis of decreasing
ordered monomials; in particular
$ \calUphiMbp{M} $  is a  $ \kqqm $--form  of  $ \uphiMbp{M} \, $.  The same
procedure yields the definition of the Hopf subalgebra  $ \calUphiMbm{M} \, $
and provides PBW bases for it.  The same integer forms can also be
constructed using modified root vectors instead of the usual ones, hence
these integer forms have also {\sl modified}  PBW bases of ordered
monomials in the  $ M_i $'s  and the modified root vectors.
Similar constructions and results hold for the algebras
$ \uphim $,  $ \uphizM{M} $,  $ \uphip $,  providing integer
forms  $ \gerUphim $,  $ \calUphip $,  and so on.  Finally, we have
decompositions
  $$  \displaylines{
   {\gerUphiMbm M } \cong \gerUphim \otimes {\gerUphizM M } \cong
{\gerUphizM M } \otimes \gerUphim \, ,  \quad  {\gerUphiMbp M } \cong \gerUphip
\otimes {\gerUphizM M } \cong {\gerUphizM M } \otimes \gerUphip  \cr
   {\calUphiMbm M } \cong \calUphim \otimes {\calUphizM M } \cong
{\calUphizM M } \otimes \calUphim \, ,  \quad  {\calUphiMbp M } \cong \calUphip
\otimes {\calUphizM M } \cong {\calUphizM M } \otimes \calUphip  \cr }  $$

\vskip7pt

   {\bf 2.6  $ \kqqm $--duality among integer forms.}   \  The very
definitions and (2.3) imply that integer forms of opposite "fonts"
(namely  $ \gerU $  or $ \calU \, $)  are  {\it  $ \kqqm $--dual\/}
of each other in the following sense: for every DRT pairing, if we take
$ \gerU $  on one side, then the form  $ \calU $  on the other side
coincides with the subset of all elements which paired with
$ \gerU $  give a value in  $ \kqqm $;  and similarly reverting
the roles of  $ \gerU $  and  $ \calU $.  For instance
  $$  \displaylines{
   {\gerUphizM M } = \Big\{\, y \in {\uphizM M } \,\Big\vert\, \pi_\varphi \Big(
{\calUphizM {M'} },y \Big) \subseteq \kqqm \,\Big\} = \Big\{\, x \in {\uphizM
M } \,\Big\vert\, \overline{\pi_\varphi} \Big( x, {\calUphizM {M'} } \Big)
\subseteq \kqqm \,\Big\}  \cr
   {\calUphizM M } = \Big\{\, y \in {\uphizM M } \,\Big\vert\, \pi_\varphi \Big(
{\gerUphizM {M'} }, y \Big) \subseteq \kqqm \,\Big\} = \Big\{\, x \in {\uphizM
M } \,\Big\vert\, \overline{\pi_\varphi} \Big( x, {\gerUphizM {M'} } \Big)
\subseteq \kqqm \,\Big\}  \cr
   \gerUphim = \Big\{\, x \in \uphim \,\Big\vert\, \pi_\varphi \Big( x,
\calUphip \Big) \subseteq \kqqm \,\Big\} = \Big\{\, y \in \uphim \,\Big\vert\,
\overline{\pi_\varphi} \Big( \calUphip, y \Big) \subseteq \kqqm \,\Big\}  \cr
   {\calUphiMbp M } = \Big\{\, x \in {\uphiMbp M } \,\Big\vert\,
\overline{\pi_\varphi} \Big( x, {\gerUphiMbm {M'} } \Big) \subseteq \kqqm
\,\Big\} = \Big\{\, y \in {\uphiMbp M } \,\Big\vert\, \pi_\varphi \Big(
{\gerUphiMbm {M'} }, y \Big) \subseteq \kqqm \,\Big\}  \cr }  $$

\vskip1,7truecm

\centerline{ \bf  \S \; 3 \,  The quantum group
$ \uqphiMg{M} $ }

\vskip10pt

   {\bf 3.1  The quantum double.}  \  Let  $ H_- $,  $ H_+ $
be two arbitrary Hopf algebras on a ground field (or ring)
$ F $,  and let  $ \, \pi \colon {\big( H_- \big)}_{op} \otimes H_+
\rightarrow F \, $  be any arbitrary Hopf pairing.  Then Drinfel'd's  {\it
quantum double \/}  $ \, D = D \big( H_-, H_+, \pi \big) \, $  is the algebra
$ \, T \big( H_- \oplus H_+ \big) \big/ {\Cal R} \, $,  where  $ {\Cal R} $
is the ideal of relations
  $$  \eqalign{
    1_{H_-} = 1 = 1_{H_+}  \,  ,  &  \qquad  x \otimes y
= x y  \qquad \qquad \;  \hbox{for} \; \, x, y \in H_+  \; \hbox{ or }
\; x, y \in H_- \, \phantom{.}  \cr
    \sum_{(x),(y)} \pi \left( y_{(2)}, x_{(2)} \right) \, x_{(1)} \otimes
y_{(1)} =  &  \sum_{(x),(y)} \pi \left( y_{(1)},x_{(1)} \right) \, y_{(2)}
\otimes x_{(2)}  \quad \;  \hbox{for}  \; \, x \in H_+ , \,  y \in H_- \, .
\cr }  $$
   \indent   Then (cf.~[DL], Theorem 3.6)  $ D $  has a
canonical structure of Hopf algebra such that  $ \, H_- $,
$ H_+ \, $  are Hopf subalgebras of it and multiplication
yields isomorphisms of coalgebras
  $$  H_+ \otimes H_- \longhookrightarrow D \otimes D {\buildrel {m} \over
\llongrightarrow} D \, ,  \qquad  H_- \otimes H_+ \longhookrightarrow D
\otimes D {\buildrel {m} \over \llongrightarrow} D \, .   \eqno (3.1)  $$
   \indent   Now consider  $ \, D_{q,\varphi}^{\scriptscriptstyle M} (\gerg) :=
D \left( {\uphiMbm Q }, {\uphiMbp M }, \pi_\varphi \right) \, $;  by definition,
 $ D_{q,\varphi}^{\scriptscriptstyle M} (\gerg) $  is generated by
$ \, K_\alpha $,  $ L_\mu $,  $ F_i $,  $ E_i $   --- identified with
$ 1 \otimes K_\alpha $,  $ L_\mu \otimes 1 $,  $ 1 \otimes F_i $,
$ E_i \otimes 1 $  via  $ \, D_{q,\varphi}^{\scriptscriptstyle M} (\gerg) \cong
\uphiMbp{M} \otimes \uphiMbm{Q} \, $ ---   ($ \alpha \in Q $,  $ \mu \in M $,
$ i = \unon $),  wile the relations defining  $ {\Cal R} $  reduce to
  $$  \eqalign{
      K_\alpha L_\mu = L_\mu K_\alpha  \; ,  \qquad  K_\alpha E_j =  &
{} \, q_i^{+(\alpha_j \vert \alpha)} E_j K_\alpha  \; ,  \qquad  L_\mu F_j =
q_i^{-(\alpha_j \vert \mu)} F_j L_\mu  \cr
      E_i F_j - F_j E_i =  &  {} \, \delta_{ij} {{L_{\alpha_i} -
K_{-\alpha_i}} \over {q_i - q_i^{-1}}}  \cr  }   \eqno (3.2)  $$
   \indent   Finally, PBW bases of quantum Borel algebras provide PBW bases
of  $ D_{q,\varphi}^{\scriptscriptstyle M} (\gerg) $.  In the sequel we shall
also use the notation  $ \, D_{\scriptscriptstyle M} :=
D_{q,\varphi}^{\scriptscriptstyle M} (\gerg) \, $.

\vskip7pt

   {\bf 3.3 The quantum algebra  $ \, \uqphiMg{M} \, $.}  \  Let
$ \gerK^{\scriptscriptstyle M}_\varphi $  be the ideal of
$ D_{q,\varphi}^{\scriptscriptstyle M} (\gerg) $  generated by the elements
$ \, L \otimes 1 - 1 \otimes L \, $,  $ \, L \in \uphizM{M} \, $;
$ \gerK^{\scriptscriptstyle M}_\varphi $  is in fact a Hopf ideal, whence
$ \, D_{q,\varphi}^{\scriptscriptstyle M} (\gerg) \big/
\gerK^{\scriptscriptstyle M}_\varphi \, $  is a Hopf algebra.  Then from
above  we get a presentation of  $ \, \uqphiMg{M} :=
D_{q,\varphi}^{\scriptscriptstyle M} (\gerg) \Big/ \gerK^{\scriptscriptstyle
M}_\varphi \, $:
            it is\break
 \eject
\noindent   the associative  $ \kq $--algebra  with 1 given by
generators  $ \, F_i $,  $ L_\mu $,  $ E_i \, $  and relations
 $$  \eqalign {
      \!\! L_0 = 1 \, ,  \quad  L_\mu L_\nu = L_{\mu + \nu} = L_\nu L_\mu \, ,
\quad  L_\mu F_i  &  = q^{-(\alpha_j | \mu)} F_i L_\mu \, , \quad  L_\mu E_i
= q^{+(\alpha_j | \mu)} E_i L_\mu  \cr
      E_i F_h - F_h E_i =  &  {} \,\, \delta_{ih} {{L_{\alpha_i} -
L_{-\alpha_i}}
\over {q_i - q_i^{-1}}}  \cr
      \!\! \sum_{k = 0}^{1-a_{ij}} \! (-1)^k \! \left[ \! { 1-a_{ij} \atop
k } \! \right]_{\!q_i} \! E_i^{1-\aij-k} E_j E_i^k = 0 , \;  &
\sum_{k=0}^{1-a_{ij}} \! (-1)^k \! \left[ \! { 1-a_{ij} \atop k } \!
\right]_{\!q_i} \! F_i^{1-\aij-k} F_j F_i^k = 0  \cr }   \eqno (3.3)  $$
(for all  $ \, \mu \in M $,  $ i, j, h = \unon \, $ with  $ \, i \neq j
\, $)  with the Hopf structure given by
  $$  \matrix
       \Delta_\varphi (F_i) = F_i \otimes L_{-\alpha_i -
\tau_i} + L_{\tau_i} \otimes F_i \, ,  &  \;\;  \epsilon_\varphi(F_i) = 0  \, ,
&  \;\;  S_\varphi (F_i) = - F_i L_{\alpha_i}  \\
       \Delta_\varphi (L_\mu) = L_\mu \otimes L_\mu \, ,  &  \;\;
\epsilon_\varphi (L_\mu) = 1 \, ,  &  \;\;  S_\varphi (L_\mu) = L_{-\mu}  \\
       \Delta_\varphi (E_i) = E_i \otimes L_{\tau_i} +
L_{\alpha_i - \tau_i} \otimes F_i \, ,  &  \;\;  \epsilon_\varphi (E_i) = 0 \, ,
 &  \;\;  S_\varphi (F_i) = - L_{-\alpha_i} E_i    \\
      \endmatrix   \eqno (3.4)  $$
   \indent   For  $ \, \varphi = 0 \, $  one recovers the usual one-parameter
quantum enveloping algebras.  Finally we let  $ \, pr_{\scriptscriptstyle M}
\colon \, D_{q,\varphi}^{\scriptscriptstyle M} (\gerg)
\llongtwoheadrightarrow D_{q,\varphi}^{\scriptscriptstyle M} (\gerg) \Big/
\gerK^{\scriptscriptstyle M}_\varphi =: \uqphiMg{M} \, $  be the canonical
Hopf algebra epimorphism; we shall also use notation  $ \, K_\alpha :=
L_\alpha $,  $ M_\mu := L_\mu $,  $ \; \forall \, \alpha \in Q, \mu \in
M \, $.

\vskip7pt

   {\bf 3.4 Integer forms of  $ \uqphiMg{M} $.} \  Let  $ \gerUphiMg{M} $  be
the  $ \kqqm $--subalgebra  of $ {\uqphiMg M } $  generated by  $ \,
\left\{\, F_i^{(\ell)}, \left( M_i; c \atop t \right), M_i^{-1}, E_i^{(m)} \,
\bigg\vert \, \ell,c,t,m \in \N; i=\unon \,\right\} \, $;  this is a Hopf
subalgebra of  $ {\uqphiMg M } $  (cf.~[DL], \S 3), with PBW basis (over
$ \kqqm $)
  $$  \left\{\, \prod_{r=N}^1 E_{\alpha^r}^{(n_r)} \cdot \prod_{i=1}^n
\left( M_i; 0 \atop t_i \right) M_i^{-Ent({t_i/2})} \cdot \prod_{r=1}^N
F_{\alpha^r}^{(m_r)} \, \bigg\vert \, n_r, t_i, m_r \in \N, \forall\, r, i
\,\right\} \; ;  $$
this is also a  $ \kq $--basis  of  $ {\uqphiMg M } $,  hence  $ {\gerUMg M } $
is a  $ \kqqm $--form  of $ {\uqphiMg M } \, $.
                                                 \par
   Let  $ \, \calUphiMg{M} \, $  be the  $ \kqqm $--subalgebra of
$ \uqphiMg{M} $  generated by (cf.~[DP], \S 12)
  $$  \left\{ \fbar_{\alpha^1}, \dots, \fbar_{\alpha^N} \right\} \cup
\left\{ M_1^{\pm 1}, \dots, M_n^{\pm 1} \right\} \cup \left\{ \ebar_{\alpha^1},
\dots, \ebar_{\alpha^N} \right\} \; ;  $$
this is a Hopf subalgebra of  $ \uqphiMg{P} $,  having a PBW basis (over
$ \kqqm $)
  $$  \left\{\, \prod_{r=N}^1 \ebar_{\alpha^r}^{\,n_r} \cdot
\prod_{i=1}^n M_i^{t_i} \cdot \prod_{r=1}^N \fbar_{\alpha^r}^{\,m_r} \,
\bigg\vert \, t_i \in \Z, n_r, m_r \in \N, \forall\, i, r \,\right\} ;  $$
the latter is also a  $ \kq $--basis  of  $ \uqphiMg{M} $,  hence
$ \calUphiMg{M} $  is a  $ \kqqm $--form  of
$ \uqphiMg{M} $.
                                                 \par
   Like for quantum Borel algebras, the same forms can also be defined using
modified root vectors, hence they have also PBW bases of ordered monomials in
the  $ M_i $'s  and the modified root vectors.

\vskip7pt

   {\bf 3.5  Specialization at roots of 1 and quantum Frobenius
morphisms.}  \  When dealing with specializations, if any scalar  $ \, c \in k
\setminus \{0\} \, $  is fixed then  $ k $  is thought of as a
$ \kqqm $--algebra  via  $ \, k \cong \kqqm \big/ (q-c) \, $.
                                               \par
   Let  $ \varepsilon $  be a primitive  $ \ell $--th
root of 1, for  $ \ell $  {\it odd\/},  $ \, \ell > d:=
\max_i {\{d_i\}}_i \, $,  or  $ \, \ell = 1 \, $.  Then we set
$ \, \gerUepsilonphiMg{M} := \, \gerUphiMg{M} \Big/ (q -
\varepsilon) \, \gerUphiMg{M} \, \cong \, \gerUphiMg{M} \otimes_{k \left[ q,
\qm \right]} k \, $.  When  $ \, \ell =1 \, $
(i.~e.~$ \, \varepsilon = 1 \, $)  it is well-known (cf.~e.~g.~[CV-2] or
[DL])\footnote{This result is more general than in  [{\it loc.~cit.\/}]:  it
can be proved on the same lines of Theorem 7.2 below.}
that  $ \, \gerUunophiMg{M} \, $  is a
                         Poisson Hopf coalgebra, and we\break
 \eject
\noindent   have a Poisson Hopf coalgebra isomorphism
  $$  \gerUunophiMg{M} \cong \ugtau \; ;   \eqno (3.5)  $$
in a word,  $ \gerUphiMg{M} $  specializes to  $ \ugtau $  for  $ \, q
\rightarrow 1 \, $:  in symbols,  $ \, \gerUphiMg{M} @>{q \rightarrow 1}>>
\ugtau \, $.
                                                \par
   When  $ \, \ell > 1 \, $,  from [CV-2], \S 3.2 (cf.~also [Lu], [DL]) we
have an epimorphism
  $$  \gerFrgtau \colon \, \gerUepsilonphiMg{M}
\llongtwoheadrightarrow \gerUunophiMg{M} \cong \ugtau   \eqno (3.6)  $$
of Hopf algebras defined by (recall that  $ \, M_i := L_{\mu_i} \, $)
  $$  {} \!\!\!\!\!  \gerFrgtau \colon
 \cases
   \! F_i^{(s)} \Big\vert_{q=\varepsilon} \!\!\!\! \mapsto F_i^{(s / \ell)}
\Big\vert_{q=1},  \left( M_i; 0 \atop s \right) \! \bigg\vert_{q=\varepsilon}
\!\!\!\! \mapsto \left( M_i; 0 \atop s / \ell \right) \! \bigg\vert_{q=1},
E_i^{(s)} \Big\vert_{q=\varepsilon} \!\!\!\! \mapsto E_i^{(s / \ell)}
\Big\vert_{q=1}  \, \hbox{ if } \; \ell \Big\vert s  \\
   \! F_i^{(s)} \Big\vert_{q=\varepsilon} \!\!\! \mapsto 0,  \quad  \left(
M_i;
0 \atop s \right) \! \bigg\vert_{q=\varepsilon} \!\!\! \mapsto 0,  \quad
E_i^{(s)} \Big\vert_{q=\varepsilon} \!\!\! \mapsto 0  \; \hbox{ \, otherwise
\, }  \\
   \! M_i^{-1} \Big\vert_{q=1} \!\!\! \mapsto 1  \\
 \endcases
   \eqno (3.7)  $$
   \indent   If  $ \, \varphi = 0 \, $   --- whence  $ \, \tau = 0 \, $ ---
and  $ \, \ell = p \, $  is prime, it is shown
in [Lu], \S 8.15, that  $ {{\frak F}r}_{\gerg^0} $  (for  $ \, M = Q \, $)
can be regarded as a  {\it lifting of the Frobenius morphism  $ \, G_{\Z_p}
\rightarrow G_{\Z_p} \, $  to characteristic zero\/};  for this reason, we
refer to  $ \gerFrgtau $  as a  {\it quantum Frobenius morphism\/}.
                                                  \par
   Similarly, we set  $ \, \calUepsilonphiMg{M} := \, \calUphiMg{M} \Big/
(q - \varepsilon) \, \calUphiMg{M} \, \cong \, \calUphiMg{M}
\otimes_{k \left[ q, \qm \right]} k \, $;  when  $ \, \ell = 1 \, $  it is
known (cf.~[DP], Theorem 12.1, and [DKP], Remark 7.7 {\it (c)\/}\,) that
  $$  \calUunophiMg{M} \cong F \left[ H^\tau_{\scriptscriptstyle M} \right]
\eqno (3.8)  $$
as Poisson Hopf algebras over  $ k \, $:  here
$ H^\tau_{\scriptscriptstyle M} $  is the connected
Poisson group with tangent Lie bialgebra  $ \gerh^\tau $   --- defined in
\S 1.2 ---   and  $ M $  the character group of a maximal torus.  In a word,
$ \calUphiMg{M} $  specializes to  $ F \left[ H^\tau_{\scriptscriptstyle M}
\right] $  as  $ \, q \rightarrow 1 \, $,  or  $ \, \calUphiMg{M} @>{q
\rightarrow 1}>> F \left[ H^\tau_{\scriptscriptstyle M} \right] \, $.  When
$ \, \ell > 1 \, $,  from [DKP], \S\S 7.6--7 we record the existence of a
Hopf
algebra monomorphism
  $$  \calFrgtau \colon \, F \left[ H^\tau_{\scriptscriptstyle M} \right]
\cong \calUunophiMg{M}
\llonghookrightarrow \calUepsilonphiMg{M}   \eqno (3.9)  $$
(cf.~also [DP] for the one-parameter case) defined by  ($ \, \alpha \in
R^+ \, $,  $ \, \mu \in M \, $)
  $$  \calFrgtau \colon \,\, \;\; \fbar_\alpha \Big\vert_{q=1} \mapsto
{\fbar_\alpha}^{\! \ell} \Big\vert_{q=\varepsilon} \, ,  \;\; L_\mu
\Big\vert_{q=1} \mapsto {L_\mu}^{\! \ell} \Big\vert_{q=\varepsilon} \, ,
\;\; \ebar_\alpha \Big\vert_{q=1} \mapsto {\ebar_\alpha}^{\! \ell}
\Big\vert_{q=\varepsilon}   \eqno (3.10)  $$
   \indent   Again, we refer to
$ \calFrgtau $  as a {\it quantum Frobenius morphism\/}:  if
$ \, \varphi = 0 \, $  and  $ \, \ell = p \, $  is prime it is a
{\it lifting of the Frobenius morphism  $ \, H_{\Z_p} \rightarrow H_{\Z_p}
\, $  to characteristic        zero\footnote{ Here  $ H_{\Z_p} $  denotes
the Chevalley-type group-scheme over  $ \Z_p $  associated to  $ H $ }}.
 \eject

\vskip1,7truecm

 \centerline{ \bf  \S \; 4 \,  Quantum function algebras }

\vskip10pt

   {\bf  4.1  The quantum function algebras  $ \fqphiMbs{M}{\pm} $.}  \  Let
$ \fqphiMbs{M}{\pm} $  be the  {\it quantum function algebra\/}  relative to
$ U_{q,\varphi}^{\scriptscriptstyle M'} (\gerb_\pm) $,  defined as the
algebra of matrix coefficients of
    {\sl positive}\footnote{Namely those having a basis on which the
                            $ L_\nu $'s  ($ \nu \in M' $)  act diagonally
                            by powers of  $ q $.}
finite dimensional representations of
$ U_{q,\varphi}^{\scriptscriptstyle M'} (\gerb_\pm) $.  Then
$ \fqphiMbs{M}{\pm} $  is a Hopf algebra, which we call dual of
$ U_{q,\varphi}^{\scriptscriptstyle M'} (\gerb_\pm) $  for there is
a perfect Hopf pairing (evaluation) among them; in fact
$ \fqphiMbs{M}{\pm} $  is a Hopf subalgebra of
$ {U_{q,\varphi}^{\scriptscriptstyle M'} (\gerb_\pm)}^\circ $  (the   ---
restricted ---   Hopf dual of  $ U_{q,\varphi}^{\scriptscriptstyle M'}
(\gerb_\pm) $,  in the sense of [SW], ch.~VI).  The DRT pairings provide Hopf
algebra isomorphisms
  $$  \fqphiMbp{M} \cong {\left( \uphiMbm{M'} \right)}_{op} \, ,  \quad \;
\fqphiMbm{M} \cong {\left( \uphiMbp{M'} \right)}_{op}   \eqno (4.1)  $$
induced by the pairing  $ \, \pi_\varphi \, $,
resp.~$ \, \overline{\pi_\varphi} \, $;  by means  of these, the DRT pairings
can be seen as natural evaluation pairings (cf.~[DL], \S 4, and [CV-2], \S\S
2--3).

\vskip7pt

   {\bf 4.2  Integer forms of  $ \, \fqphiMbs{M}{\pm} \, $.}  \  Let
  $$  \eqalign{
      \gerFphiMbs{M}{\pm} := \left\{\, f \in
\fqphiMbs{M}{\pm} \,\big\vert\, \left\langle f,
\gerU_\varphi^{\scriptscriptstyle M'} (\gerb_\pm) \right\rangle \subseteq
\kqqm \,\right\}  \cr
      \calFphiMbs{M}{\pm} := \left\{\, f \in
\fqphiMbs{M}{\pm} \,\big\vert\, \left\langle f,
{\Cal U}_\varphi^{\scriptscriptstyle M'} (\gerb_\pm) \right\rangle \subseteq
\kqqm \,\right\}  \cr }   \eqno (4.2)  $$
where \  $ \langle \  , \  \rangle \colon \fqphiMbs{M}{\pm} \otimes
U_{q,\varphi}^{\scriptscriptstyle M'} (\gerb_\pm) \rightarrow \kq $  \  is
the natural evaluation pairing; then
  $$  \gerFphiMbs{M}{\pm} \cong {\left( {\Cal U}_\varphi^{\scriptscriptstyle
M} (\gerb_\mp) \right)}_{op} \, ,  \quad  \calFphiMbs{M}{\pm} \cong {\left(
\gerU_\varphi^{\scriptscriptstyle M} (\gerb_\mp) \right)}_{op} \, ;
\eqno (4.3)  $$
because of \S 2.6 and (4.1): in particular  $ \gerFphiMbs{M}{\pm} $  and
$ \calFphiMbs{M}{\pm} $  are integer forms of  $ \fqphiMbs{M}{\pm} $.

\vskip7pt

   {\bf 4.3  The quantum function algebra  $ \fqphiMg{M} $  and its integer
forms.}  \  Like in \S 4.1, we define  the  {\it quantum
function algebra\/}  $ \fqphiMg{M} $  (relative to  $ \uqphiMg{M'} $)  to be the
algebra of matrix coefficients of  {\it positive\/}  finite dimensional
representations of  $ \uqphiMg{M'} $  (cf.~[DL], \S 4, and [CV-2], \S 2.1); it
is a Hopf subalgebra of  $ {\uqphiMg {M'} }^\circ $, perfectly paired with
$ \uqphiMg{M'} $  by the natural evaluation pairing (whence we call it
{\sl dual\/}  of  $ \uqphiMg{M'} \, $).  As for integer forms, let
  $$  \eqalign{
      \gerFphiMg{M} := \left\{\, f \in \fqphiMg{M} \,\big\vert\,
\left\langle f, \gerUphiMg{M'} \right\rangle \subseteq \kqqm \,\right\}  \cr
      \calFphiMg{M} := \left\{\, f \in \fqphiMg{M} \,\big\vert\,
\left\langle f, \calUphiMg{M'} \right\rangle \subseteq  \kqqm \,\right\}
\cr }   \eqno (4.4)  $$
where \  $ \langle \  , \  \rangle \colon \fqphiMg{M} \otimes \uqphiMg{M'}
\rightarrow \kq $  \  is the natural evaluation pairing; we shall later
prove that these are  $ \kqqm $--integer  forms (as Hopf subalgebras) of
$ \fqphiMg{M} $.

\vskip7pt

   {\bf 4.4  Specialization at roots of 1.}   \   Let
$ \varepsilon $  be a primitive  $ \ell $--th  root of 1 in  $ k $  (with the
assumptions of \S 3.5 on  $ \ell \, $),  and set  $ \, \gerFepsilonphiMg{M} :=
\, \gerFphiMg{M} \Big/ \, (q - \varepsilon) \, \gerFphiMg{M} \, \cong \,
\gerFphiMg{M} \otimes_{k \left[ q, \qm \right]} k \, $.  For  $ \, \ell
= 1 \, $,  we have  $ \, \gerFunophiMg{M} \, \cong \, F \left[
G^\tau_{\scriptscriptstyle M} \right] \, $  as Poisson Hopf  $ k $--algebras
(cf.~[CV-2], [DL]), i.~e.
  $$  \gerFphiMg{M} \, {\buildrel {q \rightarrow 1} \over \llongrightarrow}
\, F \left[ G^\tau_{\scriptscriptstyle M} \right] \; ;  $$
here  $ G^\tau_{\scriptscriptstyle M} $  is the connected Poisson group with
tangent Lie bialgebra  $ \gerg^\tau $  and  $ M $  as character group of a
maximal torus.  In fact this result arises as dual of  $ \, \gerUphiMg{M'} \,
@>{q \rightarrow 1}>> \, \ugtau \, $.  When  $ \, \ell > 1 \, $,  another
{\it quantum Frobenius morphism\/},  namely a Hopf algebra monomorphism
  $$  \gerFrGtau \colon \, F \left[ G^\tau_{\scriptscriptstyle M} \right]
\cong \gerFunophiMg{M} \llonghookrightarrow \gerFepsilonphiMg{M} \; ,
\eqno (4.5)  $$
is defined (cf.~[CV-2], \S 3.3), which is dual of  $ \, \gerFrgtau \colon
\, \gerUepsilonphiMg{M'} \longtwoheadrightarrow \gerUunophiMg{M'} \cong \ugtau
\, $.
 \eject

\vskip1,7truecm

 \centerline{ \bf  \S \; 5 \,  Quantum formal groups }

\vskip10pt

   {\bf 5.1 Formal Hopf algebras and quantum formal groups.}
\  In this subsection we introduce the notion of  {\it
quantum formal group\/}.  Recall (cf.~[Di], ch.~I) that
formal groups can be defined in a category of a special type of commutative
topological algebras, whose underlying vector space (or module) is linearly
compact; following Drinfel'd's philosophy,we
define quantum formal groups by simply dropping out any commutativity
assumption of the classical notion of formal group; thus now we quickly
outline how to modify the latter (following [Di], ch.~I) in
order to define our new quantum objects.
                                                 \par
   Let  $ E $  be any vector space over
a field  $ K $  (one can then generalize more or less
wathever follows to the case of free modules over a ring),
and let  $ E^* $  be its (linear) dual; we write  $ \, \langle x^*, x
\rangle \, $  for  $ x^*(x) $  for  $ \, x \in E \, $,  $ \, x^* \in
E^* \, $.  We consider on  $ E^* $  the  {\it weak  $ \ast $--topology},
i.~e.~the coarsest topology such
that for each  $ \, x \in E \, $  the linear map
$ \, x^* \mapsto \langle x^*, x \rangle \, $
of  $ E^* $  into  $ K $  is continuous, when  $ K $  is given the
discrete topology.  We can describe this
topology by choosing a basis  $ {\{e_i\}}_{i \in I} $  of  $ E \, $:  to
each  $ i \in I $  we associate the linear (coordinate) form  $ e_i^* $  on
$ E $  such that  $ \, \langle e_i^*, e_j \rangle = \delta_{i j} \, $,  and
we say that the family  $ {\{e^*_i\}}_{i \in I} $  is the
{\it pseudobasis\/}  of  $ E^* $  dual to  $ {\{e_i\}}_{i \in I} $;  then the
subspace  $ E' $  of  $ E $  which is (algebraically)
generated by the  $ e_i^* $  is dense in  $ E^* $,  and  $ E^* $ is nothing
but the  {\it completion\/}  of  $ E' $,  when  $ E' $  is given the topology
for which a fundamental system of
neighborhoods of  $ 0 $  consists of the vector subspaces containing almost
all the  $ e_i^* $;  thus
elements of  $ E^* $  can be described by series in the  $ e_i^* $'s  which
in the given topology are
in fact convergent.  Finally, the topological vector spaces  $ E^* $  are
characterized by the
property of linear compactness.
                                                 \par
   Let now  $ E $,  $ F $  be any two vector spaces over
$ K $,  and  $ \, u \colon E \rightarrow F \, $  a linear map; then the
dual map  $ \, u^* \colon F^* \rightarrow E^* \, $ is  continuous, and
conversely for any linear map  $ \, v \colon F^* \rightarrow E^* \, $  which
is continuous there exists a unique linear map  $ \, u \colon E
\rightarrow F \, $  such that  $ \, v = u^* \, $.
                                              \par
   The tensor product  $ E^* \otimes F^* \, $
is naturally identified to a subspace of  $ {(E \otimes
F)}^* $  by  $ \, \langle x^* \otimes y^*, x \otimes y \rangle = \langle x^*,
x \rangle \cdot \langle y^*, y \rangle \, $;  thus if
$ {\{e_i\}}_{i \in I} $  and  $ {\{f_j\}}_{j \in J} $  are bases of  $ E $
and  $ F $,  and  $ {\{e^*_i\}}_{i \in I} $  and
$ {\{f^*_j\}}_{j \in J} $  their dual pseudobases in
$ E^* $  and  $ F^* $,  then  $ {\{e^*_i \otimes f^*_j\}}_{i \in I, j
\in J} $  is the dual pseudobasis of  $ {\{e_i \otimes f_j\}}_{i \in I, j
\in J} $  in  $ {(E \otimes F)}^* $.  Thus  $ {(E \otimes F)}^* $  is the
completion of  $ E^* \otimes F^* $  for the tensor product topology,
i.~e.~the topology of  $ E^* \otimes F^* $  for which a fundamental system
of neighborhoods of  $ 0 $  consists of the sets  $ \, E^* \otimes V +
W \otimes F^* \, $  where  $ V $,  resp.  $ W $,  ranges in a fundamental
system of neighborhoods of  $ 0 $  made of vector subspaces; we denote this
completion by  $ E^* \otimeshat F^* $,  and we call it the
{\it completed\/}  (or  {\it topological\/})  {\it tensor product\/}  of
$ E^* $  and  $ F^* $;  the embedding  $ \, E^* \otimes F^*
\longhookrightarrow {(E \otimes F)}^* = E^* \otimeshat F^* \, $  is then
continuous.  Finally, if  $ \, u \colon E_1 \rightarrow  E_2
\, $,  $ \, v \colon F_1 \rightarrow F_2 \, $  are linear maps, then  $ \,
{(u \otimes v)}^* \colon {(E_2 \otimes F_2)}^* = {E_2}^* \otimeshat {F_2}^*
\longrightarrow {(E_1 \otimes F_1)}^* = {E_1}^* \otimeshat {F_1}^* \, $
coincides with the continuous extension to  $ {E_2}^* \otimeshat {F_2}^* $
of the continuous map  $ \, u^* \otimes v^* \colon {E_2}^* \otimes {F_2}^*
\rightarrow {E_1}^* \otimes  {F_1}^* \, $;  thus it is also denoted by  $ u^*
\otimeshat v^* \, $.
                                                  \par
   We define a  {\it linearly compact algebra}  to be a
topological algebra whose underlying vector space (or free
module) is linearly compact: then
linearly compact algebras form a full subcategory of the category of
topological algebras; morever, for any two objects  $ A_1 $  and
$ A_2 $  in this category, their topological tensor product  $ A_1
\otimeshat A_2 $  is defined.  Dually, within the category of linearly
compact vector spaces we define  {\it linearly compact
coalgebras}  as triplets  $ (C, \Delta, \epsilon) $  with  $ \, \Delta
\colon C \rightarrow C \otimeshat C \, $  and  $ \, \epsilon \colon C
\rightarrow K \, $  satisfying the usual coalgebra
axioms.  The arguments in [Di] (which never
require commutativity nor cocommutativity) show that
$ \, {(\ )}^* \colon (A, m, 1) \mapsto (A^*, m^*, 1^*) \, $  defines a
contravariant functor from algebras to linearly compact coalgebras, while
$ {(\ )}^* \colon (C, \Delta, \epsilon) \mapsto (C^*, \Delta^*, \epsilon^*)
\, $  defines a contravariant functor from coalgebras to linearly compact
algebras.  Finally, we define a  {\it formal Hopf algebra}  as a
datum  $ (H, m, 1, \Delta, \epsilon, S) $  such that $ (H, m, 1) $  is a
linearly compact algebra,  $ (H, \Delta, \epsilon) $  is a linearly compact
coalgebra, and the usual compatibility axioms of Hopf
algebras are satisfied.  "Usual" Hopf algebras are particular cases of formal
Hopf algebras.
                                                 \par
   We define  {\bf quantum formal group}  the  {\it spectrum\/}  of a formal
Hopf algebra (whereas  {\it classical\/}  formal groups are
spectra of  {\it commutative\/}  formal Hopf algebras:
cf.~[Di], ch.~I).
                                           \par
   Our goal is to study  $ \, {\uqphiMg M }^* \, $.  Since
$ {\uqphiMg M } $  is a Hopf algebra, its linear dual  $ {\uqphiMg M }^* $
is a  {\sl formal}  Hopf algebra.  The functor  $ {(\ \;)}^* $  turns the
natural epimorphism  $ \, {pr}_{\! \scriptscriptstyle M} \colon
D_{\scriptscriptstyle M} \longtwoheadrightarrow \uqphiMg{M} \, $ into a
monomorphism  $ \, j_{\scriptscriptstyle M'} \! :=
{\left( {pr}_{\! \scriptscriptstyle M} \right)}^* \colon {\uqphiMg M }^*
\longhookrightarrow {D_{\scriptscriptstyle M}}^* \, $  of formal Hopf
algebras: therefore we begin by studying  $ {D_{\scriptscriptstyle M}}^* $.  The
following is straightforward:

\vskip7pt

\proclaim{Proposition 5.2}  Let  $ H_- $,  $ H_+ $  be Hopf
$ F $--algebras,  let \  $ \, \pi \colon {(H_-)}_{op} \otimes H_+
\longrightarrow F \, $  \  be an arbitrary Hopf pairing, and let
$ \, D:= D(H_-,H_+,\pi) \, $  be the corresponding quantum
double.  Then there exist  $ F $--algebra isomorphisms
  $$  D^* \cong {H_+}^* \otimeshat {H_-}^* \; ,  \quad  D^*
\cong {H_-}^* \otimeshat {H_+}^*  $$
dual of the  $F$--coalgebra  isomorphisms  $ \, D \cong H_+
\otimes H_- \, $,  $ \, D \cong H_- \otimes H_+ \, $  (cf.~\S 3.1).
$ \square $
\endproclaim

\vskip7pt

   {\bf 5.3  Quantum enveloping algebras as function algebras.}  \  The DRT
pairings induce several linear embeddings, namely
  $$  \matrix
   \uphim \longhookrightarrow {\uphip}^{\! *} \, ,  &
{im}_{\scriptscriptstyle M} \colon \,
{\uphizM M } \longhookrightarrow {\uphizM {M'} }^{\, *} \, ,  &
{\uphiMbm M } \longhookrightarrow {\uphiMbp {M'} }^{\, *}  &  \quad
\hbox{(induced by  $ \, \pi_\varphi \, $)}  \\
   \uphip \longhookrightarrow {\uphim}^{\! *} \, ,  &
{\overline{im}}_{\scriptscriptstyle M} \colon \, {\uphizM M }
\longhookrightarrow {\uphizM {M'} }^{\, *} \, ,  &  {\uphiMbp M }
\longhookrightarrow {\uphiMbm {M'} }^{\, *}  &  \quad  \hbox{(induced by
$ \, \overline{\pi_\varphi} \, $)}  \\
      \endmatrix   \eqno (5.1)  $$
the right-hand-side ones being also embeddings of
formal Hopf algebras.  Therefore we identify the various quantum algebras
with their images in the corresponding dual spaces.

\vskip7pt

\proclaim{Lemma 5.4}                          \hfill\break
   \indent   (a) \  The subset  $ \, \left\{\, \prod_{r=N}^1 {(-1)}^{f_r}
q_{\alpha^r}^{-{f_r \choose 2}} {\left( \fbar^\varphi_{\alpha^r}
\right)}^{f_r} \,\bigg\vert\, f_1, \dots, f_N \in \N \,\right\}
\, $  of  $ \calUphim $  is the pseudobasis of
$ {U_{\scriptscriptstyle +}}^{\! *} $  dual of the
PBW basis of  $ \gerUp $  of decreasing ordered monomials, while the subset
$ \, \left\{\, \prod_{r=N}^1 {(-1)}^{f_r} q_{\alpha^r}^{-{f_r \choose 2}}
{\left( F^\varphi_{\alpha^r} \right)}^{(f_r)} \,\bigg\vert\, f_1, \dots, f_N
\in \N \,\right\} \, $  of  $ \gerUphim $  is the pseudobasis of
$ {U_{\scriptscriptstyle +}}^{\! *} $  dual of the PBW basis of  $ \calUp $
of decreasing ordered monomials.  A similar statement holds with the roles of
$ U_{\scriptscriptstyle -} $  and  $ U_{\scriptscriptstyle +} $  reversed.
                                              \hfill\break
   \indent   (b) \  $ {\calUphizM M } $  (hence  $ {\uphizM M } \, $)
contains the pseudobasis
$ {\Cal B}_{\scriptscriptstyle M} $  (relative to
$ {im}_{\scriptscriptstyle M} $),  resp.~$ {\overline{\Cal
B}}_{\scriptscriptstyle M} $  (relative to
$ {\overline{im}}_{\scriptscriptstyle M} $),  of  $ {\uzM {M'} }^{\, *} $
dual of the PBW basis of  $ \gerUzM{M'} $.
                                              \hfill\break
   \indent   (c) \  $ {\calUphiMbm M } \, $,
resp.~$ {\calUphiMbp M } \, $  (hence  $ {\uphiMbm M } \, $,
resp.~$ {\uphiMbp M } \, $)  contains the pseudobasis of
$ {\uphiMbp {M'} }^{\, *} $,  resp.~of  $ {\uphiMbp {M'} }^{\, *} $,  dual of
the PBW basis of  $ \gerUphiMbm{M'} \, $,  resp.~of  $ \gerUphiMbp{M'} \, $.
The elements of this pseudobasis have form  $ \, {\Cal F}^\varphi \cdot \psi
\, $,  resp.~$ \, \psi \cdot {\Cal E}^\varphi \, $,  where  $ {\Cal
F}^\varphi $,  resp.~$ {\Cal E}^\varphi $,  is an ordered monomial in the
$ \fbar^\varphi_\alpha $'s,  resp.~the  $ \ebar^\varphi_\alpha $'s,  and
$ \, \psi \in \calUphizM{M} $.
\endproclaim

\demo{Proof}  Claim  {\it (a)}  is trivial.
As for  {\it (b)}  and  {\it (c)},  let  $ \, \gerE_\eta \cdot u_\tau
\, $  be any PBW monomial of  $ \, {\gerUphiMbp {M'} } \cong
\gerUp \otimes {\gerUzM {M'} } \, $,  with  $ \, u_\tau :=
\prod_{i=1}^n \left( \Lambda_i; \, 0 \atop t_i \right)
\cdot \Lambda_i^{-Ent({t_i/2})} \, $  ($ \tau = (t_1, \dots, t_n) \in
\N^n \, $)  and  $ \, \gerE_\eta := \prod_{k=N}^1 E_{\alpha^k}^{(e_k)} \, $
($ \eta = (e_1, \dots, e_N) \in \N^N \, $).  Let also  $ \,
\calF^\varphi_\phi := \prod_{k=N}^1 {\left( \fbar^\varphi_{\alpha^k}
\right)}^{\! (f_k)} \, $
($ \, \phi = (f_1, \dots, f_N) $\break
$ \in \N^n \, $)  be any (modified) PBW monomial
of  $ \, \calUphim \, $.  Then (for all  $ \, \mu \in M $,  $ \nu \in M' $)
  $$  \pi_\varphi \left( \calF^\varphi_\phi \cdot L_{-(1+\varphi)(\mu)}, \,
\gerE_\eta \cdot L_\nu \right) = c_\eta \cdot \delta_{\phi,\eta} \cdot
q^{+\left( \mu \vert \nu \right)} \cdot q^{-\left( \mu \vert s \left(
\gerE_\eta \right) \right)}  $$
by (2.3), where  $ \, c_\eta := {(-1)}^{\sum_{k=1}^N e_k} \cdot q^{-\sum_{h<k}
\left( e_h \tau_{\alpha^h} \vert e_k \alpha^k \right)}
\cdot q^{\sum_{k=1}^N d_{\alpha^k} {e_k \choose 2}} \, $  {\sl is
independent of  $ \mu $  and $ \nu $}.  Therefore only PBW monomials
of shape  $ \, \calF^\varphi_\phi \cdot
z \, $  ($\, z \in \uphizM{M} \, $)  give non zero values when paired with
$ \, \gerE_\eta \cdot L_\nu \, $,  hence also with  $ \, \gerE_\eta \cdot
u_\tau \, $.  Now direct computation gives
  $$  {\Big\langle \calF^\varphi_\eta \cdot L_{-(1+\varphi)(\mu)},
\gerE_\eta \cdot u_\tau \Big\rangle}_{\! \pi_\varphi} = c_\eta \cdot
q^{-\left( \mu \vert s \left( \gerE_\eta \right) \right)} \cdot
\prod_{i=1}^n {m_i \choose t_i}_{\!\! q_i} \cdot q^{ - d_i m_i \cdot
Ent(t_i/2)}  \; \qquad  \forall\, \mu, \tau \in \N^n  $$
where we identify  $ \, M_+ \cong \N^n \, $  so that  $ \, M_+ \ni \mu = m_1
\mu_1 + \cdots + m_n \mu_n \cong (m_1, \dots, m_n) \in \N^n \, $.  Then
endowing  $ \N^n $  with the product ordering (of the natural ordering of
$ \N $)  we have
  $$  \eqalign{
   {\Big\langle \calF^\varphi_\eta \cdot L_{-(1+\varphi)(\mu)}, \gerE_\eta
\cdot u_\tau \Big\rangle}_{\! \pi_\varphi}  &  \neq 0  \iff  \tau \preceq \mu
\cr
   {\Big\langle \calF^\varphi_\eta \cdot L_{-(1+\varphi)(\tau)}, \gerE_\eta
\cdot u_\tau \Big\rangle}_{\! \pi_\varphi}  &  = c_\eta \cdot q^{-\left( \tau
\vert s \left( \gerE_\eta \right) \right)} \cdot q^{-T(\tau)}  \qquad \qquad
\qquad \forall\, \tau \in \N^n  \cr }  $$
where  $ \, T(\tau) := \sum_{i=1}^n d_i t_i Ent(t_i/2) \, $;  in particular
$ \, C_{\eta,\tau} := c_\eta \cdot q^{-\left( \mu \vert s \left( \gerE_\eta
\right) \right)} \cdot q^{-T(\tau)} \, $  {\sl is invertible in}  $ \kqqm $.
Thus we have formulas  (for all  $ \, \tau \in \N^n \, $)
  $$  {} \;  \calF^\varphi_\eta \cdot L_{-(1+\varphi)(\tau)} = C_{\eta,\tau}
\cdot {\left( \gerE_\eta \cdot u_\tau \right)}^* + \sum_{\tau' \prec \tau}
{\Big\langle \calF^\varphi_\eta \cdot L_{-(1+\varphi)(\tau)}, \gerE_\eta
\cdot u_{\tau'} \Big\rangle}_{\! \pi_\varphi} \! \cdot {\left( \gerE_\eta
\cdot u_{\tau'} \right)}^*  $$
which tell us that  $ \, \big\{\, \calF^\varphi_\eta \cdot
L_{-(1+\varphi)(\tau)} \,\big\vert\, \tau \in \N^n \,\big\} \, $  is obtained
from  $ \, \big\{\, {\left( \gerE_\eta \cdot u_\tau \right)}^* \,\big\vert\,
\tau \in \N^n \,\big\} \, $  by means of the matrix  $ \, {\Bbb M} :=  {\Big(
{\big\langle \calF^\varphi_\eta \cdot L_{-(1+\varphi)(\tau)}, \, \gerE_\eta
\cdot u_{\tau'} \big\rangle}_{\! \pi_\varphi} \Big)}_{\tau, \tau' \in \N^n} $
which has lower triangular shape, all entries in  $ \kqqm $,  and
diagonal entries  {\sl invertible}  in  $ \kqqm $;  then
the inverse matrix  $ {\Bbb M}^{-1} $  has the same properties, whence
{\it (c)}  follows for  $ \calUphiMbm{M} $.  The same proof applies for
$ \calUphiMbp{M} $  with  $ \overline{\pi_\varphi} $
instead of  $ \pi_\varphi $,  and also gives  {\it (b)}  for
$ \, \eta = 0 \, $.   $ \square $
\enddemo

\vskip7pt

   {\bf 5.5  Remark.}  \;  Since  $ \, D_{\scriptscriptstyle {M'} } \cong
{\uMbp {M'} } \otimes {\uMbm Q } \cong U_{\scriptscriptstyle +} \otimes
{\uzM {M'} } \otimes {\uzM Q } \otimes U_{\scriptscriptstyle -} \, $,  we
have  $ \, {D_{\scriptscriptstyle {M'} }}^* \cong {\uMbp {M'} }^* \otimeshat
{\uMbm Q }^* \cong {U_{\scriptscriptstyle +}}^* \otimeshat {\uzM {M'} }^*
\otimeshat {\uzM Q }^* \otimeshat {U_{\scriptscriptstyle -}}^* \, $;  hence
from Lemma 5.4 we deduce that
                                                    \par
   {\it  Every element  $ \, f \in {D_{\scriptscriptstyle {M'} }}^* \, $  has
a unique expression as formal series
  $$  f = \sum_{{\Cal F}, {\Cal M}, {\Cal L}, {\Cal E}} a_{{\Cal F},
{\Cal M}, {\Cal L}, {\Cal E}} \cdot {\Cal F} \cdot {\Cal M } \cdot {\Cal L}
\cdot  {\Cal E}  $$
in which  $ \, a_{{\Cal F}, {\Cal M}, {\Cal L}, {\Cal E}} \in \kq \, $,
$ \, {\Cal M} \in {\Cal B}_{\scriptscriptstyle M} \, $,  $ \, {\Cal L} \in
{\overline{\Cal B}}_{\scriptscriptstyle P} \, $,  and the
$ {\Cal F}^\varphi $'s,  resp.~the  $ {\Cal E}^\varphi $'s,  are ordered
monomials in the  $ F^\varphi_\alpha $'s,  resp.~in the
$ E^\varphi_\alpha $'s.
                                                    \par
   In particular, every  $ \, f \in {D_{\scriptscriptstyle {M'} }}^*
\, $  can be uniquely expressed as a formal series in the
$ F^\varphi_{\alpha^1}, \dots, F^\varphi_{\alpha^N}, E^\varphi_{\alpha^1},$
\dots, $ E^\varphi_{\alpha^N} $  \, with coefficients in  $ \, {\left( {\uzM
{{M'}} } \otimes {\uzM Q } \right)}^* \cong {\uzM {{M'}} }^* \otimeshat
{\uzM Q }^* \, $.}
                                                    \par
   Similarly the triangular decompositions
$ \, U_{\scriptscriptstyle +} \otimes {\uzM {M'} } \otimes
U_{\scriptscriptstyle -} \cong {\uqMg {M'} } \cong U_{\scriptscriptstyle -}
\otimes {\uzM {M'} } \otimes U_{\scriptscriptstyle +} \, $  give
$ \, {U_{\scriptscriptstyle +}}^{\! *} \otimeshat {\uzM {M'} }^{\! *}
\otimeshat {U_{\scriptscriptstyle -}}^{\! *} \cong {\uqMg {M'} }^{\! *}
\cong {U_{\scriptscriptstyle -}}^{\! *} \otimeshat {\uzM {M'} }^{\! *}
\otimeshat {U_{\scriptscriptstyle +}}^{\! *} \, $,  whence Lemma 5.4 implies
that
                                                    \par
   {\it  Every  $ \, f \in {\uqMg {M'} }^* \, $  can be uniquely expressed as
a formal series in the  $ F^\varphi_{\alpha^1},
                                \dots, F^\varphi_{\alpha^N}, $\break
$ E^\varphi_{\alpha^1}, \dots, E^\varphi_{\alpha^N} \, $  with
coefficients in  $ \, {\uzM {M'} }^* \, $.}
                                                    \par
   In the sequel when considering the composed embedding  $ \, {\uzM M }
\hookrightarrow {\uzM {M'} }^* \hookrightarrow {\uqMg {M'} }^* \, $  we
shall always mean that the first embedding is induced by
$ \overline{\pi_\varphi} $  (cf.~(5.1)).

\vskip7pt

\proclaim{Proposition 5.6}  The monomorphism  $ \; j_{\scriptscriptstyle M}
\colon \, {\uqphiMg {M'} }^* \longhookrightarrow {D_{\scriptscriptstyle
M'}}^{\! *} \; $  (cf.~\S 5.1) is given by
  $$   j_{\scriptscriptstyle M} \colon  \; \quad  F^\varphi_i \mapsto
F^\varphi_i \otimes 1 \, ,  \quad  L_\mu \mapsto L_{-(1+\varphi)(\mu)}
\otimes L_{(1-\varphi)(\mu)} \, ,  \quad  E^\varphi_i \mapsto 1 \otimes
E^\varphi_i  \; \quad  \forall \, i, \mu \, ;   \eqno (5.2)  $$
in particular the image of  $ j_{\scriptscriptstyle M} $  is the closure of
the subalgebra generated by the set
  $$   \Big\{\, F^\varphi_i \otimes 1, \, L_{-(1+\varphi)(\mu)} \otimes
L_{(1-\varphi)(\mu)}, \, 1 \otimes E^\varphi_i \,\Big\vert\, i= \unon,
\, \mu \in M \,\Big\} \; .  $$
\endproclaim

\demo{Proof}  For PBW monomials we have  $ \, {pr}_{\scriptscriptstyle M}
\Big( E \cdot L \otimes K \cdot F \Big) =  E \cdot L \cdot K \cdot F \, $;
therefore (5.2) comes out of the definition  $ \, j_{\scriptscriptstyle M} :=
{\left( {pr}_{\scriptscriptstyle M} \right)}^{\! *} \, $.  As an example
  $$  \displaylines{
   \Big\langle j_{\scriptscriptstyle M} \big( L_\mu \big), E \cdot L_\nu
\otimes K_\alpha \cdot F \Big\rangle = {\Big\langle L_\mu, E \cdot L_\nu
\cdot K_\alpha \cdot F \Big\rangle}_{\overline{\pi}_\varphi} = \delta_{E,1}
\cdot \delta_{E,1} \cdot q^{(\mu \vert \nu + \alpha)}  \cr
   {\Big\langle L_{-(1+\varphi)(\mu)} \otimes L_{(1-\varphi)(\mu)}, E \cdot
L_\nu \otimes K_\alpha \cdot F \Big\rangle}_{\! \pi_\varphi \otimes
\overline{\pi}_\varphi} = \delta_{E,1} \cdot \delta_{F,1} \cdot
q^{(\mu \vert \nu + \alpha)}  \cr }  $$
whence  $ \, j_{\scriptscriptstyle M} \big( L_\mu \big) =
L_{-(1+\varphi)(\mu)}
\otimes L_{(1-\varphi)(\mu)} \, $.  Since  $ \, j_{\scriptscriptstyle M}
:= {\left( {pr}_{\scriptscriptstyle M'} \right)}^{\! *} \, $  is continuous
(cf.~\S 1.1), by Lemma 5.4 and Remark 5.5 it is uniquely determined by (5.2).
$ \square $
\enddemo

\vskip7pt

{\bf  Remark 5.7.}  \  Now we can identify  $ j_{\scriptscriptstyle M}
\left( {\uqMg {M'} }^* \right) $  with
the space of formal series in the  $ F^\varphi_{\alpha^1}, \dots,
F^\varphi_{\alpha^N}, E^\varphi_{\alpha^1}, \dots, E^\varphi_{\alpha^N} $
\, with coefficients in  $ \, j_{\scriptscriptstyle M} \left( {\uzM {M'} }^*
\right) \, $.  In order to locate the image   --- under
$ j_{\scriptscriptstyle M} $  ---   of the pseudobasis of
$ j_{\scriptscriptstyle M} \left( {\uqphiMg {M'} }^* \right) $
dual of the PBW basis of  $ {\gerUphiMg {M'} } $,  let
  $$  \displaylines{
   \gerE_\eta := \prod_{k=N}^1 E_{\alpha^k}^{(e_k)} \, ,  \quad \,
u_\tau := \prod_{i=1}^n \left( \Lambda_i; 0 \atop t_i \right) \cdot
\Lambda_i^{-Ent({t_i/2})} \, ,  \quad \, \gerF_\phi := \prod_{k=1}^N
F_{\alpha^k}^{(f_k)}  \cr
   X_{\scriptscriptstyle \eta, \tau, \phi} :=
\gerE_\eta \cdot u_\tau \cdot \gerF_\phi \, ,  \quad \, \calF_\phi :=
\prod_{k=N}^1 {\left( \fbar^\varphi_{\alpha^k} \right)}^{f_k} \, ,
\quad \, \calE_\eta := \prod_{k=1}^N {\left( \ebar^\varphi_{\alpha^k}
\right)}^{e_k}  \cr }  $$
and  $ \, L^{\varphi,\otimes}_\mu := L_{-(1+\varphi)(\mu)}
\otimes L_{(1-\varphi)(\mu)} \, $.  Then (2.3) gives (for some
$ \, a, b \in \Z \, $  and  $ \, \varepsilon = \pm 1 \, $)
  $$  \displaylines{
   {\ } \quad   {\left\langle \calF^\varphi_\phi \cdot
L^{\varphi,\otimes}_\mu
\cdot \calE^\varphi_\eta, \, \gerE_{\bar\eta}
\otimes u_\tau \cdot \gerF_{\bar\phi} \right\rangle}_{\! \pi_\varphi \otimes
\overline{\pi_\varphi}} =   \hfill {\ }  \cr
   {\ } \hfill   = \delta_{\phi, \bar\eta} \, \delta_{\eta, \bar\phi} \cdot
\varepsilon \, q^{ a + b - \left( \mu \vert s \left( \gerE_{\bar\eta} \right)
\right) - \left( \mu \vert s \left( \gerF_{\bar\phi} \right) \right)} \cdot
\prod_{i=1}^n {m_i \choose t_i}_{\!\! q_i} q^{- d_i m_i \cdot Ent(t_i/2)}
\quad {\ } \cr }  $$
   \indent   Thus among the elements of the form  $ \, \calF^\varphi_\phi
\cdot L^{\varphi,\otimes}_\mu \cdot \calE^\varphi_\eta \, $  only those with
$ \, (\phi, \eta) = (\bar\eta, \bar\phi) \, $  and  $ \, \mu \preceq
\tau \, $  takes non-zero value on  $ \, \gerE_{\bar\eta} \otimes
u_\tau \cdot \gerF_{\bar\phi} \, $.  Therefore
  $$  \displaylines{
   {\ }  \calF^\varphi_{\bar\eta} \cdot L^{\varphi,\otimes}_\tau
\cdot \calE^\varphi_{\bar\phi} = \varepsilon \, q^z \cdot
{X_{\scriptscriptstyle \bar\eta, \tau,
\bar\phi}}^{\! *} + \sum_{\tau' \prec \tau} \Big\langle
\calF^\varphi_{\bar\eta} \cdot L^{\varphi,\otimes}_\tau
\cdot \calE^\varphi_{\bar\phi}, \,
X_{\scriptscriptstyle \bar\eta, \tau', \bar\phi} \Big\rangle \cdot
{X_{\scriptscriptstyle \bar\eta, \tau', \bar\phi}}^{\! *} =  \cr
   {} \hfill   = \varepsilon \, q^z \cdot {X_{\scriptscriptstyle \bar\eta,
\tau, \bar\phi}}^{\! *} + \sum_{\tau' \prec \tau} c_{\tau,\tau'} \cdot
{X_{\scriptscriptstyle \bar\eta, \tau', \bar\phi}}^{\! *}  {\ }  \cr }  $$
(with  $ \, z \in \Z \, $,  $ \, c_{\tau,\tau'} \in \kqqm \, $;  we set also
$ \, c_{\tau,\tau}:= \varepsilon \, q^z \, $  and  $ \, c_{\tau,\tau'}:= 0 \, $
for  $ \, \tau' \not\prec \tau \, $);  then we turn from
$ \, \big\{\, {X_{\scriptscriptstyle \bar\eta, \tau',
\bar\phi}}^{\! *} \,\big\vert\, \tau' \in \N^n \,\big\} \, $  to
$ \, \big\{\, \calF^\varphi_{\bar\eta} \cdot L^{\varphi,\otimes}_\tau
\cdot \calE^\varphi_{\bar\phi} \,\big\vert\, \tau \in \N^n \,\big\} \, $
by means of a lower triangular matrix
$ \, {\Bbb M}_{\scriptscriptstyle \bar\eta, \bar\phi} :=  {\big( c_{\tau,\tau'}
\big)}_{\tau, \tau' \in \N^n} \, $,  whose entries belong to  $ \kqqm $  and
whose diagonal entries are  {\sl invertible\/}  in  $ \kqqm $;  then letting
$ \, {\left( {\Bbb M}_{\scriptscriptstyle \bar\eta, \bar\phi} \right)}^{\! -1}
= {\left( c'_{\tau, \tau'} \right)}_{\tau, \tau' \in \N^n} \, $  we find
that  $ \, {X_{\scriptscriptstyle \eta, \tau, \phi}}^{\! *} =
\sum_{\tau' \preceq \tau} c'_{\tau, \tau'} \cdot \calF^\varphi_{\bar\eta}
\cdot L^{\varphi,\otimes}_{\tau'} \cdot \calE^\varphi_{\bar\phi} \, $.  Now
let  $ \, B^{\varphi, \otimes}_{\scriptscriptstyle \bar\eta, \tau, \bar\phi} :=
\sum_{\tau' \preceq \tau} c'_{\tau, \tau'} \cdot
L^{\varphi,\otimes}_{\tau'} \, $:  then
  $$  {X_{\scriptscriptstyle \bar\eta, \tau, \bar\phi}}^{\! *} =
\calF^\varphi_{\bar\eta} \cdot B^{\varphi, \otimes}_{\scriptscriptstyle
\bar\eta, \tau, \bar\phi} \cdot \calE^\varphi_{\bar\phi}   \eqno (5.3)  $$
thus  $ \, \big\{\, \calF^\varphi_{\bar\eta} \cdot B^{\varphi,
\otimes}_{\scriptscriptstyle \bar\eta, \tau, \bar\phi} \cdot
\calE^\varphi_{\bar\phi} \,\big\vert\, \bar\eta \in \N^N, \, \tau \in \N^n,
\, \bar\phi \in \N^N \,\big\} \, $  is the image pseudobasis
(of  $ j_{\scriptscriptstyle M} \left( {\uqphiMg {M'\!} } \right) $)  we were
looking for; in particular we stress the fact that
                                                      \par
   {\it  The pseudobasis of  $ \, j_{\scriptscriptstyle M} \left( {\uqphiMg
{M'\!} }^* \right) \, $  dual of the PBW basis of  $ {\gerUphiMg {M'\!} } $
                     is contained in\break
$ \, j_{\scriptscriptstyle M} \left( {\uqphiMg {M'\!} }^* \right) \cap \left(
\calUphim \otimes {\calUphizM M } \otimes \calUphip \right) \, $.}

\vskip7pt

   {\bf 5.8  Integer forms.}  \;  We want to study the subspaces of linear
functions on  $ {\uqphiMg {M'} } $  which are "integer-valued" on its
integer forms.  Thus we define
  $$  \displaylines {
   {\gerUphiMg {M'} }^* := \Big\{\, f \in {\uqphiMg
{M'} }^* \,\Big\vert\, \left\langle f, {\gerUphiMg {M'} } \right\rangle
\subseteq \kqqm \,\Big\}  \cr
   {\calUphiMg {M'} }^* := \Big\{\, f \in {\uqphiMg {M'} }^* \,\Big\vert\,
\left\langle f, {\calUphiMg {M'} } \right\rangle \subseteq \kqqm \,\Big\}  \cr
   {\frak I}_\varphi^{\scriptscriptstyle M} := \Big\{\, f \in
j_{\scriptscriptstyle M} \Big({\uqphiMg {M'} }^* \Big) \,\Big\vert\,
\left\langle f, {\gerUphiMg {M'} } \right\rangle \subseteq \kqqm \,\Big\}  \cr
   {\Cal I}_\varphi^{\scriptscriptstyle M} := \Big\{\, f \in
j_{\scriptscriptstyle M} \Big({\uqphiMg {M'} }^* \Big) \,\Big\vert\,
\left\langle f, {\calUphiMg {M'} } \right\rangle \subseteq \kqqm \,\Big\}  \cr
   {\gerUphizM {M'} }^* := \Big\{\, f \in {\uphizM {M'} }^* \,\Big\vert\,
\left\langle f, {\gerUphizM {M'} } \right\rangle \subseteq \kqqm \,\Big\}  \cr
   {\calUphizM {M'} }^* := \Big\{ f \in {\uphizM{M'} }^* \,\Big\vert\,
\left\langle f, {\calUphizM {M'} } \right\rangle \subseteq \kqqm \,\Big\}
\cr }  $$
notice that  $ j_{\scriptscriptstyle M} $  restricts
to isomorphisms  $ \; j_{\scriptscriptstyle M} \colon \, {\gerUphiMg
{M'} }^* \, {\buildrel \cong \over \loongrightarrow} \, {\frak
I}_\varphi^{\scriptscriptstyle M} \, $,  $ \; j_{\scriptscriptstyle M}
\colon \, {\calUphiMg {M'} }^* \, {\buildrel \cong \over \loongrightarrow}
\, {\Cal I}_\varphi^{\scriptscriptstyle M} \, $.

\vskip7pt

\proclaim{Proposition 5.9}
                                      \hfill\break
   \indent   (a)  $ {\gerUphiMg {M'} }^* $  is the
$ \kqqm $--submodule  (of  $ {\uqphiMg {M'} }^* \, $)
                        of formal series (cf.~\S 5.5)\break
$ \; \sum_{{\Cal F}^\varphi, \psi, {\Cal E}^\varphi} {\Cal F}^\varphi \cdot
\psi \cdot {\Cal E}^\varphi \; $  in which  $ \, \psi \in {\gerUphizM {M'} }^*
\, $  and the  $ {\Cal F}^\varphi $'s,  resp.~the
$ {\Cal E}^\varphi $'s,  are monomials of the PBW basis of
$ \calUphim $,  resp.~of  $ \calUphip $.
                                     \hfill\break
   \indent   In particular  $ {\gerUphiMg {M'} }^* $  is a formal Hopf
subalgebra of  $ {\uqphiMg {M'} }^* $.
                                      \hfill\break
   \indent   (b)  $ {\calUphiMg {M'} }^* $  is the  $ \kqqm $--submodule  (of
$ {\uqphiMg {M'} }^* \, $)
                        of formal series (cf.~\S 5.5)\break
$ \; \sum_{{\frak F}^\varphi, \phi, {\frak E}^\varphi} {\frak F}^\varphi
\cdot \phi \cdot {\frak E}^\varphi \; $  in which  $ \, \phi \in
{\calUphizM {M'} }^* \, $  and the  $ {\frak F}^\varphi $'s,  resp.~the
$ {\frak E}^\varphi $'s,  are monomials of the PBW basis of
$ \gerUphim $,  resp.~of  $ \gerUphip $.
                                     \hfill\break
   \indent   In particular  $ {\calUphiMg {M'} }^* $  is a formal Hopf
subalgebra of  $ {\uqphiMg {M'} }^* $.
\endproclaim

\demo{Proof}  Let us prove  {\it (b)}.  Let  $ \, f \in {\uqphiMg {M'} }^*
\, $  be given, and expand
it as a series  $ \; f = \sum_{\phi, \eta \in \N^N} \gerF^\varphi_\phi \cdot
\Phi^\varphi_{\phi,\eta} \cdot \gerE^\varphi_\eta \; $  in which the
$ \gerF^\varphi_\phi $'s,  resp.~the
$ \gerE^\varphi_\eta $'s,  are PBW monomials of
$ \gerUphim $,  resp.~of  $ \gerUphip $,  and  $ \, \Phi^\varphi_{\phi,\eta}
\in {\uphizM {M'} }^* \, $.  Let  $ \, \Phi^\varphi_{\phi,\eta} = \sum_{\tau
\in \N^n} a_{\phi,\eta}^\tau B^\varphi_{\scriptscriptstyle \phi, \tau, \eta}
\, $  and  $ \, j_{\scriptscriptstyle M} \big( B^\varphi_{\scriptscriptstyle
\phi, \tau, \eta} \big) = \sum_{\mu \preceq \tau}
c_\tau^\mu \cdot L^{\varphi,\otimes}_\mu = B^{\varphi,
\otimes}_{\scriptscriptstyle \phi, \tau, \eta} \, $.  For all monomials
$ \, \calE_{\bar \eta} \cdot L_\nu \cdot \calF_{\bar\phi} \, $  of a PBW
basis
of  $ {\calUphiMg {M'} } $  we have
  $$  \displaylines{
   \Big\langle f, \, \calE_{\bar \eta} \cdot L_\nu \cdot \calF_{\bar\phi}
\Big\rangle = \sum_{\phi, \tau, \eta} a_{\phi,\eta}^\tau \cdot \left\langle
\gerF^\varphi_\phi \cdot B^\varphi_\tau \cdot
\gerE^\varphi_\eta, \, \calE_{\bar \eta} \cdot L_\nu \cdot \calF_{\bar\phi}
\right\rangle =  \cr
   = \sum_{\phi, \tau, \eta} a_{\phi,\eta}^\tau \cdot \sum_{\mu \preceq \tau}
c_\tau^\mu \cdot {\left\langle \gerF^\varphi_\phi \cdot
L_{-(1+\varphi)(\mu)}, \, \calE_{\bar \eta} \cdot L_\nu
\right\rangle}_{\! \pi_\varphi} \! \cdot {\left\langle
L_{(1-\varphi)(\mu)} \cdot \gerE^\varphi_\eta, \,
\calF_{\bar\phi} \right\rangle}_{\! \overline{\pi_\varphi}} =  \cr
   = \pm q^{a+b} \cdot \sum_\tau a_{\bar\eta,
\bar\phi}^\tau \cdot \left\langle B^\varphi_\tau,
L_{\nu - \left( s \left( \gerE_{\bar\phi} \right) +
s \left( \gerF_{\bar\eta} \right) \right)} \right\rangle
= \pm q^{a+b} \cdot \left\langle \Phi^\varphi_{\bar\eta, \bar\phi},
L_{\nu - \left( s \left( \gerE_{\bar\phi} \right) + s \left(
\gerF_{\bar\eta} \right) \right)} \right\rangle  \cr }  $$
for some  $ \, a, b \in \Z \, $  depending only respectively on
$ \bar\eta $  and  $ \bar\phi $.  Then if
$ \, \Phi^\varphi_{\bar\eta,\bar\phi} \in {\calUphizM {M'} }^* \, $  we
have  $ \, \left\langle f, \calE_{\bar\eta} \cdot L_\nu \cdot
\calF_{\bar\phi} \right\rangle \in \kqqm \, $  for all  $ \bar\eta $,
$ \nu $,  $ \bar\phi $,  hence  $ \, f \in {\calUphiMg {M'} }^* \, $;
conversely, the latter gives  $ \, \left\langle \Phi^\varphi_{\bar\eta,
\bar\phi}, L_{\nu'} \right\rangle \in \kqqm \, $  for all
$ \, \nu' \in M' \, $,  hence  $ \, \Phi^\varphi_{\bar\eta,\bar\phi}
\in {\calUphizM {M'} }^* \, $.
                                                   \par
   Now consider the Hopf structure.  Let  $ \, f \in
{\calUphiMg {M'} }^* \, $,  and expand  $ \Delta(f) $
                                        as a series\break
$ \; \Delta(f) = \sum_\sigma \left( \gerF^\varphi_\sigma \cdot \phi_\sigma
\cdot \gerE^\varphi_\sigma \right) \otimes \left( {\gerF^\varphi_\sigma}' \cdot
\phi'_\sigma \cdot {\gerE^\varphi_\sigma}' \right) \; $  so that  $ \,
\phi_\sigma \otimes \phi'_\sigma \neq \phi_\tau \otimes \phi'_\tau \, $  for
all  $ \sigma $,  $ \tau $,  such that  $ \, \left( \gerF^\varphi_\sigma,
\gerE^\varphi_\sigma, {\gerF^\varphi_\sigma}', {\gerE^\varphi_\sigma}' \right)
\neq \left( \gerF^\varphi_\tau, \gerE^\varphi_\tau, {\gerF^\varphi_\tau}',
{\gerE^\varphi_\tau}' \right) \, $  (this is always possible).
As  $ \, f \in {\calUphiMg {M'} } \, $,  then  $ \Delta(f) $  is integer-valued
on  $ \, {\calUphiMg {M'} } \otimes {\calUphiMg {M'} } $.
Fix any  $ \bar\sigma $:  exploiting (2.3) we get the existence
of unique (non-modified) PBW monomials  $ \calE_{\bar \sigma} $,
$ \calF_{\bar \sigma} $,  $ \calE'_{\bar \sigma} $,
$ \calF'_{\bar\sigma} $  such that
  $$  \left\langle \Delta(f), \big( \calE_{\bar\sigma} \otimes
\calE'_{\bar\sigma} \big) \cdot \big( L_\nu \otimes L_{\nu'} \big) \cdot
\big( \calF_{\bar\sigma} \otimes \calF'_{\bar\sigma} \big) \right\rangle =
\pm q^{c_{\bar\sigma}} \cdot \big\langle \phi_\sigma \otimes \phi'_\sigma,
L_{\nu + \xi} \otimes L_{\nu' + \xi'} \big\rangle  $$
for all  $ \, \nu, \nu' \in M' \, $  (for some  $ \, c_{\bar\sigma} \in
\Z \, $  and  $ \, \xi, \xi' \, \in Q \, (\subseteq M') \, $
independent of  $ \nu $,  $ \nu' \, $);
since  $ \Delta(f) $  is integer-valued,  $ \phi_{\bar \sigma} \otimes
\phi'_{\bar \sigma} $  is integer-valued on  $ {\calUphizM {M'} } \otimes
{\calUphizM {M'} } $,  that is  $ \, \phi_{\bar \sigma} \otimes
\phi'_{\bar \sigma} \in {\left( {\calUphizM {M'} } \otimes {\calUphizM {M'} }
\right)}^* = {\calUphizM {M'} }^* \otimeshat {\calUphizM {M'} }^* \, $;
but  $ \, \phi_{\bar \sigma} \otimes \phi'_{\bar
\sigma} \in {\uphizM {M'} }^* \otimes {\uphizM {M'} }^* \, $,
thus  $ \, \phi_{\bar \sigma}, \phi'_{\bar
\sigma} \in {\calUphizM {M'} }^* $,  \, q.e.d.
                                               \par
   Finally, we have  $ \, 1 \in {\calUphiMg {M'} }^* \, $,  because  $ \, 1
:= \epsilon \, $,  $ \, \epsilon \big( {\calUphiMg {M'} }^*
\big) \subseteq \kqqm \, $  because  $ \, \epsilon := 1^* \, $  and  $ \, 1
\in {\calUphiMg {M'} } \, $,  and  $ \, S \big( {\calUphiMg {M'} }^* \big) =
{\calUphiMg {M'} }^* \, $  because  $ \, S := S^* \, $  and  $ \, S \big(
{\calUphiMg {M'} } \big) = {\calUphiMg {M'} } \, $.  Thus  $ {\calUphiMg
{M'} }^* $  is a formal Hopf subalgebra of  $ {\uqphiMg {M'} }^* $,
q.e.d.   $ \square $
\enddemo

\vskip7pt

\proclaim{Definition 5.10} \  We call
$ A_\varphi^{\scriptscriptstyle M} $  the subalgebra of
$ \, {\uphiMbm M } \otimes {\uphiMbp P } \; \big( \, \subset
{D_\varphi^{\scriptscriptstyle M'}}^* \, \big) \, $  generated by
$ \; \big\{ \, F^\varphi_i \otimes 1, \, L^{\varphi,\otimes}_\mu,
\, 1 \otimes E^\varphi_i \,\big\vert\, i=\unon; \, \mu \in M \,\big\} \, $.
Then we set
  $$  \eqalign{
   {\frak A}_\varphi^{\scriptscriptstyle M}  &  := \big\{\, f \in
A_\varphi^{\scriptscriptstyle M} \,\big\vert\, \big\langle f,
{\gerUphiMg {M'} } \big\rangle \subseteq \kqqm \,\big\} =
A_\varphi^{\scriptscriptstyle M} \cap {\frak I}_\varphi^{\scriptscriptstyle M}
\cr
   {\Cal A}_\varphi^{\scriptscriptstyle M}  &  := \big\{\, f \in
A_\varphi^{\scriptscriptstyle M} \,\big\vert\, \big\langle f,
{\calUphiMg {M'} } \big\rangle \subseteq \kqqm \,\big\} =
A_\varphi^{\scriptscriptstyle M} \cap {\Cal I}_\varphi^{\scriptscriptstyle M}
\, .  \cr }  $$
\endproclaim

\vskip7pt

\proclaim{Lemma 5.11}
                                               \hfill\break
   \indent   (a) \  $ {\frak A}_\varphi^{\scriptscriptstyle M} $  is a
$ \kqqm $--integer  form of  $ A_\varphi^{\scriptscriptstyle M} $,
generated as a  $ \kqqm $--subalgebra by
  $$  \left\{\, \fbar^\varphi_{\alpha^h} \otimes 1, \,
L^{\varphi,\otimes}_\mu, \, 1 \otimes \ebar^\varphi_{\alpha^k}
\,\Big\vert\, h, k =1, \dots, N; \, \mu \in M \,\right\} \, .  $$
   \indent   (b) \  $ {\Cal A}_\varphi^{\scriptscriptstyle M} $  is a
$ \kqqm $--integer  form of  $ A_\varphi^{\scriptscriptstyle M} $,
generated as a  $ \kqqm $--subalgebra by
  $$  \displaylines{
   \bigg\{\, {\left( F^\varphi_{\alpha^h} \right)}^{(a)} \otimes 1, \,
\left( L^{\varphi,\otimes}_{\mu_i}; \, c \atop t \right),
\, L^{\varphi,\otimes}_{-\mu_i}, \, 1 \otimes {\left( E^\varphi_{\alpha^k}
\right)}^{(d)} \,\bigg\vert\, h, k, i= \unon; \, a, t, d \in \N; \, c \in
\Z \,\bigg\} \, .  \cr }  $$
\endproclaim

\demo{Proof}  Definitions yield a linear isomorphism  $ \,
\Phi_{\scriptscriptstyle M} \colon A_\varphi^{\scriptscriptstyle M} \,
{\buildrel \cong \over \longrightarrow} \, \uphim \otimes {\uphizM M }
\otimes \uphip \, $
                                       given by\break
$ \; \Phi_{\scriptscriptstyle M} \colon  \; F^\varphi_i \otimes 1
\mapsto F^\varphi_i \otimes 1 \otimes 1 \, ,  \; L^{\varphi,\otimes}_\mu
\mapsto 1 \otimes L_\mu \otimes 1 \, ,  \; 1 \otimes E^\varphi_i \mapsto 1
\otimes 1 \otimes E^\varphi_i \; $;  but this restricts to  $ \, \Phi \colon
\, {\Cal A}_\varphi^{\scriptscriptstyle M} \, {\buildrel \cong \over
\rightarrow} \, \calUphim \otimes {\calUphizM M } \otimes \calUphip \, $,
$ \, \Phi \colon \, {\frak A}_\varphi^{\scriptscriptstyle M} \, {\buildrel
\cong \over \rightarrow} \, \gerUphim \otimes {\gerUphizM M } \otimes
\gerUphip \, $,  so \S 3.4 gives the claim.   $ \square $
\enddemo

\vskip7pt

   The following result stems from [C-V2], Lemma 2.5 (which extends [D-L],
Lemma 4.3), relating our quantum formal groups to quantum function algebras;
in particular we prove that  $ {\gerFphiMg M } $  and  $ {\calFphiMg M } $
are integer forms  (over  $ \kqqm $)  of  $ {\fqphiMg M } $  {\sl as Hopf
algebras}.

\vskip7pt

\proclaim{Proposition 5.12}
                                               \hfill\break
   \indent   (a) \  The monomorphism of formal Hopf algebras
$ \,  j_{\scriptscriptstyle M} \colon \, {\uqphiMg {M'} }^*
\longhookrightarrow {D_\varphi^{\scriptscriptstyle M'}}^* \, $
restricts to an embedding  $ \, \mu_{\scriptscriptstyle M} \colon \,
{\fqphiMg M } \llonghookrightarrow {D_\varphi^{\scriptscriptstyle M'}}^* \, $
whose image is contained in  $ A_\varphi^{\scriptscriptstyle M} $.
                                               \hfill\break
   \indent   (b) \  The embedding in (a) preserves integer forms, namely
$ \, {\gerFphiMg M } = {\mu_{\scriptscriptstyle M}}^{\! -1}
\big( {\frak A}_\varphi^{\scriptscriptstyle M} \big)
\, $,  $ \, {\calFphiMg M } = {\mu_{\scriptscriptstyle M}}^{\! -1}
\big( {\Cal A}_\varphi^{\scriptscriptstyle M} \big) \, $,  so that
restriction provides embeddings of      $ \kqqm $--algebras\break
$ \, \mu_{\scriptscriptstyle M} \colon \, {\gerFphiMg M }
\llonghookrightarrow {\frak A}_\varphi^{\scriptscriptstyle M} \, $,
$ \, \mu_{\scriptscriptstyle M} \colon \, {\calFphiMg M }
\llonghookrightarrow {\Cal A}_\varphi^{\scriptscriptstyle M} \, $.
It follows that
                                               \hfill\break
   \centerline{ $ {\gerFphiMg M } $  is a Hopf subalgebra of
$ {\fqphiMg M } $,  and a  $ \kqqm $--integer  form of it, }
                                               \hfill\break
   \centerline{ $ {\calFphiMg M } $  is a Hopf subalgebra of
$ {\fqphiMg M } $,  and a  $ \kqqm $--integer  form of it. }
\endproclaim

\demo{Proof} \, {\it (a)} \  The first part is obvious.  As for the second,
recall that the identification  $ \, D_{\scriptscriptstyle M'} =
{\uphiMbp {M'} } \otimes {\uphiMbm Q } \, $  is given by
$ \; {\uphiMbp {M'} } \otimes {\uphiMbm Q } \, {\buildrel j_+ \otimes j_-
\over \lllonghookrightarrow} \, D_{\scriptscriptstyle M'} \otimes
D_{\scriptscriptstyle {M'}} \, {\buildrel {m_D} \over \longrightarrow} \,
D_{\scriptscriptstyle M'} \; $  where  $ \, j_+ \colon \, {\uphiMbp {M'} }
\hookrightarrow D_{\scriptscriptstyle M'} \, $  and
$ \, j_- \colon \, {\uphiMbm Q } \hookrightarrow
D_{\scriptscriptstyle {M'}} \, $  are the natural Hopf
algebra embeddings,  $ m_D $  is the multiplication of
$ D_{\scriptscriptstyle {M'}} \, $,  and we look this composition
as a Hopf algebra isomorfism; then the identification
$ \, {D_{\scriptscriptstyle M'}}^* =
{\uphiMbp {M'} }^* \otimeshat {\uphiMbm Q }^*  \, $  is given by  $ \,
{\big( m_D \smallcirc (j_+ \otimes j_-) \big)}^* =
\left( j_+^* \otimeshat j_-^* \right) \smallcirc m_D^* \, $.
If  $ m_U $  is the multiplication of
$ {\uqphiMg {M'} } $,  we have  $ \, m_U \smallcirc
(pr_{\scriptscriptstyle M'} \otimes pr_{\scriptscriptstyle M'}) =
pr_{\scriptscriptstyle M'} \smallcirc m_D \, $,  hence dualizing
yields  $ \, {\big( pr_{\scriptscriptstyle M'} \smallcirc m_D \smallcirc
(j_+ \otimes j_-) \big)}^* = {\left( pr_{\scriptscriptstyle M'} \smallcirc
j_+ \right)}^* \otimeshat {\left( pr_{\scriptscriptstyle M'} \smallcirc j_-
\right)}^* \smallcirc m_U^* \, $;  but
          $ \, pr_{\scriptscriptstyle M'} \smallcirc j_\pm = $\break
\noindent   $ = i_\pm \colon U_{q,\varphi}^{\scriptscriptstyle {M'\!}}
(\gerb_\pm) \longhookrightarrow {\uqMg {M'\!} } \, $  (the natural
embedding), thus  $ \; {\big( pr_{\scriptscriptstyle M'}
\smallcirc m_U \smallcirc (j_+ \otimes j_-) \big)}^* = \left( i_+^*
\otimeshat i_-^* \right) \smallcirc m_U^* \, $.  Now  $ m_U^* $  is
the comultiplication  $ \Delta $ of  $ {\uqphiMg {M'} }^* $,  which restricts
to  $ {\fqphiMg M } $,  while  $ \, \rho_\pm :=
i_\pm^* \colon {\uqphiMg {M'} }^* \rightarrow
{U_{q,\varphi}^{\scriptscriptstyle {M'\!}} (\gerb_\pm)}^* \, $
is the "restriction" map, which maps  $ {\fqphiMg M } $  onto
$ \fqphiMbs {M} {\pm} $;  using also (4.1), we obtain
  $$   \displaylines{
   {\big( pr_{\scriptscriptstyle {M'}} \smallcirc m_U \smallcirc
(j_+ \otimes j_-) \big)}^* \left( \fqphiMg{M} \right) =
\left( \rho_+ \otimeshat \rho_- \right) \Big( \Delta
\left( {\fqphiMg M } \right) \Big) \subseteq {\uphiMbm M } \otimes
{\uphiMbp M } \; ;  \quad {}  \cr }  $$
in other words,  $ j_{\scriptscriptstyle M} $  maps
$ {\fqphiMg M } $  into  $ \, {\uphiMbm M } \otimes {\uphiMbp P } \, $.
From the very definition we get that  $ \mu_{\scriptscriptstyle M}
\left( {\fqphiMg M } \right) $  vanishes on  $ \gerK_{\scriptscriptstyle
M'}^\varphi $,  hence   --- by Proposition 5.6 ---
$ \, \mu_{\scriptscriptstyle M} \left( {\fqphiMg M } \right)
\subseteq A_{\scriptscriptstyle M} \, $,  q.e.d.
                                                  \par
   {\it (b)} \  The first two claims are obvious by definition.  Let now, for
instance,  $ \, f \in
{\calFphiMg M } \, $:  then  $ \, \mu_{\scriptscriptstyle M} \big( S(f)
\big) \in A_\varphi^{\scriptscriptstyle M} \, $  by  {\it (a)},
$ \, \mu_{\scriptscriptstyle M} \big( S(f) \big) = S \big(
\mu_{\scriptscriptstyle M}(f) \big) \, $,  and
$ \, \Big\langle S \big( \mu_{\scriptscriptstyle M}(f) \big), {\calUphiMg
{M'} } \Big\rangle = \Big\langle \mu_{\scriptscriptstyle M}(f), S \big(
{\calUphiMg {M'} } \big) \Big\rangle = \Big\langle
\mu_{\scriptscriptstyle M}(f), {\calUphiMg {M'} } \Big\rangle \subseteq
\kqqm \, $,  hence  $ \, \mu_{\scriptscriptstyle M} \big( S(f) \big) \in{\frak A}_\varphi^{\scriptscriptstyle M} \, $,
thus  $ \, S(f) \in {\calFphiMg M } \, $;  similarly,
$ \, \Delta \big( \mu_{\scriptscriptstyle M} (f) \big) \in
A_\varphi^{\scriptscriptstyle M} \otimes
A_\varphi^{\scriptscriptstyle M} \, $,
and  $ \, \Big\langle \Delta \big( \mu_{\scriptscriptstyle M}(f) \big),
{\calUphiMg {M'} } \otimes {\calUphiMg {M'} }
\Big\rangle \subseteq \kqqm \, $,
hence  $ \, (\mu_{\scriptscriptstyle M} \otimes \mu_{\scriptscriptstyle M})
\big( \Delta(f) \big) \in {\frak A}_\varphi^{\scriptscriptstyle M} \otimes
{\frak A}_\varphi^{\scriptscriptstyle M} \, $:  we have only to remark that
  $$  \displaylines{
   \big( A_\varphi^{\scriptscriptstyle M} \otimes
A_\varphi^{\scriptscriptstyle M} \big) \cap
\left\{\, \phi \in j_{\scriptscriptstyle M} \left( {\left(
{\left( \uqphiMg{M'} \right)}^{\otimes 2} \right)}^* \right)
\,\Big\vert\, \left\langle \phi, {\big( \calUphiMg{M'} \big)}^{\otimes 2}
\right\rangle \subseteq \kqqm \right\} = {\Cal
A}_\varphi^{\scriptscriptstyle M} \otimes
{\Cal A}_\varphi^{\scriptscriptstyle M} \, ;  \cr }  $$
we conclude that  $ \, \Delta(f) \in {\calFphiMg M } \otimes
{\calFphiMg M } \, $.  Therefore  $ {\calFphiMg M } $  is a  $ \kqqm $--Hopf
subalgebra of  $ {\fqphiMg M } $.  Finally, let  $ \, f \in {\fqphiMg M }
\, $;  then  $ \, \mu_{\scriptscriptstyle M} \big( c(q) f \big) = c(q)
\cdot \mu_{\scriptscriptstyle M} (f) \in
{\Cal A}_\varphi^{\scriptscriptstyle M} \, $  for some  $ \, c(q) \in
\kqqm \, $.  Thus  $ \, c(q) f \in {\mu_{\scriptscriptstyle M}}^{\! -1}
\big( {\Cal A}_\varphi^{\scriptscriptstyle M} \big) =
{\calFphiMg M } \, $,  and  $ \, f = {{\, 1 \,} \over {\, c(q) \,}} \cdot
\big( c(q) f \big) \, $  with  $ \, c(q) f \in {\gerFphiMg M } \, $:  hence
$ \, \kq \otimes_{k \left[ q, \qm \right]} {\gerFphiMg M } = {\fqphiMg M }
\, $,  i.~e.~$ {\gerFphiMg M } $  is  $ \kqqm $--integer  form of
$ {\fqphiMg M } $,  q.e.d.  The same procedure works for
$ {\gerFphiMg M } $  too, so the proof is complete.   $ \square $
\enddemo

\vskip7pt

   {\bf 5.13  Matrix coefficients.}  \;  The result above can be
refined, extending embeddings to isomorphisms.
Let  $ \, \mu \in M_+ := M \cap P_+ \, $,  and let
$ V_{-\mu} $  be an irreducible  $ \uphizM{M'} $--module  of lowest weight
$ -\mu $  (recall that  $ \, \uqphiMg{M'} \cong
U_{q,0}^{\scriptscriptstyle M'} (\gerg) \, $  as algebras, hence their
representation theory is the same).  Let  $ v_{-\mu} \neq 0 $  be a lowest
weight vector of  $ V_{-\mu} $,  and let  $ \, \phi_{-\mu} \in
{V_{-\mu}}^{\! *} \, $  be the linear functional on  $ V_{-\mu} $  defined
by \  {\it (a)\/}~$ \phi_{-\mu} (v_{-\mu}) = 1 $  \  and  \
{\it (b)\/}~$ \phi_{-\mu} $   vanishes on the unique
$ \uphizM{M'} $--invariant  complement of  $ \kq.v_{-\mu} $  in
$ V_{-\mu} \, $;  let  $ \, \psi_{-\mu} := c_{\phi_{-\mu},v_{-\mu}} \, $  be
the corresponding matrix coefficient, i.~e.  $ \, \psi_{-\mu} \colon x
\mapsto \phi_{-\mu} (x.v_{-\mu}) \, $  for all  $ \, x  \in \uqphiMg{M'}
\, $.  The following refines Proposition 5.12, improving [DL], Theorem 4.6,
and [CV-2], Lemma 2.5:

\vskip7pt

\proclaim{Theorem 5.14}  Let  $ \, \rho := \sum_{i=1}^n \mu_i \, $
($ \{\mu_1, \dots, \mu_n\} $  being our fixed  $ \Z $--basis  of  $ M $,
cf.~\S 1.1).
                                           \hfill\break
   \indent    The algebra monomorphisms  $ \; \mu_{\scriptscriptstyle M}
\colon \, {\fqphiMg M } \llonghookrightarrow A_{\scriptscriptstyle M} \; $,
$ \; \mu_{\scriptscriptstyle M} \colon \, {\gerFphiMg M }
\llonghookrightarrow {\frak A}_\varphi^{\scriptscriptstyle M} \; $
                                    and\break
$ \; \mu_{\scriptscriptstyle M} \colon \, {\calFphiMg M }
\llonghookrightarrow {\Cal A}_\varphi^{\scriptscriptstyle M} \; $
respectively extend to algebra isomorphisms
   $$  \mu_{\scriptscriptstyle M} \colon {\fqphiMg M } \! \left[
\psi_{-\rho}^{-1} \right] {\buildrel \cong \over \longrightarrow} \,
A_\varphi^{\scriptscriptstyle M},  \quad  \mu_{\scriptscriptstyle M} \colon
{\gerFphiMg M } \! \left[ \psi_{-\rho}^{-1} \right] {\buildrel \cong \over
\longrightarrow} \, {\frak A}_\varphi^{\scriptscriptstyle M},  \quad
\mu_{\scriptscriptstyle M} \colon {\calFphiMg M } \! \left[ \psi_{-\rho}^{-1}
\right] {\buildrel \cong \over \longrightarrow} \,
{\Cal A}_\varphi^{\scriptscriptstyle M}  \, ;  \;  $$
moreover,  $ \, \mu_{\scriptscriptstyle M} \big( {\fqphiMg M } \big) \, $
and  $ A_\varphi^{\scriptscriptstyle M} $,  resp.~$ \,
\mu_{\scriptscriptstyle M} \big( {\gerFphiMg M } \big) \, $  and
$ {\frak A}_\varphi^{\scriptscriptstyle M} $,  resp.~$ \,
\mu_{\scriptscriptstyle M} \big( {\calFphiMg M } \big) \, $  and
$ {\Cal A}_\varphi^{\scriptscriptstyle M} $,  are dense in
$ \, j_{\scriptscriptstyle M} \left( {{\uqphiMg {M'} }}^* \right) \, $,
resp.~$ \, {\frak I}_\varphi^{\scriptscriptstyle M} \, $,
resp.~$ \, {\Cal I}_\varphi^{\scriptscriptstyle M} \, $.
\endproclaim

\demo{Proof}  It is proved in [DL], Theorem 4.6, that
$ \, \mu_{\scriptscriptstyle P} \colon \,
{\frak F}_0^{\scriptscriptstyle P} [G] \longhookrightarrow
{\frak A}_0^{\scriptscriptstyle P} \, $  extends to an isomorphism of
$ \kqqm $--algebras  $ \, \mu_{\scriptscriptstyle P} \colon \,
{\frak F}_0^{\scriptscriptstyle P} [G] \! \left[ \psi_{-\rho}^{-1} \right] \,
{\buildrel \cong \over \longrightarrow} \,
{\frak A}_0^{\scriptscriptstyle P} \, $:  in particular scalar extension
gives  $ \, \mu_{\scriptscriptstyle P} \colon \,
F_{q,0}^{\scriptscriptstyle P} [G] \! \left[ \psi_{-\rho}^{-1} \right] \,
{\buildrel \cong \over \longrightarrow} \, A_0^{\scriptscriptstyle P} \, $.
This is easily extended to general  $ \varphi $  and  $ M $.
                                                  \par
  Now, computations like in [DL] give also  $ \, \mu_{\scriptscriptstyle M}
(\psi_{-\mu}) = L^{\varphi,\otimes}_{-\mu} \, $  for all
$ \, \mu \in M_+ \, $;
therefore  $ \, \mu_{\scriptscriptstyle M} \left( \psi_{-\rho}^{-1} \right) =
L^{\varphi,\otimes}_{\rho} \, $.  Again from the proof in [DL] we get
$ \, F^\varphi_i \, L^{\varphi,\otimes}_{-\mu_i} $,
$ L^{\varphi,\otimes}_{-\mu_i} \, E^\varphi_i \in \mu_{\scriptscriptstyle M}
\left( {\fqphiMg M } \right) \, $,  hence
$ \, {\left( F^\varphi_i \right)}^{(f)} L^{\varphi,\otimes}_{- f \mu_i} $,
$ L^{\varphi,\otimes}_{- e \mu_i} {\left( E^\varphi_i \right)}^{(e)}
\in \mu_{\scriptscriptstyle M} \left( {\fqphiMg M }
                             \right) \, $  too; then Proposition
%
%\break
%
5.12{\it (b)} gives  $ \, {\left( F^\varphi_i \right)}^{(f)}
L^{\varphi,\otimes}_{- f \mu_i}, L^{\varphi,\otimes}_{- e \mu_i}
{\left( E^\varphi_i \right)}^{(e)} \in \mu_{\scriptstyle M}
\big( \calFphiMg{M} \big) \, $;  similarly we find that
$ \, L^{\varphi,\otimes}_{-\mu_i} = \mu_{\scriptscriptstyle M}(z_i)
\in \mu_{\scriptscriptstyle M} \big( \calFphiMg{M} \big) \, $,  with
$ \, z_i := \psi_{-\mu_i} \in \calFphiMg{M} \, $.  Then
  $$  L^{\varphi,\otimes}_{\mu_i} = \Bigg( \prod_{j=1}^n
{\phantom{\big\vert}}^{\!\!\widehat{i}} L^{\varphi,\otimes}_{-\mu_j} \Bigg)
\cdot L^{\varphi,\otimes}_{\rho} \in \mu_{\scriptscriptstyle M}
\Big( {\calFphiMg M } \! \left[ \psi_{-\rho}^{-1} \right] \Big)  $$
hence  $ \, {\left( F^\varphi_i \right)}^{(f)} \otimes 1 $,  $ 1 \otimes
{\left( E^\varphi_i \right)}^{(e)} \in \mu_{\scriptscriptstyle M}
\Big( {\calFphiMg M } \! \left[ \psi_{-\rho}^{-1}
\right] \Big) \, $;  moreover,  $ \, \left( L^{\varphi,\otimes}_{-\mu_i};
\, c \atop t \right) = \left( \mu_{\scriptscriptstyle M}(z_i); \, c \atop
t \right) = \mu_{\scriptscriptstyle M} \left( \left( z_i; c \atop t \right)
\right) \, $,  and  $ \left( z_i; \, c \atop t \right) \in {\calFphiMg M } $,
thus  $ \, \left( L^{\varphi,\otimes}_{-\mu_i}; \, c \atop t \right)
\in \mu_{\scriptscriptstyle M} \left( {\calFphiMg M }
\left[ \psi_{-\rho}^{-1} \right] \right) \, $.  Then Lemma 5.11
gives  $ \, \mu_{\scriptscriptstyle M} \left( {\calFphiMg M } \! \left[
\psi_{-\rho}^{-1} \right] \right) = {\Cal A}_\varphi^{\scriptscriptstyle M}
\, $.
The same can be done for the other integer form.
                                                   \par
   Now let  $ v_\tau $  be the image of   $ u_\tau $  (cf.~the proof
of Lemma 5.4) in the  $ \kq $--algebra  isomorphism  $ \, \theta \colon \,
{\uzM {M'} } @>{\cong}>> {\uzM {M'} } \, $  given by  $ \, L_\nu \mapsto
L_{-\nu} \, $  ($ \, \nu \in M' \, $):  then  $ \, \big\{\, v_\tau
\,\big\vert\,
\tau \in M'_+ \cong \N^n \,\big\} \, $  is a basis of  $ \uzM{M'} \, $;
a quick review of the proof of Lemma 5.4 shows that  $ v_\tau^* $
(with respect to  $ \, {\overline{im}}_{\scriptscriptstyle M} \colon \,
{\uphizM M } \longhookrightarrow {\uzM {M'} }^* \, $)  is a linear
combination of elements  $ \, L_{-\mu} \, $  ($ \mu \in M_+ $).  Then
$ \, j_{\scriptscriptstyle M} \big( L_{-\mu} \big) =
L^{\varphi,\otimes}_{-\mu} \, $  (cf.~(5.2)) and
$ \, L^{\varphi,\otimes}_{-\mu} \in \mu_{\scriptscriptstyle M}
\big( {\fqphiMg M } \big) \, $  imply  $ \, j_{\scriptscriptstyle M}
(v_\tau^*) \in \mu_{\scriptscriptstyle M} \big( {\fqphiMg M } \big) \, $,
for all  $ \, \tau \in M_+ \, $;  since  $ \, L_\mu \in {\uzM {M'} }^*  \, $
is a series of  $ v_\tau^* $  ($ \tau \in M'_+ $)  with
coefficients in  $ \kqqm $,  then  $ \, L^{\varphi,\otimes}_\mu =
j_{\scriptscriptstyle M} \big( L_\mu \big) \, $  lies in the (topological)
closure of  $ \mu_{\scriptscriptstyle M} \big( {\fqphiMg M } \big) $,  for
all  $ \, \mu \in M_+ \, $,  so the same is true for
$ \, L^{\varphi,\otimes}_\rho = \mu_{\scriptscriptstyle M}
\left( \psi_{-\rho}^{-1} \right) \, $:  this
proves the denseness claim for  $ {\fqphiMg M } $.  As
$ \, L^{\varphi,\otimes}_{-\mu} \in \mu_{\scriptscriptstyle M}
\big( {\gerFphiMg M } \big) $,  $ L^{\varphi,\otimes}_{-\mu} \in
\mu_{\scriptscriptstyle M} \big( {\calFphiMg M } \big) $,  \, this argument
works for integer forms too.   $ \square $
\enddemo

\vskip7pt

   {\bf 5.15  Gradings.}  \  Recall that  $ {\uphiMbp M } $  has a
$ Q_+ $--grading  $ \, {\uphiMbp M } = \oplus_{\alpha \in Q_+}
{\left( {\uMbp M } \right)}_{\! \alpha} \, $  given by decomposition in
direct sum of weight spaces for the adjoint action of  $ {\uphizM M } $;
also  $ {\uphiMbm M } $  has an analogous
$ Q_- $--grading.  These are gradings of Hopf algebras
(in the usual obvious sense), inherited by the integer forms, and
DRT pairings respect them, that is e.~g.
$ \, \pi \left({\left( {\uphiMbm M } \right)}_{\! \beta}, {\left(
{\uphiMbp {M'} } \right)}_{\! \gamma} \right) = 0  \, $
for all  $ \, \beta \in Q_- $,  $ \gamma \in Q_+ \, $
such that  $ \, \beta + \gamma \neq 0 \, $.
                                                    \par
   The gradings of quantum Borel subalgebras induce a
$ Q $--grading of the Hopf algebra  $ \, D_{\scriptscriptstyle M} :=
{\uphiMbp M } \otimes {\uphiMbm Q } \, $  (inherited by its quotient
Hopf algebra  $ {\uqphiMg M } $),  where the subspace
$ \, {\left( {\uphiMbp M } \right)}_{\! \beta} \otimes
{\left( {\uphiMbm Q } \right)}_{\! \gamma} \, $  has degree
$ \beta + \gamma \, $,  and also a  $ Q $--grading  of the
subalgebra  $ \, {\uphiMbp M } \otimes {\uphiMbm Q } \, $
of  $ {D_{\scriptscriptstyle M'}}^{\! *} \, $;  since
$ {D_{\scriptscriptstyle M'}}^{\! *} $  is a
completion (via formal series) of this subalgebra, it inherits on
its own sort of a "pseudograding", in the sense that every element of
$ {D_{\scriptscriptstyle M'}}^{\! *} $  is a
(possibly infinite) sum of terms each of whom has a well-defined degree:
namely, given  $ \, f \in {D_{\scriptscriptstyle M'}}^{\! *} \, $
with formal series expansion (cf.~Remark 5.5)  $ \, f =
\sum_{{\Cal F}^\varphi, \phi, {\Cal E}^\varphi} {\Cal F}^\varphi \cdot \phi
\cdot {\Cal E}^\varphi \, $  (where  $ \, \phi \in {\left( {\uphizM {M'} }
\otimes {\uphizM Q } \right)}^* \, $,  and  $ {\Cal F}^\varphi $'s  and
$ {\Cal E}^\varphi $'s  are PBW monomials),  we define the degrees
of its various summands as given by
  $$  deg \big( {\Cal F}^\varphi \cdot \phi \cdot {\Cal E}^\varphi \big)
:= deg \big( {\Cal F}^\varphi \big) + deg \big( {\Cal E}^\varphi \big)  $$
where  $ \, deg \left( \prod_{r=N}^1 {\left( F^\varphi_{\alpha^r}
\right)}^{f_r} \right) := - \sum_{r=1}^N
f_r \alpha^r \, $,  $ \, deg \left( \prod_{r=1}^N {\left(
E^\varphi_{\alpha^r} \right)}^{e_r} \right) :=
\sum_{r=1}^N e_r \alpha^r \, $  (this degree is again a weight for a
suitable action of  $ {\uphizM M } $  on  $ \, {\uphiMbm M }
\otimes {\uphiMbp P } \, $).  Now  $ \, {\uphiMbm M } \otimes
{\uphiMbp P } \, $  is dense in  $ {D_{\scriptscriptstyle M'}}^{\! *} $,
and the restriction of the pairing
$ \, {D_{\scriptscriptstyle M'}}^{\! *} \otimes D_{\scriptscriptstyle M'}
\rightarrow \kq \, $  to  $ \, \left( {\uphiMbm M } \otimes {\uphiMbp P }
\right) \otimes \left( {\uphiMbp {M'} } \otimes {\uphiMbm Q } \right)
\, $  is nothing but  $ \, \big( \pi_\varphi \otimes \overline{\pi_\varphi}
\big) \smallcirc \tau_{2,3} \, $  (with  $ \, \tau_{2,3} \colon \, x \otimes
y \otimes z \otimes w \mapsto x \otimes z \otimes y \otimes w \, $;)
therefore, since  $ \pi_\varphi $  and  $ \overline{\pi_\varphi} $  respect
the gradings, also the pairing   $ \, {D_{\scriptscriptstyle M'}}^{\! *}
\otimes D_{\scriptscriptstyle M'} \rightarrow \kq \, $  respects the
pseudogradings we are dealing with.
                                                        \par
   Finally, the pseudograding of
$ {D_{\scriptscriptstyle M'}}^{\! *} $  is compatible with the formal Hopf
structure.  For example, look at  $ S(x) $,  for homogeneous  $ \, x \in
{D_{\scriptscriptstyle M'}}^{\! *} \, $:  given homogeneous
$ \, y \in D_{\scriptscriptstyle M'} \, $,  we have
$ \, \big\langle S(x), y \big\rangle = \big\langle x, S(y) \big\rangle =
\big\langle x, y' \big\rangle \, $  where  $ \, y' := S(y) \, $  is
homogeneous on its own of degree  $ \, deg(y') = deg(y) \, $
(for the grading of  $ D_{\scriptscriptstyle M'} $  is
compatible with the Hopf structure); therefore
  $$  \big\langle S(x), y \big\rangle \neq 0 \Longrightarrow deg(y) =
deg(y') = deg(x) \Longrightarrow S(x) \in
{\left( {D_{\scriptscriptstyle M'}}^{\! *} \right)}_{deg(x)}  $$
that is  $ \, deg\big(S(x)\big) = deg(x) \, $,  q.e.d.

\vskip7pt

   {\bf 5.16  Umbral calculus.}  \;  In this section we provide concrete
information about the Hopf structure of our quantum formal groups.  This
will be especially important for defining integer forms and specializing
them at roots of 1.
                                                 \par
   The counit  $ \, \epsilon \colon \, {D_{\scriptscriptstyle
M'}}^{\! *} \rightarrow \kq \, $  is  $ \, \epsilon := 1^* \, $,
hence  $ \, \epsilon (x^*) := \big\langle x^*, 1 \big\rangle \, $  for all
$ \, x^* \in {D_{\scriptscriptstyle M'}}^{\! *} \, $;  thus
  $$  \epsilon \big( F^\varphi_i \otimes 1 \big) = 0 \, ,  \qquad  \epsilon
\big( L^{\varphi,\otimes}_\mu \big) = 1 \, ,  \qquad  \epsilon
\big(1 \otimes E_i \big) = 0 \; ;   \eqno (5.4)  $$
the elements above generate the algebra
$ \, j_{\scriptscriptstyle M} \big( {\uqMg {M'\!} }^* \big) \, $  (in
topological sense, cf.~Theorem 5.14), hence (5.4) uniquely determines
$ \, \epsilon \colon \, j_{\scriptscriptstyle M} \big( {\uqphiMg {M'\!} }^*
\big)  \longrightarrow \kq \, $.
                                                  \par
   The antipode of  $ {D_{\scriptscriptstyle M'}}^{\! *} \, $
is by definition the dual of the antipode of  $ D_{\scriptscriptstyle
M'} $,  hence it is characterized by  $ \, \big\langle S(x^*), x \big\rangle
= \big\langle x^*, S(x) \big\rangle \, $,  for all  $ \, x^* \in
{D_{\scriptscriptstyle M'}}^{\! *} $,  $ x \in D_{\scriptscriptstyle M'} $.
Now consider  $ \, {\left( F^\varphi_i \right)}^f \otimes 1 \in {\uphiMbm M }
\otimes {\uphiMbp P } \leq {D_{\scriptscriptstyle M'}}^{\! *} \, $,
$ \, f \in \N \, $:  it is homogeneous of degree  $ \, - f \alpha_i \, $,
whence  $ S \left( {\left( F^\varphi_i \right)}^f \otimes 1 \right) $
has the same degree.  Thus writing  $ S \left( {\left(
F^\varphi_i \right)}^f \otimes 1 \right) $  as a series
  $$  S \left( {\left( F^\varphi_i \right)}^f \otimes 1 \right) =
\sum_\sigma F_\sigma \cdot \varPhi_\sigma \cdot E_\sigma  $$
we have  $ \, deg \left( F^\varphi_\sigma \cdot \varPhi_\sigma \cdot
E^\varphi_\sigma \right) := deg(F^\varphi_\sigma) + deg(E^\varphi_\sigma)
= - f \alpha_i \, $.  Now, the pseudograding of
$ {D_{\scriptscriptstyle M'}}^{\! *} $  induces a pseudograding of
$ {\frak I}_\varphi^{\scriptscriptstyle M} $  too; hence, since
$ {\frak I}_\varphi^{\scriptscriptstyle M} $  is a formal Hopf subalgebra
of  $ {D_{\scriptscriptstyle M'}}^{\! *} $  (Proposition 5.9), we
can apply the same procedure and get
  $$  S \Big( {\big( \fbar^\varphi_i \big)}^f \otimes 1 \Big)
= \sum_\sigma \calF^\varphi_\sigma \cdot
\varphi_\sigma \cdot \calE^\varphi_\sigma   \eqno (5.5)  $$
where  $ \, \varphi_\sigma \in {\gerUphizM {{M'}} }^{\! *} \, $
and the  $ \calF^\varphi_\sigma $'s,  resp.~$ \calE^\varphi_\sigma $'s,
are PBW monomials of  $ \calUphim $,  resp.~$ \calUphip $,  such that
$ \, deg \! \left( \calF^\varphi_\sigma \right) + deg \! \left(
\calE^\varphi_\sigma \right) = - f \alpha_i \, $.  An entirely similar
argument yields
  $$  S \Big( {\big( F^\varphi_i \big)}^{(f)} \otimes 1 \Big)
= \sum_\sigma \gerF^\varphi_\sigma \cdot
\phi_\sigma \cdot \gerE^\varphi_\sigma   \eqno (5.6)  $$
where  $ \, \phi_\sigma \in {\calUzM {M'} }^* \, $  and the
$ \gerF^\varphi_\sigma $'s,  resp.~$ \gerE^\varphi_\sigma $,
are PBW monomials of  $ \gerUphim $,  resp.~$ \gerUphip $,   such that
$ \, deg \! \left( \gerF^\varphi_\sigma \right) + deg \! \left(
\gerE^\varphi_\sigma \right) = - f \alpha_i \, $.  Now remark that
$ {\frak I}_\varphi^{\scriptscriptstyle M} $  and
$ {\Cal I}_\varphi^{\scriptscriptstyle M} $  can be compared through the
natural embedding  $ \, {\Cal I}_\varphi^{\scriptscriptstyle M} \cong
{\gerUphiMg {M'} }^* \hookrightarrow {\calUphiMg {M'} }^* \cong
{\frak I}_\varphi^{\scriptscriptstyle M} \, $  (dual of
$ \, {\calUphiMg {M'} } \hookrightarrow {\gerUphiMg {M'} } \, $);  since
$ \, {\Big( \fbar^\varphi_{\alpha^h} \Big)}^f = \prod_{s=1}^f
\left( q_{\alpha^h}^s - q_{\alpha^h}^{-s} \right) \cdot {\left(
F^\varphi_{\alpha^h} \right)}^{(f)} \, $,  $ \, {\Big(
\ebar^\varphi_{\alpha^k} \Big)}^e = \prod_{s=1}^f \left( q_{\alpha^k}^s -
q_{\alpha^k}^{-s} \right) \cdot {\left( E^\varphi_{\alpha^k} \right)}^{(e)}
\, $,  comparing (5.5) and (5.6) we find
  $$  \gerF^\varphi_\sigma \cdot \phi_\sigma \cdot \gerE^\varphi_\sigma \in
\prod_{h,k=1}^n \prod_{r=1}^{f_h} \prod_{s=1}^{e_k} \left( q_{\alpha^h}^r -
q_{\alpha^h}^{-r} \right) \cdot \left( q_{\alpha^k}^s - q_{\alpha^k}^{-s}
\right) \cdot \prod_{u=1}^f {\left( q_i^u - q_i^{-u}
\right)}^{-1} \cdot {\Cal I}_\varphi^{\scriptscriptstyle M}  $$
for  $ \, \gerF^\varphi_\sigma = \prod_{h=N}^1 {\left( F^\varphi_{\alpha^h}
\right)}^{(f_h)} \, $,  $ \, \gerE^\varphi_\sigma = \prod_{k=1}^N
{\left( E^\varphi_{\alpha^k} \right)}^{(e_k)} \, $.  Therefore
  $$  S \left( {\left( F^\varphi_i \right)}^{(f)} \otimes 1 \right) =
\sum_\sigma \prod_{h,k=1}^n { \prod_{r=1}^{f_h} \prod_{s=1}^{e_k}
\left( q_{\alpha^h}^r - q_{\alpha^h}^{-r} \right) \cdot \left( q_{\alpha^k}^s
- q_{\alpha^k}^{-s} \right) \over \prod_{u=1}^f {\left( q_i^u - q_i^{-u}
\right)} } \cdot \gerF^\varphi_\sigma \cdot \phi'_\sigma \cdot
\gerE^\varphi_\sigma   \eqno (5.7)  $$
in particular from every coefficient in (5.7) we can pick out a factor
                                           of type\break
$ \, \prod_{h=1}^N \prod_{r=1}^{a_h} \prod_{s=1}^{b_h} \left( q^r - q^{-r}
\right) \cdot \left( q^s - q^{-s} \right) \, $  with  $ \, \sum_{h=1}^N
(a_h + b_h) = \sum_{h=1}^N (f_h + e_k) - f \, $;  then we can rearrange the
terms of the series (5.7) and write it again as
  $$  S \left( {\left( F^\varphi_i \right)}^{(f)} \otimes 1 \right) =
\sum_{n=0}^{+\infty} \, \sum_{\sum_h \! a_h + b_h = n} \, \prod_{h=1}^N
\prod_{r=1}^{a_h} \prod_{s=1}^{b_h} \left( q^r -
q^{-r} \right) \cdot \left( q^s - q^{-s} \right) \cdot X_n   \eqno (5.8)  $$
where  $ \, X_n \in \gerUphim \otimes {\calUphizM {M'} }^* \otimes \gerUphip
\, $.  Similarly occurs for the other generators of
$ \, {\Cal A}_\varphi^{\scriptscriptstyle M} \, $:  thus
                                                        \par
  {\it  For any root of unity  $ \varepsilon $,  the series
$ S \left( {\left( F^\varphi_i \right)}^{(f)} \otimes 1 \right) $,
$ S \left( \left( L^{\varphi,\otimes}_{\mu_i} ; \, c \atop t \right)
\right) $,  $ S \left( L^{\varphi,\otimes}_{-\mu_i} \right) $
                                  and\break
$ S \left( 1 \otimes {\left(  E^\varphi_i \right)}^{(e)} \right) $  are
{\sl finite}  sums modulo   $ (q - \varepsilon) $.}
                                                        \par
   In principle, one can compute all the terms of these series up to any
fixed order  $ n $;  actually, we need to know them only up to  $ \, n = 0
\, $.  For  $ S \big( F^\varphi_i \otimes 1 \big) $  the first
term (call it  $ \hbox{\bf F}_1 $),  with index  $ \, n = 0 \, $  in (5.8),
corresponds to the terms  $ \, \gerF^\varphi_\sigma \cdot \phi_\sigma \cdot
\gerE_\sigma \, $  in (5.6) such that  $ \, \sum_{s=1}^N (f_s + e_s) = 1
\, $;  but these must have degree  $ \, deg \! \left( \gerF^\varphi_\sigma
\right) + deg \! \left( \gerE^\varphi_\sigma \right) = -\alpha_i \, $  too,
whence it  $ \, \gerF^\varphi_\sigma = F_i \, $  and
$ \, \gerE^\varphi_\sigma = 1 \, $.  Now,  $ \hbox{\bf F}_1 $  takes
non-zero values only on the free  $ {\calUzM {M'} } $--module  with basis
$ \big\{ \ebar_i \big\} $,  call it  $ V_{1,i} \, $:  direct computation shows
that  $ \, \hbox{\bf F}_1 + q_i^{(\alpha_i \vert \alpha_i + \tau_i)} \cdot
F^\varphi_i L^{\varphi,\otimes}_{-\alpha_i} \, $  is zero in
$ {V_{1,i}}^{\!\!*} \, $,  therefore  $ \, \hbox{\bf F}_1 = - q_i^{(\alpha_i
\vert \alpha_i + \tau_i)} \cdot F^\varphi_i L^{\varphi,\otimes}_{-\alpha_i}
\, $,  whence
   $$  {\ } \hskip1,25cm  S \left( F^\varphi_i \otimes 1 \right) \equiv -
q^{- (\alpha_i \vert \alpha_i + \tau_i)}  \cdot F^\varphi_i L^{\varphi,
\otimes}_{-\alpha_i}   \eqno \mod \, \left( q-\qm \right)  $$
   \indent   Similar arguments give
  $$  \eqalignno{
   \hskip-1,5cm
   S \left( \left( L^{\varphi,\otimes}_{\mu_i} ; \, 0
\atop 1 \right) \right)  &  \equiv - L^{\varphi,\otimes}_{-\mu_i} \cdot
\left( L^{\varphi,\otimes}_{\mu_i} ; \, 0 \atop 1 \right)  &   \mod \,
\left( q-\qm \right)  \cr
   \hskip-1,5cm
   S \big( 1 \otimes E^\varphi_i \big) \equiv -  &
q^{+ (\alpha_i \vert \alpha_i - \tau_i)} \cdot L^{\varphi,\otimes}_{-\alpha_i}
E^\varphi_i  &  \mod \, \left( q - \qm \right)  \cr }  $$
   \indent   As for the coproduct  $ \, \Delta \colon
\, {D_{\scriptscriptstyle M'}}^{\! *} \rightarrow
{D_{\scriptscriptstyle M'}}^{\! *} \otimeshat {D_{\scriptscriptstyle
M'}}^{\! *} \, $,  it is the dual of the product of
$ D_{\scriptscriptstyle M'} $,  hence it is characterized by
$ \, \big\langle \Delta(x^*), y \otimes z \big\rangle = \big\langle x^*,
y \cdot z \big\rangle \, $.  Mimicking the procedure used for  $ S $,  we
find that  $ \, \Delta \left( {\left( F^\varphi_i \right)}^{(f)} \otimes 1
\right) \, $  is given by a series of type
  $$  \eqalign{
   \Delta \left( {\left( F^\varphi_i \right)}^{(f)} \otimes 1 \right) =
\sum_{n=0}^{+\infty} \; \sum_{\Sb  {\scriptscriptstyle \sum}_h (a_h + a'_h +
\\  + b_h + b'_h) = n \endSb} \;  &  \prod_{h=1}^N \prod_{r=1}^{a_h}
\prod_{s=1}^{b_h} \left( q^r - q^{-r} \right) \cdot \left( q^s - q^{-s}
\right) \cdot  \cr
   &  \cdot \prod_{r'=1}^{a'_h} \prod_{s'=1}^{b'_h} \left( q^{r'} - q^{-r'}
\right) \cdot \left( q^{s'} - q^{-s'} \right)  \cdot Y_n  \cr }
\eqno (5.9)  $$
in which  $ \, Y_n \in {\left( \gerUphim \otimes {\calUphizM
{M'} }^* \otimes \gerUphip \right)}^{\otimes 2} $.  Similar formulas exist
for all the generators of  $ \, {\Cal A}_\varphi^{\scriptscriptstyle M} \, $
(which are topological generators of  $ j_{\scriptscriptstyle M}
\Big( {\uqphiMg {M'} }^* \Big) $): in particular this implies
                                                      \par
  {\it  For any root of unity  $ \varepsilon $,  the series
$ \Delta \left( {\left( F^\varphi_i \right)}^{(f)} \otimes 1 \right) $,
$ \Delta \left( \left( L^{\varphi, \otimes}_{\mu_i} ; \, c \atop t \right)
\right) $,  $ \Delta \left( L^{\varphi,\otimes}_{-\mu_i} \right) $
and  $ \Delta \left( 1 \otimes {\left( E^\varphi_i \right)}^{(e)} \right) $
are  {\sl finite}  sums modulo  $ (q - \varepsilon) $.}
                                                      \par
   Direct computation gives us the following  {\sl congruences modulo}
$ {\left( q - \qm \right)}^2 $  (using notation  $ \, F^{\varphi,\otimes}_i, :=
F^\varphi_i \otimes 1 \, $,  $ \, L^{\varphi,\otimes}_\mu := L_{-(1+\varphi)\mu}
\otimes L_{(1-\varphi)(\mu)} \, $,  $ \, E^{\varphi,\otimes}_i := 1 \otimes
E^\varphi_i \, $,  and so on)
  $$  \displaylines{
   {\ }  \Delta \big( F^{\varphi,\otimes}_i \big) \equiv F^{\varphi,\otimes}_i
\otimes 1^\otimes + 1^\otimes \otimes F^{\varphi,\otimes}_i + (q_i - 1) \cdot
\bigg( {L^{\varphi,\otimes}_{\alpha_i} ; \, 0 \atop 1} \bigg) \otimes
F^{\varphi,\otimes}_i +   \hfill  \cr
   \hfill   + {\left( q_i - q_i^{\,-1} \right)}^{-1} \cdot \! \sum_{\alpha,
\beta \in R^+} C^{i,+}_{\alpha,\beta} \big( q_\alpha - q_\alpha^{\,-1} \big)
\big( q_\beta - q_\beta^{\,-1} \big) \cdot L^{\varphi,\otimes}_{\alpha_i}
E^{\varphi,\otimes}_\alpha \otimes F^{\varphi,\otimes}_\beta  \cr }  $$
  $$  \displaylines{
   \Delta \bigg( \bigg( {L^{\varphi,\otimes}_{\mu_i} ; \, 0 \atop 1} \bigg)
\bigg) \equiv \left( L^{\varphi,\otimes}_{\mu_i} ; \, 0 \atop 1 \right) \otimes
1^\otimes + 1^\otimes \otimes \left( L^{\varphi,\otimes}_{\mu_i} ; \, 0 \atop 1
\right) + (q_i - 1) \cdot \left( L^{\varphi,\otimes}_{\mu_i} ; \, 0 \atop 1
\right) \otimes \left( L^{\varphi,\otimes}_{\mu_i} ; \, 0 \atop 1 \right) +
\hfill  \cr
   + \; {(2)}_{\qm}^{\,2} {(d_i)}_q^{\,-1} \cdot \! \sum_{\gamma \in R^+}
(q - 1) \, {[d_\gamma]}_q {\big[ (\mu_i \vert \gamma) \big]}_q \cdot
L^{\varphi,\otimes}_{\mu_i} E^{\varphi,\otimes}_\gamma \otimes
F^{\varphi,\otimes}_\gamma L^{\varphi,\otimes}_{\mu_i}  \cr
   {\ }  \Delta \big( E^{\varphi,\otimes}_i \big) \equiv 1^\otimes \otimes
E^{\varphi,\otimes}_i + E^{\varphi,\otimes}_i \otimes 1^\otimes + (q_i - 1)
\cdot E^{\varphi,\otimes}_i \otimes \bigg( {L^{\varphi,\otimes}_{\alpha_i} ;
\, 0 \atop 1} \bigg) -   \hfill  \cr
   \hfill   - {\left( q_i - q_i^{\,-1} \right)}^{-1} \cdot \! \sum_{\alpha,
\beta \in R^+} C^{i,-}_{\alpha,\beta} \big( q_\alpha - q_\alpha^{\,-1} \big)
\big( q_\beta - q_\beta^{\,-1} \big) E^{\varphi,\otimes}_\alpha \otimes
F^{\varphi,\otimes}_\beta L^{\varphi,\otimes}_{\alpha_i}  \cr }  $$
where the  $ C^{i,\pm}_{\alpha,\beta} $'s  are given by the equations
$ \; \pi_i^- \big( [ F_\alpha, E_\beta ] \big) = C^{i,-}_{\alpha,\beta} \cdot
F_i \, $,  $ \, \pi_i^+ \big( [ F_\alpha, E_\beta ] \big) =
C^{i,+}_{\alpha,\beta} \cdot E_i \, $  ($ \, \pi_i^- : \uqMg{Q}
\twoheadrightarrow \kq \cdot F_i \, $  and  $ \, \pi_i^+ : \uqMg{Q}
\twoheadrightarrow \kq \cdot E_i \, $  being the canonical maps).

\vskip1,7truecm

 \centerline{ \bf  \S \; 6 \,  The quantum group  $ \uqphiMh{M} $ }

\vskip10pt

   {\bf 6.1 The quantum enveloping algebra  $ \, \uqphiMh{M} \, $.}  \  The
results of \S 5 can be given an axiomatic form: to this end, we introduce
a new object  $ {\uqphiMh M } $  which is with
respect to  $ \uhtau $  what  $ {\uqphiMg M } $  is for
$ \ugtau $.  Here  $ M $  is a fixed lattice as in \S 2.2.
                                                       \par
   We define  $ \hbox{\bf H}_\varphi^{\scriptscriptstyle M} $  to be the
associative  $ \kq $--algebra  with 1 with generators
  $$  F^\varphi_i \, ,  \  L^\varphi_\mu \, ,  \  E^\varphi_i  \quad
(\lambda \in M; \, i=\unon)  $$
and relations
  $$  \eqalign {
   L^\varphi_0 = 1 \, ,  \; \qquad  L^\varphi_\mu L^\varphi_\nu =
L^\varphi_{\mu + \nu}  &  \, ,  \; \qquad  E^\varphi_i F^\varphi_j -
F^\varphi_j E^\varphi_i = 0  \cr
   L^\varphi_\mu F^\varphi_j = q^{( \alpha_j | (1+\varphi)(\mu) )}
F^\varphi_j L^\varphi_\mu \, ,  &  \; \qquad  L^\varphi_\mu E^\varphi_j =
q^{( \alpha_j | (1-\varphi)(\mu) )} E^\varphi_j L^\varphi_\mu  \cr
   \sum_{k = 0}^{1-a_{ij}} (-1)^k q^{+c_{i j}^k} {\left[ {1-a_{ij}
\atop k } \right]}_{\! q_i}  &  \! {\left( E^\varphi_i \right)}^{1-\aij-k}E^\varphi_j {\left( E^\varphi_i \right)}^k = 0   \qquad \qquad  \forall \;
i \neq j  \cr
   \sum_{k = 0}^{1-a_{ij}} (-1)^k q^{-c_{i j}^k} {\left[ {1-a_{ij}
\atop k } \right]}_{\! q_i}  &  \! {\left( F^\varphi_i \right)}^{1-\aij-k}
F^\varphi_j {\left( F^\varphi_i \right)}^k = 0   \qquad \qquad  \forall \;
i \neq j  \cr }   \eqno (6.1)  $$
where  $ \, c_{i j}^k := - \big( k \alpha_i \,\big\vert\, \tau_j +
(1 - \aij - k) \, \tau_i \big) - \big( \alpha_j \,\big\vert\, (1 -
\aij - k) \, \tau_i \big) \, $  for all  $ i $,  $ j $,  $ k $.  We also
use notation  $ \, M^\varphi_i := L^\varphi_{\mu_i} \, $  ($ i= \unon $),
$ \{ \mu_1, \dots, \mu_n \} $  being a fixed  $ \Z $--basis  of  $ M $,
cf.~\S 1.1.
                                                       \par
   Now consider  $ F^\varphi_{\alpha^1}, \dots, F^\varphi_{\alpha^N} \, $
in  $ \, \uphim \; (\, \subseteq \hbox{\bf H}_\varphi^{\scriptscriptstyle M}
\,) \, $,  the elements  $ \, B^\varphi_{\scriptscriptstyle \eta, \tau, \phi}
:= \sum_{\tau' \preceq \tau} c'_{\tau, \tau'} \cdot L^\varphi_{\tau'} \, $
(cf.~\S 5.6)  in $ \, {\uphizM M } \; (\, \subseteq
\hbox{\bf H}_\varphi^{\scriptscriptstyle M} \,) \, $,  and
$ E^\varphi_{\alpha^1}, \dots, E^\varphi_{\alpha^N} $  in  $ \, \uphip \;
(\, \subseteq \hbox{\bf H}_\varphi^{\scriptscriptstyle M} \,) \, $.
                                                     \par
   {\it  We define  $ {\uqphiMh M } $  to be the completion of
$ \hbox{\bf H}_{\scriptscriptstyle M} $  by means of formal series, with
coefficients in  $ \kq $,  in the elements of the set}
  $$  {\Bbb B}^\varphi_{\scriptscriptstyle M} := \left\{\, \prod_{r=N}^1
{\left( F^\varphi_{\alpha^r} \right)}^{f_r} \cdot
B^\varphi_{\scriptscriptstyle \eta, \tau, \phi} \cdot
\prod_{r=1}^N {\left( E^\varphi_{\alpha^r} \right)}^{e_r}  \,\bigg\vert\,
\phi = {(f_r)}_r, \eta = {(e_r)}_r \in \N^N \, ;  \, \tau \in \N^n
\,\right\} \, .  $$
   \indent   Thus  $ {\uqphiMh M } $  is the completion of
$ \hbox{\bf H}_\varphi^{\scriptscriptstyle M} $  with respect to the topology
(of  $ \hbox{\bf H}_\varphi^{\scriptscriptstyle M} \, $)  for which a
fundamental system of neighborhoods of  $ 0 $  is the set of vector subspaces
of $ \hbox{\bf H}_\varphi^{\scriptscriptstyle M} $  which contain almost all
the elements of  $ {\Bbb B}^\varphi_{\scriptscriptstyle M} \, $,  and the set
$ \, {\Bbb B}^\varphi_{\scriptscriptstyle M} \, $  is a  {\it pseudobasis}
of  $ {\uqphiMh M } $.  Roughly speaking,  $ {\uqphiMh M } $  is an algebra
of (non-commutative) formal series with (6.1) as commutation rules.  Finally,
thanks to Lemma 5.3, we can identify  $ {\uqphiMh M } $  with the space of
formal series in the  $ F^\varphi_{\alpha^h} $'s,  $ E^\varphi_{\alpha^k} $'s
with coefficients in  $ \, {\uphizM {M'} }^* \, $.
                                                    \par
   From \S 5 we can explicitely realize  $ {\uqphiMh M } $  and endow it with
a Hopf structure: in fact, the definition of  $ \uqphiMh{M} $  is nothing
but a presentation of  $ {\uqphiMg{M'} }^* $,  as the following shows:

\vskip7pt

\proclaim{Theorem 6.2}  There exists an isomorphism of topological
$ \kq $--algebras
  $$  \nu^\varphi_{\scriptscriptstyle M} \colon \, {\uqphiMh M } \,
{\buildrel \cong \over \llongrightarrow}
\, j_{\scriptscriptstyle M} \left( {\uqphiMg {M'} }^* \right)  $$
given by: \  $ F^\varphi_i \mapsto F^\varphi_i \otimes 1 \, ,
\; L^\varphi_\mu \mapsto L^{\varphi,\otimes}_\mu \, ,  \; E^\varphi_i
\mapsto 1 \otimes E^\varphi_i \, . \; $  Then the pull-back of the formal
Hopf structure of  $ \, j_{\scriptscriptstyle M} \left( {\uqphiMg {M'} }^*
\right) $  uniquely defines a formal Hopf structure on
$ {\uqphiMh M } $,  so that  $ \nu^\varphi_{\scriptscriptstyle M} $  and
$ {j_{\scriptscriptstyle M}}^{\! -1} \smallcirc
\nu^\varphi_{\scriptscriptstyle M} \, $  are formal Hopf algebra
isomorphisms.
\endproclaim

\demo{Proof}  By construction  $ \, \hbox{\bf H}_\varphi^{\scriptscriptstyle M}
\cong \uphim \otimes {\uphizM M } \otimes \uphip \cong
A_\varphi^{\scriptscriptstyle M} \; (\, \subseteq j_{\scriptscriptstyle M}
\left( {\uqphiMg {M'} }^* \right) \,) \, $  as vector spaces;  now
$ \, F^\varphi_i \otimes 1, L^{\varphi,\otimes}_\mu, 1 \otimes E^\varphi_i
\in {\uphiMbm M } \otimes {\uphiMbp P } \, $,  hence comparing (6.1) and (2.1)
we see that formulas above gives a well-defined isomorphism of algebras
$ \, \nu^\varphi_{\scriptscriptstyle M} \colon \,
\hbox{\bf H}_\varphi^{\scriptscriptstyle M} @>{\cong}>>
A_\varphi^{\scriptscriptstyle M} \, $.  Moreover,
$ A_\varphi^{\scriptscriptstyle M} $  contains a pseudobasis  $ \hbox{\bf
B}^\varphi_{\scriptscriptstyle M} $  of  $ j_{\scriptscriptstyle M}
\left( {\uqphiMg {M'} }^* \right) $  (cf.~Lemma 5.4, Proposition 5.6, and
Remark 5.7)  such that  $ \, \nu^\varphi_{\scriptscriptstyle M}({\Bbb
B}^\varphi_{\scriptscriptstyle M}) = \hbox{\bf B}^\varphi_{\scriptscriptstyle
M} \, $,  hence  $ \nu^\varphi_{\scriptscriptstyle M} $  continuosly extends,
in a unique way, to an isomorphism of topological algebras
$ \, \nu_{\scriptscriptstyle M} \colon \, \uqphiMh{M} @>\cong>>
j_{\scriptscriptstyle M} \left( {\uqphiMg {M'} }^* \right) \, $,  q.e.d.
$ \square $
\enddemo

\vskip7pt

   {\bf Remark 6.3.} \; Notice that. setting
$ \, Y^\varphi_{\scriptscriptstyle \eta, \tau, \phi} :=
\calF^\varphi_\eta \cdot B^\varphi_{\scriptscriptstyle
\eta, \tau, \phi} \cdot \calE^\varphi_\phi \, $,
(notations of \S 5), Theorem 6.2 and definitions give
$ \, \nu^\varphi_{\scriptscriptstyle M} \left(
Y^\varphi_{\scriptscriptstyle \eta, \tau, \phi} \right)
= {X_{\scriptscriptstyle \eta, \tau, \phi}}^{\! *} \,
$  for all  $ \, \eta \in \N^N $,  $ \tau \in \N^n $,
$ \phi \in \N^N $  (cf.~\S 5.7).

\vskip7pt

\proclaim{Lemma 6.4}  The subset  $ \, \Omega_\varphi^{\scriptscriptstyle M}
:= \Big\{\, x = \sum_\sigma F^\varphi_\sigma \cdot \varPhi^\varphi_\sigma
\cdot E^\varphi_\sigma \in {\uqphiMh M } \,\Big\vert\,
\varPhi^\varphi_\sigma \in {\uphizM M }, \; \forall\, \sigma \,\Big\} \, $
(where  $ \, x = \sum_\sigma F^\varphi_\sigma \cdot \varPhi^\varphi_\sigma
\cdot E^\varphi_\sigma \, $  is the expansion of  $ \, x \in {\uqphiMh M }
\, $  as a series with coefficients in  $ \, {\uphizM {M'} }^* \, $)
is a formal Hopf subalgebra of  $ {\uqphiMh M } $.
\endproclaim

\demo{Proof}  It is clear that
$ \Omega_\varphi^{\scriptscriptstyle M} $  is a subalgebra of
$ {\uqphiMh M } $.  Now let  $ \, x = \sum_\tau {F^\varphi_\tau}'
\cdot {\varPhi^\varphi_\tau}' \cdot {E^\varphi_\tau}' \in
\Omega_\varphi^{\scriptscriptstyle M} \, $:  then
$ \, {\varPhi^\varphi_\tau}' = \sum_{\mu \in M } c_{\tau,\mu}
L^\varphi_\mu \, $  with  $ \, c_{\tau,\mu} \neq 0 \, $  for finitely
many  $ \mu $.
                                                   \par
   Let  $ \, S(x) = \sum_\sigma F^\varphi_\sigma \cdot \varPhi^\varphi_\sigma
\cdot E^\varphi_\sigma \, $:  for any fixed  $ \bar\sigma $,  we must
prove that  $ \, \varPhi^\varphi_{\bar\sigma} \in {\uphizM M } \;
\left(\, \subseteq {\uzM {M'} }^* \,\right) \, $,  so that
$ \, S \big( \Omega_\varphi^{\scriptscriptstyle M} \big) =
\Omega_\varphi^{\scriptscriptstyle M} \, $;  to this end, we
use the identification  $ \, \uqphiMh{M} \cong {\uqphiMg {M'} }^* \, $
(cf.~Theorem 6.2).  For (2.3) there exist two PBW monomials
$ {\Cal E}_{\bar\sigma} $  and  $ {\Cal F}_{\bar\sigma} $  such that
  $$  \big\langle S(x), {\Cal E}_{\bar\sigma} \cdot y \cdot
{\Cal F}_{\bar\sigma} \big\rangle = \big\langle F^\varphi_{\bar\sigma}
\cdot \varPhi^\varphi_{\bar\sigma} \cdot E^\varphi_{\bar\sigma}, {\Cal
E}_{\bar\sigma} \cdot y \cdot {\Cal F}_{\bar\sigma} \big\rangle = \big\langle
F_{\bar\sigma}, {\Cal E}_{\bar\sigma} \big\rangle \cdot \big\langle
E_{\bar\sigma}, {\Cal F}_{\bar\sigma} \big\rangle \cdot
\varPhi^\varphi_{\bar\sigma} \big( y \cdot L_\alpha \big)  $$
for all  $ \, y \in {\uzM {M'} } \, $,  with  $ \, \alpha := s \big(
F_{\bar\sigma} \big) + s \big( E_{\bar\sigma} \big) \, $  and
$ \, c_{\bar\sigma} := \big\langle F_{\bar\sigma}, {\Cal E}_{\bar\sigma}
\big\rangle \cdot \big\langle E_{\bar\sigma}, {\Cal F}_{\bar\sigma}
\big\rangle \neq 0 \, $:  in other words,  $ \, \varPhi^\varphi_{\bar\sigma} =
{c_{\bar\sigma}}^{\! -1} \cdot \Big( \big( L_{-\alpha} \cdot
{\Cal F}_{\bar\sigma} \big) \triangleright S(x) \triangleleft
{\Cal E}_{\bar\sigma} \Big){\Big\vert}_{{\uzM {M'} }} \, $  (where
$ \triangleleft $  and  $ \triangleright $  denote standard left and right
action, cf.~[DL], \S 1.4),  hence we have to study  $ \, \big\langle S(x),
{\Cal E}_{\bar\sigma} \cdot y \cdot L_{-\alpha} {\Cal F}_{\bar\sigma}
\big\rangle \, $  as a function of  $ \, y \in {\uphizM {M'} } \, $;  by
linearity we can assume  $ \, y = L_\nu \, $,  $ \, \nu \in M' \, $.  By
definition,  $ \, \big\langle S(x), {\Cal E}_{\bar\sigma} \cdot y \cdot
L_{-\alpha} {\Cal F}_{\bar\sigma} \big\rangle = \big\langle x,
S({\Cal E}_{\bar\sigma} \cdot y L_{-\alpha} \cdot {\Cal F}_{\bar\sigma})
\big\rangle \, $;  in order to compute the latter we have to "straighten"
$ S \big( {\Cal E}_{\bar\sigma} \cdot y \cdot L_{-\alpha} {\Cal
F}_{\bar\sigma} \big) $,  i.~e.~to express it in terms of a PBW basis
of  $ \, U_{\scriptscriptstyle -} \otimes {\uzM {M'} } \otimes
U_{\scriptscriptstyle +} \, $.  Since  $ \, S \big( {\Cal E}_{\bar\sigma}
\cdot y \cdot L_{-\alpha} {\Cal F}_{\bar\sigma} \big) =
S \big( L_{-\alpha} {\Cal F}_{\bar\sigma} \big) \cdot S(y)
\cdot S \big( {\Cal E}_{\bar\sigma} \big) \, $,  let us
consider the various factors.
                                                  \par
   First,  $ S \big( L_{-\alpha} {\Cal F}_{\bar\sigma} \big) \in
{\uMbm {M'} } \, $,  and  $ S \big( L_{-\alpha} {\Cal F}_{\bar\sigma}
\big) $  {\it does not depend on}  $ y $.  Second,  $ S(y) = S \big( L_\nu \big)
= L_{-\nu} \, $.  Third,  $ S \big( {\Cal E}_{\bar\sigma} \big) \in
{\uMbp {M'} } \, $,  and  $ S \big( {\Cal E}_{\bar\sigma} \big) $
{\it does not depend on}  $ y $.
                                                  \par
   Now we straighten the product.  Commuting  $ S \big( L_{-\alpha} {\Cal
F}_{\bar\sigma} \big) $  and  $ \, S(y) = L_{-\nu} \, $  produces a coefficient
$ \, q^{-(\nu \vert \beta_{\bar\sigma})} = {\big\langle
L^\varphi_{-\beta_{\bar\sigma}}, L_\nu \big\rangle}_{\! \overline{\pi}}
\, $,  where  $ \, \beta_{\bar\sigma} \in Q_- \, $  is the weight of
$ S \big( {\Cal F}_{\bar\sigma} \big) $.  Straightening the product
$ \, S \big( L_{-\alpha} {\Cal F}_{\bar\sigma} \big) \cdot S \big( {\Cal
E}_{\bar\sigma} \big) \, $  produces a sum  $ \, \sum_k x_k \, $  of terms
which  {\it do not depend on}  $ y \, $.  Straightening the product
$ \, S \big( y \big) = L_{-\nu} \cdot \sum_k x_k \, $  produces for each term
$ x_k $  a coefficient  $ \, q^{-(\nu \vert \gamma_{ \bar\sigma, k})} =
{\big\langle L^\varphi_{-\gamma_{ \bar\sigma, k}}, L_\nu \big\rangle}_{\!
\overline{\pi}} \, $,  where  $ \, \gamma_{\bar\sigma, k} \in Q_+ \, $  is the
weight of the "positive" part  $ x_k^+ $  of  $ x_k $  (with respect to the
triangular decomposition).
                                                  \par
   Therefore  $ \, \big\langle x, S \big( {\Cal E}_{\bar\sigma} \cdot y
\cdot L_{-\alpha} {\Cal F}_{\bar\sigma} \big) \big\rangle \, $  depends on
$ y $  according to the functions  $ L^\varphi_{-\beta_{ \bar\sigma}} $,
$ L^\varphi_{-\gamma_{ \bar\sigma, k}} $,  and  $ \, {\varPhi^\varphi_\tau}'
\smallcirc S \, $:  to be precise,  $ \, \varPhi^\varphi_{\bar\sigma} = \big(
L_{-\alpha} {\Cal F}_{\bar\sigma} \triangleright S(x) \triangleleft {\Cal
E}_{\bar\sigma} \big) {\Big\vert}_{{\uzM {M'} }} \, $  is a linear combination
of functions of type  $ \, L^\varphi_{-\beta_{\bar\sigma}} \cdot \big(
{\varPhi^\varphi_\tau}' \smallcirc S \big) \cdot L^\varphi_{-\gamma_{\bar\sigma,
k}} = \sum_{\mu \in M} c_{\tau,\mu} L^\varphi_{- \mu - \beta_{\bar\sigma} -
\gamma_{\bar\sigma, k}} \, $,  so  $ \, \varPhi^\varphi_{\bar\sigma} \in
{\uphizM M } \, $,  q.e.d.  An entirely analogous procedure   --- slightly
simpler indeed ---   works for comultiplication, thus proving that  $ \, \Delta
\big( \Omega_\varphi^{\scriptscriptstyle M} \big) \subseteq
\Omega_\varphi^{\scriptscriptstyle M} \otimeshat
\Omega_\varphi^{\scriptscriptstyle M} \, $.  The thesis follows.   $ \square $
\enddemo

\vskip7pt

   Now we introduce integer forms of  $ {\uqphiMh M } $  and prove their
first properties.  We freely use the term  {\sl pseudobasis}  to mean a
topological basis of a topological module, so that any element in the module
has a unique expansion as a series in the elements of the basis.

\vskip7pt

\proclaim{Definition 6.5}  We define
$ {\Cal H}_\varphi^{\scriptscriptstyle M} $  to
be the  $ \kqqm $--subalgebra  of  $ {\uqphiMh M } $
                                       generated by\break
$ \, \big\{\, \fbar^\varphi_{\alpha^r}, \, L^\varphi_\mu,
\, \ebar^\varphi_{\alpha^r} \,\big\vert\, r= 1,
\dots, N \,; \, \mu \in M \,\big\} \, $,  and  $ \, { \calUphiMh M } \, $
to be its closure in  $ {\uqphiMh M } $.
\endproclaim

\vskip7pt

\proclaim{Theorem 6.6}  $ { \calUphiMh M } $  is a
$ \kqqm $--integer  form (in topological sense) of
$ {\uqphiMh M } $,  as a formal Hopf algebra, with
$ \kqqm $--pseudobasis
  $$  {\widetilde{\Bbb B}}^\varphi_{\scriptscriptstyle M} :=
\Big\{\, Y^\varphi_{\scriptscriptstyle \eta, \tau, \phi} \,\Big\vert\,
\tau \in \N^n; \eta, \phi \in \N^N \,\Big\} =  \Big\{\, \calF^\varphi_\eta
\cdot B^\varphi_{\scriptscriptstyle \eta, \tau, \phi} \cdot
\calE^\varphi_\phi \,\Big\vert\, \tau \in \N^n; \eta, \phi \in \N^N \,\Big\} \;
;   \eqno (6.2)  $$
in particular  $ \, \nu^\varphi_{\scriptscriptstyle M}
\big( {\calUphiMh M } \big) = j_{\scriptscriptstyle M}
\left( {\gerUphiMg {M'\!} }^* \right) =:
{\frak I}_\varphi^{\scriptscriptstyle M} \, $.
\endproclaim

\demo{Proof}  By construction
$ \, {\widetilde{\Bbb B}}^\varphi_{\scriptscriptstyle M}  \subseteq
{ \calUphiMh M } \, $,  so the claim follows
from \S 6.1 or Remark 6.3.   $ \square $
\enddemo

\vskip7pt

   Let  $ \, \widehat{\Omega}_\varphi^{\scriptscriptstyle M} :=
\Omega_\varphi^{\scriptscriptstyle M} \cap {\left(
\nu^\varphi_{\scriptscriptstyle M} \right)}^{-1}
({\Cal I}_\varphi^{\scriptscriptstyle M}) \, $;  notice that
(cf.~Proposition 5.9{\it (b)\/})
  $$  \widehat{\Omega}_\varphi^{\scriptscriptstyle M} = \bigg\{\, x =
\sum_\sigma \gerF^\varphi_\sigma \cdot \phi^\varphi_\sigma \cdot
\gerE^\varphi_\sigma \in {\uqphiMh M } \,\bigg\vert\, \gerF^\varphi_\sigma
\in \gerUphim, \phi^\varphi_\sigma \in {\gerUphizM M }, \gerE^\varphi_\sigma
\in \gerUphip, \; \forall\, \sigma \,\bigg\} \, .  $$

\vskip7pt

\proclaim{Definition 6.7}  We call
$ {\frak H}_\varphi^{\scriptscriptstyle M} $  the  $ \kqqm $--subalgebra  of
$ {\uqphiMh M } $
             generated by\break
$ \, \Big\{\, {\big( F^\varphi_i \big)}^{(f)}, \, \left( M^\varphi_i ;
c \atop t \right), \, {\big( M^\varphi_i \big)}^{-1}, \,  {\big( E^\varphi_i
\big)}^{(e)} \, \Big\vert \, f, c, t, e \in \N; i=\unon \,\Big\} \, $,  and
$ \, { \gerUphiMh M } \, $  the set
  $$  \Bigg\{ x \in
\widehat{\Omega}_\varphi^{\scriptscriptstyle M} \,\Bigg\vert\, x =
\sum_{n=0}^{+\infty} x_n, \; \, x_n \in \! \sum_{{\scriptscriptstyle \sum}_h
(\! a_h + b_h \!) = n} \prod_{h=1}^N \prod_{r,s=1}^{a_h,b_h} \! \left( q^r -
q^{-r} \right) \! \cdot \! \left( q^s - q^{-s} \right) \cdot
{\frak H}_\varphi^{\scriptscriptstyle M} \; \, \forall \, n \Bigg\}
\eqno (6.3)  $$
\endproclaim

\vskip7pt

\proclaim{Theorem 6.8}  $ { \gerUphiMh M } $  is a
$ \kqqm $--integer  form of  $ {\uqphiMh M } $  and
$ \Omega_\varphi^{\scriptscriptstyle M} $.
\endproclaim

\demo{Proof}  By construction  $ { \gerUphiMh M } $  is a
$ \kqqm $--subalgebra  of  $ {\uqphiMh M } $  and
$ \Omega_\varphi^{\scriptscriptstyle M} $;  moreover Theorem 6.2 and
Proposition 5.9{\it (b)}  ensure that
$ \widehat{\Omega}_\varphi^{\scriptscriptstyle M} $  is a
$ \kqqm $--integer  form (in topological sense) of
$ \Omega_\varphi^{\scriptscriptstyle M} $  (as an algebra), hence also
$ { \gerUphiMh M } $ is.  Proposition 5.9{\it (b)}  and Lemma 6.4 imply
that  $ \widehat{\Omega}_\varphi^{\scriptscriptstyle M} $  is a formal
Hopf subalgebra of  $ \Omega_\varphi^{\scriptscriptstyle M} $.  Finally
the analysis in \S 5.16 (especially (5.8) and (5.9)) via
$ {\nu_{\scriptscriptstyle M}}^{\! -1} $  gives
$ \, S\big({\gerUphiMh M }\big) = {\gerUphiMh M } \, $  and
$ \, \Delta\big({ \gerUphiMh M }\big) \subseteq {\gerUphiMh M } \otimeshat
{\gerUphiMh M } \, $.   $ \square $
\enddemo

\vskip7pt

   {\bf 6.9  Presentation of  $ {\gerUphiMh M } $.}  \;  By the similar result
available for  $ \, {\gerUphiMg M } \cong \gerU_{q,0}^{\scriptscriptstyle M}
(\gerg) \, $  (cf.~[DL], \S 3.4) we get a presentation of  $ {\gerUphiMh M } $
by (topological) generators and relations.  The algebra  $ {\frak
H}_\varphi^{\scriptscriptstyle M} $  of \S 6.7 is the associative  $ \kqqm
$--algebra  with 1 with generators
  $$  M^\varphi_i \, ,  \; {\big( M^\varphi_i \big)}^{-1} \, ,  \;
\left( M^\varphi_i; c \atop t \right) \, ,  \; {\big( E^\varphi_i
\big)}^{(r)} \,  \; {\big( F^\varphi_i \big)}^{(s)}  $$
($ i=\unon $;  $ c \in \Z $,  $ t, r, s \in \N $;  here we set
$ \, M^\varphi_i := L^\varphi_{\mu_i} \, $),  and relations
  $$  \displaylines{
   M^\varphi_i {\left( M^\varphi_i \right)}^{-1} = 1 =
{\left( M^\varphi_i \right)}^{-1} M^\varphi_i \, , \; \qquad
{\left( M^\varphi_i \right)}^{\pm 1} {\left( M^\varphi_j \right)}^{\pm 1}
= {\left( M^\varphi_j \right)}^{\pm 1} {\left( M^\varphi_i
\right)}^{\pm 1}  \cr
   {\left( M^\varphi_i \right)}^{\pm 1} \left( M^\varphi_j; c  \atop t
\right) = \left( M^\varphi_j; c \atop t \right) {\left( M^\varphi_i
\right)}^{\pm 1} \, ,  \; \left( M^\varphi_i; c \atop 0 \right) = 0 \, ,
\;  (q_i - 1) \left( M^\varphi_i; 0 \atop 1 \right) = M^\varphi_i - 1  \cr
   \left( M^\varphi_i; c \atop t \right) \left( M^\varphi_i; c-t \atop
s \right) = {\left( t+s \atop t \right)}_{\! q} \left( M^\varphi_i;
c \atop t+s \right) \, ,   \; \qquad \forall \, t, s  \cr
   \left( M^\varphi_i; c+1 \atop t \right) - q^t \left( M^\varphi_i;
c \atop t \right) = \left( M^\varphi_i; c \atop t-1 \right) \, ,
\; \qquad \forall \, t \geq 1  \cr
   \left( M^\varphi_i; c \atop t \right) = \sum_{p \geq 0}^{p \leq c, t}
q^{(c-p)(t-p)} \left( c \atop p \right)_{\! q} \left( M^\varphi_i; 0 \atop
t-1 \right) \, ,   \; \qquad \forall \, c \geq 0  \cr
   \left( M^\varphi_i; -c \atop t \right) = \sum_{p=0}^t {(-1)}^p
q^{-t(c+p) + p(p+1)/2} {\left( p+c-1 \atop p \right)}_{\! q}
\left( M^\varphi_i; 0 \atop t-p \right) \, ,   \; \qquad \forall \,
c \geq 1  \cr
   \left( M^\varphi_i; c+1 \atop t \right) - \left( M^\varphi_i; c \atop
t \right) = q^{c-t+1} M^\varphi_i \left( M^\varphi_i; c \atop t-1 \right)
\, , \; \qquad \forall \, t \geq 1  \cr
   M^\varphi_i {\left( E^\varphi_j \right)}^{(p)} = q^{p \, (\alpha_j
\vert (1+\varphi)(\mu_i))} {\left( E^\varphi_j \right)}^{(p)} M^\varphi_i
\, , \; \;  M^\varphi_i {\left( F^\varphi_j \right)}^{(p)} = q^{p \,
(\alpha_j \vert (1-\varphi)(\mu_i))} {\left( F^\varphi_j \right)}^{(p)}
M^\varphi_i  \cr
   \left( M^\varphi_i; c \atop t \right) {\left( E^\varphi_j
\right)}^{(p)} = {\left( E^\varphi_j \right)}^{(p)} \left( M^\varphi_i;
c + p \, (\alpha_j \vert (1+\varphi)(\mu_i)) \atop t \right)  \cr
   \left( M^\varphi_i; c \atop t \right) {\left( F^\varphi_j
\right)}^{(p)} = {\left( F^\varphi_j \right)}^{(p)} \left( M^\varphi_i; c + p \,
(\alpha_j \vert (1-\varphi)(\mu_i)) \atop t \right)  \cr }  $$
  $$  \displaylines{
   {\left( E^\varphi_i \right)}^{(r)} {\left( E^\varphi_i \right)}^{(s)} =
{\left[ r+s \atop r \right]}_{q_i} {\left( E^\varphi_i \right)}^{(r+s)}
\, , \qquad  {\left( F^\varphi_i \right)}^{(r)} {\left( F^\varphi_i
\right)}^{(s)} = {\left[r+s \atop r \right]}_{q_i} {\left( F^\varphi_i
\right)}^{(r+s)}  \cr
   \sum_{r+s=1-\aij} {(-1)}^s {\left( E^\varphi_i \right)}^{(r)}
E^\varphi_j {\left( E^\varphi_i \right)}^{(s)} = 0 \, ,  \;
\sum_{r+s=1-\aij} {(-1)}^s {\left( F^\varphi_i \right)}^{(r)}
F^\varphi_j {\left( F^\varphi_i \right)}^{(s)} = 0 \, ,  \quad
\forall \; i \neq j  \cr
   {\left( E^\varphi_i \right)}^{(0)} = 1 \, ,  \qquad
{\left( E^\varphi_i \right)}^{(r)} {\left( F^\varphi_j \right)}^{(s)} =
{\left( F^\varphi_j \right)}^{(s)} {\left( E^\varphi_i \right)}^{(r)} \, ,
\qquad  {\left( F^\varphi_i \right)}^{(0)} = 1  \cr }  $$
   \indent   Then  $ { \gerUphiMh M } $  is the
completion of  $ {\frak H}_\varphi^{\scriptscriptstyle M} $  obtained by
taking formal series in the PBW monomials of  $ \gerUphim $  and
$ \gerUphip $,  with coefficients in  $ {\gerUphizM M } $,  which satisfy
the condition in (6.3).  Finally, formulas in \S 5.16 yield   --- via
$ \nu^\varphi_{\scriptscriptstyle M} $  ---   the following (where
$ \, K^\varphi_i := L^\varphi_{\alpha_i} \, $):
  $$  \displaylines{
   {\ }  \Delta \big( F^\varphi_i \big) \equiv F^\varphi_i \otimes
1 + 1 \otimes F^\varphi_i + (q_i - 1) \cdot \bigg( {K^\varphi_i ; \, 0 \atop 1}
\bigg) \otimes F^\varphi_i +   \hfill  \cr
   \hfill   + {\left( q_i - q_i^{\,-1} \right)}^{-1} \cdot \! \sum_{\alpha,
\beta \in R^+} C^{i,+}_{\alpha,\beta} \big( q_\alpha - q_\alpha^{\,-1} \big)
\big( q_\beta - q_\beta^{\,-1} \big) K^\varphi_i E^\varphi_\alpha \otimes
F^\varphi_\beta   \hfill   \mod\, {\left( q - \qm \right)}^2  \cr    \Delta \bigg( \left( M^\varphi_i ; \, 0 \atop 1 \right) \bigg) \equiv
\left( M^\varphi_i ; \, 0 \atop 1 \right) \otimes 1 + 1 \otimes
\left( M^\varphi_i ; \, 0 \atop 1 \right) + (q_i - 1) \cdot
\left( M^\varphi_i ; \, 0 \atop 1 \right) \otimes
\left( M^\varphi_i ; \, 0 \atop 1 \right) +   \hfill  \cr
   \hfill   + \; {(2)}_{\qm}^{\,2} {(d_i)}_q^{\,-1} \cdot \! \sum_{\gamma \in
R^+} (q - 1) \, {[d_\gamma]}_q {\big[ (\mu_i \vert \gamma) \big]}_q \cdot
M^\varphi_i E^\varphi_\gamma \otimes F^\varphi_\gamma M^\varphi_i   \hfill
\mod\, {\left( q - \qm \right)}^2  \cr
     {\ }  \Delta \big( E^\varphi_i \big) \equiv 1 \otimes E^\varphi_i +
E^\varphi_i \otimes 1 + (q_i - 1) \cdot E^\varphi_i \otimes \left( K^\varphi_i ;
\, 0 \atop 1 \right) -   \hfill  \cr
   \hfill   - {\left( q_i - q_i^{\,-1} \right)}^{-1} \cdot \! \sum_{\alpha,
\beta \in R^+} C^{i,-}_{\alpha,\beta} \big( q_\alpha - q_\alpha^{\,-1} \big)
\big( q_\beta - q_\beta^{\,-1} \big) E^\varphi_\alpha \otimes F^\varphi_\beta
K^\varphi_i   \hfill   \mod\, {\left( q - \qm \right)}^2  \cr
   \hfill   S \left( F^\varphi_i \right) \equiv - q_i^{-2} \cdot F^\varphi_i
{\left( K^\varphi_i \right)}^{-1} \; ,  \qquad  S \left( E^\varphi_i \right)
\equiv - q_i^{+2} \cdot {\left( K^\varphi_i \right)}^{-1} E^\varphi_i
\hfill   \mod\, \left( q-\qm \right)  \cr
   \hfill   S \left( \left( M^\varphi_i ; \, 0 \atop 1 \right) \right) \equiv -
{\left( M^\varphi_i \right)}^{-1} \cdot \left( M^\varphi_i ; \, 0 \atop
1 \right)   \hfill  \mod\, \left( q-\qm \right)  \cr
   \epsilon \left( F^\varphi_i \right) = 0 \; ,  \qquad
\epsilon \left( \left( M^\varphi_i ; 0 \atop 1 \right) \right) = 0 \; ,  \qquad
\epsilon \left( E^\varphi_i \right) = 0  \; .  \hskip30pt  \cr }  $$

\vskip7pt

\proclaim{Definition 6.10}  We call  $ \, \xi^\varphi_{\scriptscriptstyle M}
\, $  the embedding of formal Hopf algebras
  $$  \xi^\varphi_{\scriptscriptstyle M} :=
{\left( \nu^\varphi_{\scriptscriptstyle M} \right)}^{-1}
\smallcirc \mu^\varphi_{\scriptscriptstyle M} \colon \, {\fqphiMg M }
\llonghookrightarrow {\uqphiMh M } \;  $$
\endproclaim

\vskip7pt

\proclaim{Theorem 6.11}  The embedding  $ \, \xi^\varphi_{\scriptscriptstyle
M} \colon \, {\fqphiMg M } \hookrightarrow {\uqphiMh M } \, $  induces
algebra monomorphisms  $ \, \xi^\varphi_{\scriptscriptstyle M} \colon \,
{\fqphiMg M } \longhookrightarrow \hbox{\bf H}_\varphi^{\scriptscriptstyle M}
\, $,  $ \, \xi^\varphi_{\scriptscriptstyle M} \colon \, {\gerFphiMg M }
\longhookrightarrow {\Cal H}_\varphi^{\scriptscriptstyle M} \, $,
$ \, \xi^\varphi_{\scriptscriptstyle M} \colon \, {\calFphiMg M }
\longhookrightarrow {\frak H}_\varphi^{\scriptscriptstyle M} \, $
and algebra isomorphisms  $ \, \xi^\varphi_{\scriptscriptstyle M} \colon \,
{\fqphiMg M } \! \left[ \psi_{-\rho}^{-1} \right] \, {\buildrel \cong \over
\longrightarrow} \, \hbox{\bf H}_\varphi^{\scriptscriptstyle M} \, $,
$ \, \xi^\varphi_{\scriptscriptstyle M} \colon \, {\gerFphiMg M } \!
\left[ \psi_{-\rho}^{-1} \right] \, {\buildrel \cong \over \longrightarrow}
\, {\Cal H}_\varphi^{\scriptscriptstyle M} \, $,
$ \, \xi^\varphi_{\scriptscriptstyle M} \colon \, {\calFphiMg M } \!
\left[ \psi_{-\rho}^{-1} \right] \, {\buildrel \cong \over \longrightarrow}
\, {\frak H}_\varphi^{\scriptscriptstyle M} \, $  whose images are dense
respectively in  $ {\uqphiMh M } $,  in  $ {\calUphiMh M } $,  in
$ {\gerUphiMh M } $.   $ \square $
\endproclaim

\vskip7pt

   {\bf  6.12  Quantum Poisson pairing.}  \;  In this section we define
perfect Hopf pairings \  $ \uqphiMh{M} \otimes \uqphiMg{M'} \longrightarrow \kq
$  \ which provide quantizations of the Hopf pairings  $ \, \fgtau \otimes
\ugtau \rightarrow k \, $  (or  $ \, \finfgtau \otimes \ugtau \rightarrow k \,
$)  and  $ \, \uhtau \otimes \fhtau \rightarrow k \, $  and of the Lie bialgebra
pairing  $ \, \gerh^\tau \otimes \gerg^\tau \rightarrow k \, $:  therefore we
call them "(multiparameter) quantum Poisson pairings";
moreover they also provide new interesting pairings between
function algebras.
 \eject
   Since  $ \, {j_{\scriptscriptstyle M}}^{\! -1} \smallcirc
\nu^\varphi_{\scriptscriptstyle M} \colon \, {\uqphiMh M } {\buildrel \cong
\over \longrightarrow} {\uqphiMg {M'} }^* \, $,  evaluation gives a perfect
Hopf pairing
  $$  \pi_{q,\varphi}^{\scriptscriptstyle M} \colon \, {\uqphiMh M } \otimes
{\uqphiMg {M'} } \llongrightarrow \kq  $$
defined by \  $ \pi_{q,\varphi}^{\scriptscriptstyle M} (h,g) := \left\langle
{j_{\scriptscriptstyle M}}^{\! -1} \big( \nu^\varphi_{\scriptscriptstyle M}
(h) \big), g \right\rangle $  \  for all  $ \, h \in
{\uqphiMh M } $,  $ g \in {\uqphiMg {M'} } \, $.
                                                      \par
   We call  $ \pi_{q,\varphi}^{\scriptscriptstyle M} $  {\bf
(multiparameter)
quantum Poisson pairing}.
                                                      \par
   By previous analysis, the integer forms of quantum enveloping algebras are
$ \kqqm $--dual  of each other (cf.~\S 2.6) with respect to
$ \pi_{q,\varphi}^{\scriptscriptstyle M} \, $;  so the latter
restrict to perfect pairings
  $$  \pi^\varphi_{q,H^\tau_{\scriptscriptstyle M}} \colon \, {\gerUphiMh
{M'} } \otimes {\calUphiMg M } \longrightarrow \kqqm \, ,  \quad
\pi^\varphi_{q, G^\tau_{\scriptscriptstyle M}} \colon \, {\calUphiMh M }
\otimes {\gerUphiMg {M'} } \longrightarrow \kqqm \; ;  $$
same symbols will also denote the Hopf pairings
$ \, \pi^\varphi_{q,H^\tau_{\scriptscriptstyle M}} \colon \,
{\calFphiMg {M'} } \otimes {\calUphiMg M } \longrightarrow \kqqm \, $,
resp.  $ \, \pi^\varphi_{q, G^\tau_{\scriptscriptstyle M}} \colon \,
{\gerFphiMg M } \otimes {\gerUphiMg {M'} } \longrightarrow \kqqm \, $,
got by restriction of the previous ones: hereafter we identify
$ {\fqphiMg M } $  with its image in  $ {\uqphiMh M } $  via  $
\xi^\varphi_{\scriptscriptstyle M} $,  and similarly for integer forms.

\vskip1,7truecm

 \centerline{ \bf  \S \; 7 \,  Specialization at roots of 1 }

\vskip10pt

   {\bf 7.1  The case  $ \, q \rightarrow 1 \, $: specialization of
$ \gerUphiMh{M} $  to  $ \uhtau $  and consequences.}  \;  Recall (cf.~\S 2.1)
that  $ \, \tau = (\tau_1, \ldots, \tau_n) := {\,1\, \over \,2\,}
(\varphi(\alpha_1), \dots, \varphi(\alpha_n)) \, $.  Now set
  $$  \gerUunophiMh{M} := \, {\gerUphiMh M } \Big/ (q-1) \, {\gerUphiMh M }
\, \cong \, {\gerUphiMh M } \otimes_{k \left[ q, \qm \right]} k  $$
let  $ \, p^\varphi_1 \colon \, {\gerUphiMh M } \rightarrow
\gerUunophiMh{M} \, $  be the canonical projection, and
set  $ \, \text{f}^{\,\tau}_i := p^\varphi_1 \left( {\left( F^\varphi_i
                                      \right)}^{(1)} \right) $,\break
$ \text{m}^\tau_i  := p^\varphi_1 \left( \left( M^\varphi_i ; \, 0 \atop 1
\right) \right) $,  $ \text{e}^\tau_i := p^\varphi_1 \left( {\left( E^\varphi_i
\right)}^{(1)}\right) $,  (where  $ \, M_i := L_{\mu_i}^\varphi \, $)  for
all  $ i=\unon $.

\vskip7pt

\proclaim{Theorem 7.2}  For  $ \, q \rightarrow 1 \, $,
$ \, {\gerUphiMh M } $  specializes to the Poisson Hopf coalgebra
$ \uhtau \, $;  in other words, there exists an isomorphism of Poisson Hopf
coalgebras
  $$  \gerUunophiMh{M} \cong \uhtau \, .  $$
\endproclaim

\demo{Proof}  The proof mimick that for  $ \, \gerUunophiMg{M} \cong \ugtau
\, $.  From the presentation of  $ {\gerUphiMh M } $  we get  $ \,
\gerUunophiMh{M} = {\frak H}_\varphi^{\scriptscriptstyle M}{\big\vert}_{q=1} :=
{\frak H}_\varphi^{\scriptscriptstyle M} \big/ (q-1) \,
{\frak H}_\varphi^{\scriptscriptstyle M} \, $,  hence we are reduced to study
$ \, {\frak H}_\varphi^{\scriptscriptstyle M}{\big\vert}_{q=1} \, $;
moreover the presentation of  $ {\frak H}_\varphi^{\scriptscriptstyle M} $
provides one of  $ \, {\frak H}_\varphi^{\scriptscriptstyle
M}{\big\vert}_{q=1} \, $.  Now the definition of  $ {\frak
H}_\varphi^{\scriptscriptstyle M}{\big\vert}_{q=1} $  and the explicit form of
the PBW basis of  $ {\gerUphizM M } $  (cf.~\S 2.5) imply that the elements
$ \, {\left( F^\varphi_i \right)}^{(r)} $,  $ \left( M^\varphi_i ; \, 0 \atop t
\right) $,  $ {\left( M^\varphi_i \right)}^{-1} $,  $ {\left( E^\varphi_i
\right)}^{(s)} \, $  ($ \, i= \unon; \, r, t, s \in \N \, $)  are enough to
generate  $ {\frak H}_\varphi^{\scriptscriptstyle M} \, $;  finally,
straightforward computation gives  $ \; p^\varphi_1 \left( {\left( F^\varphi_i
\right)}^{(r)}\right) = {{\ {\left( \text{f}^{\,\tau}_i \right)}^{\! r} \ }
\over {\ r! \ }} \, $,  $ \; p^\varphi_1 \left( \left( M^\varphi_i ; 0 \atop t
\right) \right) = {\text{m}^\tau_i \choose t} \, $,  $ \; p^\varphi_1 \left(
{\left( M^\varphi_i \right)}^{-1} \right) = 1 \, $,  $ \; p^\varphi_1 \left(
{\left( E^\varphi_i \right)}^{(s)}\right) = {{\ {\left( \text{e}^\tau_i
\right)}^s \ } \over {\ s! \ }} \; $  (where  $ \, {\text{m}^\varphi_i \choose
t} := {\text{m}^\varphi_i (\text{m}^\varphi_i - 1) (\text{m}^\varphi_i - 2)
\cdots (\text{m}^\varphi_i - t + 1) \over t!} \, $),  hence  $ \,
\gerUunophiMh{M} = {\frak H}_\varphi^{\scriptscriptstyle M} {\Big\vert}_{q=1}
\, $  is generated by the  $ \text{f}^{\,\tau}_i $'s,  $ \text{m}^\tau_i $'s,
$ \text{e}^\tau_i $'s,  with some relations.
                                                 \par
  When  $ \, M = Q \, $  this presentation is exactly the same of  $ \uhtau $
(cf.~(1.2)),  with  $ \, \hbox{h}^\tau_i = \hbox{m}^\tau_i \, $;  comparing
(1.3) with formulas in \S 6.9  (for  $ q=1 $)  shows that also the Hopf
structure is the same.  In particular  $ \gerUunophiMh{Q} $  is cocommutative,
hence has a canonical co-Poisson structure given by  $ \, \delta := \big({\,
\Delta - \Delta^{op} \, \over \, q - 1 \,}\big){\big\vert}_{q=1} \, $,
described by formulas   --- deduced from those in \S 8.9 ---   which do coincide
with (1.4), as a straightforward checking shows.
                                                 \par
   Finally, for  $ \, M \neq Q \, $  we prove that
$ \, \gerUunophiMh{M} \cong \gerUunophiMh{Q} \, $
as Poisson Hopf coalgebras: since  $ \, \gerUphiMh{M} \supseteq
\gerUphiMh{Q} \, $  by definition, it is
enough to check that  $ \, {\frak H}_\varphi^{\scriptscriptstyle M}
{\big\vert}_{q=1} = {\frak H}_\varphi^{\scriptscriptstyle Q}
{\big\vert}_{q=1} \, $  as  $ k $--vector  spaces.  Assume we are in the simply
laced case.  Since  $ \, M^\varphi_i := L^\varphi_{\mu_i} \, $  and  $ \,
K^\varphi_j := L^\varphi_{\alpha_j} \, $,  it is  $ \, K^\varphi_j :=
\prod_{i=1}^n {\left( M^\varphi_{\alpha_i} \right)}^{c_{ij}} \, $,  where  $ \,
c_{ij} \in \Z \, $  are such that  $ \, \alpha_j = \sum_{i=1}^n c_{ij} \mu_i
\, $.  Then  $ \; \left( K^\varphi_j ; 0 \atop t \right) {\Big\vert}_{q=1} =
\sum_{i=1}^n c_{ij} \cdot \left( M^\varphi_i ; 0 \atop t \right)
{\Big\vert}_{q=1} \; $  so that  $ \, {\gerUphizM M }{\Big\vert}_{q=1} =
{\gerUphizM Q }{\Big\vert}_{q=1} \, $  follows, whence  $ \, {\frak
H}_\varphi^{\scriptscriptstyle M}{\big\vert}_{q=1} =
{\frak H}_\varphi^{\scriptscriptstyle Q} {\big\vert}_{q=1} \, $,  q.e.d.  In the
other cases  $ \, M = P \, $,  and this argument still works,  {\sl mutatis
mutandis},  because  $ \, \alpha_j = \sum_{i=1}^n \aij \omega_i \, $,  hence
$ \, K^\varphi_j := \prod_{i=1}^n {\left( L^\varphi_{\alpha_i} \right)}^{\aij}
\, $,  so that  $ \; \left( K^\varphi_j ; 0 \atop t \right) {\Big\vert}_{q=1} =
\sum_{i=1}^n a_{ji} \cdot \left( L^\varphi_i ; 0 \atop t \right)
{\Big\vert}_{q=1} \; $  and we are done again.   $ \square $
\enddemo

\vskip7pt

   {\sl Remark:}  {\it Thus  $ {\gerUphiMh M } $  provides the announced
infinitesimal quantization of}  $ H^\tau $.  This can be partially
explained as follows.  $ {\gerUphiMh M } $  is a subspace of
$ {\uqphiMg {M'} }^* $  made of series satisfying a certain "growth
condition" (cf.~(6.3)): then specializing  $ q $  at  $ 1 $  one gets an
isomorphism  {\sl of Hopf algebras}  $ \; \gerUunophiMh{M} \cong \Big\{\, f \in
{F \left[ H^\tau_{\scriptscriptstyle M} \right]}^* \,\Big\vert\, \exists\,
n \in \N : f \left( {\frak e}^n \right) = 0 \,\Big\} \; $  where
$ \, {\frak e} := Ker \left( \epsilon \colon \, F \left[
H^\tau_{\scriptscriptstyle M} \right] \rightarrow k \right) \, $,  and
$ \, {\frak e} = {\frak m}_e \, $,  where  $ {\frak m}_e $  is the maximal
ideal of  $ F \left[ H^\tau_{\scriptscriptstyle M} \right] $  associated to
$ \, e \in H^\tau_{\scriptscriptstyle M} \, $.  Since  $ \; \Big\{\, f \in
{F \left[ H^\tau_{\scriptscriptstyle M} \right]}^* \,\Big\vert\, \exists\,
n \in \N : f \left( {{\frak m}_e}^{\! n} \right) = 0 \,\Big\} \cong \uhtau \; $
as Hopf algebras (cf.~for instance [On], Part I, Ch.~3, \S 2), we conclude
that there exists a Hopf algebra isomorphism  $ \; \gerUunophiMh{M} \cong
\uhtau \, . $  But regarding co-Poisson structure, such an analysis gives no
information, thus the proof of Theorem 7.2 given above is really necessary.

\vskip7pt

   The previous theorem has two interesting consequences.  As for the first, set
  $$  \calFunophiMg{M} := \, {\calFphiMg M } \Big/ (q-1) \, {\calFphiMg M }
\, \cong \, {\calFphiMg M } \otimes_{k \left[ q, \qm \right]} k \; .  $$

\vskip7pt

\proclaim{Theorem 7.3}  The Hopf algebra  $ \, {\calFphiMg
M } \, $  specializes to the Poisson Hopf coalgebra
$ \uhtau $  for  $ \, q \rightarrow 1 \, $;  in other words, there exists an
isomorphism of Poisson Hopf coalgebras
  $$  \calFunophiMg{M} \cong \uhtau \, .  $$
\endproclaim

\demo{Proof}  Consider the monomorphism  $ \, \xi^\varphi_{\scriptscriptstyle M}
\colon \, {\calFphiMg M } \llonghookrightarrow {\gerUphiMh M } \, $
and compare it with the isomorphism  $ \, \xi^\varphi_{\scriptscriptstyle M}
\colon \, {\calFphiMg M } \! \left[ \psi_{-\rho}^{-1} \right] \, {\buildrel
\cong \over \longrightarrow} \, {\frak H}_\varphi^{\scriptscriptstyle M}
\subseteq \gerUphiMh{M} \, $.  When  $ \, q \rightarrow 1 \, $  we get
$ \, \gerUunophiMh{M} = {\frak H}_\varphi^{\scriptscriptstyle
M}{\big\vert}_{q=1} = \Big( {\calFphiMg M } \! \left[ \psi_{-\rho}^{-1}
\right] \Big){\Big\vert}_{q=1} = \calFunophiMg{M} \! \left[
\psi_{-\rho}^{-1}{\big\vert}_{q=1} \right] \, $;
but  $ \, \xi^\varphi_{\scriptscriptstyle M} \left( \psi_{-\rho}^{-1} \right)
= L^\varphi_\rho = \prod_{i=1}^n M^\varphi_i \, $  (cf.~Theorem 5.14),  hence
$ \, \xi^\varphi_{\scriptscriptstyle M} \left( \psi_{-\rho}^{-1}
\right){\big\vert}_{q=1} = \prod_{i=1}^n M^\varphi_i {\big\vert}_{q=1} = 1 \, $
because  $ \, M^\varphi_i = 1 + (q_i - 1) \cdot \left( M^\varphi_i ; \, 0 \atop
                                1 \right) \equiv 1 $\break
$ \mod (q - 1) \, \, $.  Therefore  $ \, \gerUunophiMh{M} \cong
\calFunophiMg{M} \! \left[ \psi_{-\rho}^{-1}{\big\vert}_{q=1} \right] =
\calFunophiMg{M} \, $,  whence the thesis.   $ \square $
\enddemo

\vskip7pt

   {\sl Remark:}  thus  $ \fqphiMg{M} $  too yields an infinitesimal
quantization of  $ H^\tau \, $;  compared with  $ \uqphiMh{M} $  the
advantage is that  $ \fqphiMg{M} $  is  {\sl usual}  Hopf algebra, whereas
$ \uqphiMh{M} $  (or  $ \gerUphiMh{M} $)  is a  {\sl topological}  Hopf
algebra.  Furthermore, for the classical groups there exists a presentation
of  $ F_{q,0}^{\scriptscriptstyle P} \left[ G \right] $  by generators and
relations, hence   --- at least in principle ---   one can study
$ \calF_{\varphi,0}^{\scriptscriptstyle P} \left[ G \right] $  exploiting
such a presentation.  For  $ \, G = SL(n+1) \, $  this is done in [Ga].

\vskip7pt

   Theorem 7.3 gives  $ {\calFphiMg M } @>{q \rightarrow 1}>> \uhtau $,
\  which is the dual result (in the sense of Poisson duality) \
$ {\calUphiMg M' } @>{q \rightarrow 1}>> F[H^\tau_{\scriptscriptstyle M}] \, $.
The original proof of the latter result in [DKP] (see also [DP]) is lenghty
involved and complicated, requiring very hard computations; on the contrary, we
can deduce it as an easy consequence of Theorem 7.2:

\vskip7pt

\proclaim{Theorem 7.4}  The Hopf algebra  $ \, {\calUphiMg M } \, $
specializes to the Poisson Hopf algebra  $ F \left[ H^\tau_{\scriptscriptstyle
M} \right] $  for  $ \, q \rightarrow 1 \, $,  in other words, there exists
an isomorphism of Poisson Hopf algebras
  $$  \calUunophiMg{P} := \, {\calUphiMg M } \Big/ (q-1) \, {\calUphiMg M }
\, \cong \, F \left[ H^\tau_{\scriptscriptstyle M} \right] \; .  $$
\endproclaim

\demo{Proof}  Since  $ \, {\calUphiMg M } \, $  is perfectly paired with
$ \, {\gerUphiMh M' } \, $,  we have that  $ \, \calUunophiMg{M} \, $  is
perfectly paired with  $ \, \gerUunophiMh{M'} \cong \uhtau \, $:  the latter
is cocommutative, hence the former is commutative.  Then $ \calUunophiMg{M} $
is a finitely generated commutative Hopf algebra over  $ k $,  hence it is
the algebra of (regular) functions of an affine algebraic group, say
$ H' \, $;  moreover  $ \, \calUunophiMg{M} = F[H'] \, $
inherits from  $ {\calUphiMg M } $  a Poisson structure, so  $ H' $
is a Poisson group.  Like in [DP] it is clear from the presentation of
$ {\calUphiMg M } $  that  $ \, F[H'] \, \big( = \calUunophiMg{M} \,\big)
\cong F \left[ H^\tau_{\scriptscriptstyle M} \right] \, $  as Hopf algebras,
hence  $ \, H' = H^\tau_{\scriptscriptstyle M} \, $  as algebraic groups
(the non-trivial part in [DP] is that dealing with Poisson structures).
Now the Hopf pairing among  $ \, \gerUunophiMh{M'} \cong \uhtau \, $  and
$ \, \calUunophiMg{M} = F[H'] = F \left[ H^\tau_{\scriptscriptstyle M}
\right] \, $  is compatible with Poisson and co-Poisson structures, that is
$ \, \big\langle h, \{f,g\} \big\rangle = \big\langle \delta(h), f \otimes g
\big\rangle \, $,  where  $ \delta $  is the Poisson cobracket
of  $ \, \gerUunophiMh{M'} = \uhtau \, $  and  $ \{\ ,\ \} $
is either the Poisson bracket  $ {\{\ ,\ \}}_\star $  of
$ H^\tau_{\scriptscriptstyle M} $  or the Poisson bracket
$ {\{\ ,\ \}}_\circ $  of  $ H' \, $:  since the pairing is
perfect, we must have  $ \, {\{\ ,\ \}}_\star =
{\{\ ,\ \}}_\circ \, $,  whence the thesis.   $ \square $
\enddemo

\vskip7pt

   {\bf  7.5  The case  $ \, q \rightarrow 1 \, $:  specialization of  $
{\calUphiMh M } $  to  $ F^\infty \left[ G^\tau_{\scriptscriptstyle M}
\right] $.}  \  We are going to show that  $ {\calUphiMh M } $  is a
quantization of  $ F^\infty \left[ G^\tau_{\scriptscriptstyle M} \right]
(=\finfgtau) $;  such a result can be seen as (Poisson) dual counterpart of  $
\, {\calUphiMg M } \, @>{q \rightarrow 1}>> \, F \left[
H^\tau_{\scriptscriptstyle M} \right] \, $  (cf.~Theorem 7.4).  As usual, we set
  $$  \calUunophiMh{M} := \, {\calUphiMh M } \Big/ (q-1) \, {\calUphiMh M }
\, \cong \, {\calUphiMh M } \otimes_{k \left[ q, \qm \right]} k \; .  $$

\vskip7pt

\proclaim{Theorem 7.6}  The formal Hopf algebra
$ {\calUphiMh M } $  specializes to the formal Poisson Hopf algebra
$ F^\infty \left[ G^\tau_{\scriptscriptstyle M} \right] (=\finfgtau) $
for  $ \, q \rightarrow 1 \, $;  in other words, there exists an
isomorphism of formal Poisson Hopf algebras
  $$  \calUunophiMh{M} \cong F^\infty \left[ G^\tau_{\scriptscriptstyle M}
\right] \, .  $$
\endproclaim

\demo{Proof}  Recall that  $ \, F^\infty \left[ G^\tau_{\scriptscriptstyle M}
\right] = \finfgtau \, $  is isomorphic to the linear dual of  $ \ugtau $,
that is  $ \, F^\infty \left[ G^\tau_{\scriptscriptstyle M} \right]
\cong {\ugtau}^* \, $.  On the other hand, we have a formal Hopf
algebra isomorphism  $ \, {j_{\scriptscriptstyle M}}^{\! -1}
\smallcirc \nu^\varphi_{\scriptscriptstyle M} \colon \, {\uqphiMh M }
@>\cong>> {\uqphiMg {M'} }^* \, $,  and Theorem 6.6 ensures that this
restricts to
  $$  {j_{\scriptscriptstyle M}}^{\! -1} \smallcirc
\nu^\varphi_{\scriptscriptstyle M} \colon \, {\calUphiMh M } {\buildrel \cong
\over \llongrightarrow} {\gerUphiMg {M'} }^* \, .   \eqno (7.3)  $$
   \indent   When  $ \, q \rightarrow 1 \, $,  we have that
$ {\gerUphiMg {M'} } $  specializes to  $ \ugtau $,
                              therefore (7.3) implies\break
$ \, \calUunophiMh{M} \cong {\gerUphiMg {M'} }^*
\otimes_{k \left[ q, \qm \right]} k = {\gerUunophiMg{M'} }^* \cong {\ugtau}^*
= \finfgtau = F^\infty \left[ G^\tau_{\scriptscriptstyle M} \right] \, $,
q.e.d.   $ \square $
\enddemo

\vskip7pt

   {\bf 7.7  The case  $ \, q \rightarrow \varepsilon \, $:
quantum Frobenius morphisms.}  \;  Let  $ \varepsilon $  be a primitive
$ \ell $--th  root of 1 in  $ k $,  for  $ \ell $  {\it odd\/},  $ \, \ell
> d:= \max_i \{d_i\} \, $,  and set
  $$  \gerUepsilonphiMh{M} := \, {\gerUphiMh M } \Big/ (q - \varepsilon) \,
{\gerUphiMh M } \, \cong \, {\gerUphiMh M }
\otimes_{k \left[ q, \qm \right]} k  $$
   \indent   First of all we remark that
  $$  \hbox{\it $ \gerUepsilonphiMh{M} $  is a  {\rm usual}  Hopf algebra
over  $ k $,  isomorphic to  $ {\frak H}_\varphi^{\scriptscriptstyle M}
{\Big\vert}_{q=\varepsilon} $}   \eqno (7.4)  $$
for every element of  $ \gerUepsilonphiMh{M} $  is a formal series of terms
in  $ {\frak H}_\varphi^{\scriptscriptstyle M} $  {\sl which is a finite sum
modulo}  $ (q - \varepsilon) $,  and \S 5.16 and Theorem 6.2 tell us that
$ \, \Delta \big( {\frak H}_\varphi^{\scriptscriptstyle
M}{\big\vert}_{q=\varepsilon} \big) \subseteq {\frak
H}_\varphi^{\scriptscriptstyle M}{\big\vert}_{q=\varepsilon} \otimes {\frak
H}_\varphi^{\scriptscriptstyle M}{\big\vert}_{q=\varepsilon} \, $,  and  $ \, S
\big( {\frak H}_\varphi^{\scriptscriptstyle M}{\big\vert}_{q=\varepsilon} \big)
= {\frak H}_\varphi^{\scriptscriptstyle M}{\big\vert}_{q=\varepsilon} \, $.  Now
we are ready for next result, the analogue for  $ {\uqphiMh M } $  of (3.6).

\vskip7pt

\proclaim{Theorem 7.8}  There exists a Hopf algebras epimorphism
  $$  \gerFrhtau : \gerUepsilonphiMh{M} \llongtwoheadrightarrow
\gerUunophiMh{M} \cong \uhtau  $$
defined (for all  $ \, i= 1, \dots, n \, $)  by
  $$  {} \!\!  \gerFrhtau \colon
 \cases
   {F^\varphi_i}^{(s)} \Big\vert_{q=\varepsilon} \!\!\!\!\! \mapsto
\! {F^\varphi_i}^{(s / \ell)} \Big\vert_{q=1},  \, \left( M^\varphi_i;
\, 0 \atop s \right) \! \bigg\vert_{q=\varepsilon} \!\!\!\!\! \mapsto \!
\left( M^\varphi_i; \, 0 \atop s / \ell \right) \! \bigg\vert_{q=1},
\, {E^\varphi_i}^{(s)} \Big\vert_{q=\varepsilon} \!\!\!\!\! \mapsto \!
{E^\varphi_i}^{(s / \ell)} \Big\vert_{q=1}  \;\; \hbox{ if } \; \ell
\Big\vert s  \\
   {F^\varphi_i}^{(s)} \Big\vert_{q=\varepsilon} \!\!\! \mapsto 0 \, ,  \quad
\left( M_i; \, 0 \atop s \right) \! \bigg\vert_{q=\varepsilon} \!\!\! \mapsto 0
\, ,  \quad {E^\varphi_i}^{(s)} \Big\vert_{q=\varepsilon} \!\!\! \mapsto 0
\quad  \hbox{ \, otherwise \, }  \\
   \big({M^\varphi_i}^{-1}\big) \Big\vert_{q=1} \!\!\! \mapsto 1  \\
 \endcases  $$
which is adjoint of  $ \calFrgtau $  (cf.~(3.9)) with respect to the
quantum Poisson pairings, that is
  $$  \pi^\varphi_{1, H^\tau_{\scriptscriptstyle {M'}}} \! \big( \gerFrhtau(h),
g \big) = \pi^\varphi_{\varepsilon, H^\tau_{\scriptscriptstyle {M'}}} \! \big(
h, \calFrgtau(g) \big)   \eqno \forall \;\, h \in {\gerUepsilonphiMh M }, \, g
\in {\calUunophiMg M' } \, .  \qquad  $$
\endproclaim

\demo{Proof}  The formulas above uniquely determine an epimorphism
$ \gerFrhtau $   --- if any ---   because
$ \, {\left( F^\varphi_i \right)}^{(s)} \Big\vert_{q=\varepsilon} $,
$ \left( M^\varphi_i; 0 \atop s \right) \Big\vert_{q = \varepsilon} $,
$ {\left( M^\varphi_i \right)}^{-1}\Big\vert_{q=\varepsilon} $,
$ {\left( E^\varphi_i \right)}^{(s)} \Big\vert_{q=\varepsilon} \, $  are
(algebraic) generators of  $ \, {\frak H}_\varphi^{\scriptscriptstyle
M}{\Big\vert}_{q=\varepsilon}
                       \!\! = $\break
$ = \gerUepsilonphiMh{M} \, $  (cf.~(7.4)).
Consider the embedding  $ \, \calFrgtau \colon \,
F \left[ H^\tau_{\scriptscriptstyle M'} \right] \! \cong \! \calUunophiMg{M'}
\longhookrightarrow \calUepsilonphiMg{M'} \, $  of Hopf algebras (cf.~(3.9)):
its linear dual is an epimorphism of formal Hopf algebras  $ \,
{\calUepsilonphiMg{M'} }^* \longtwoheadrightarrow {\calUunophiMg{M'} }^* \, $.
On the other hand we have an embedding  $ \, \gerUepsilonphiMh{M}
\longhookrightarrow {\calUepsilonphiMg{M'} }^* \, $  provided by
the specialized quantum Poisson pairing
$ \; \pi^\varphi_{\varepsilon,H^\tau_{\scriptscriptstyle {M'}}}
\colon \, \gerUepsilonphiMh{M} \otimes \calUepsilonphiMg{M'} \longrightarrow
k \; $.
                                                  \par
   Now composition yields a morphism  $ \, \gerFrhtau \colon
\, \gerUepsilonphiMh{M} \llongrightarrow {\calUunophiMg{M'} }^* \, $;
the very construction then gives  $ \, \left\langle \gerFrhtau(h), g
\right\rangle = \pi^\varphi_{1, H^\tau_{\scriptscriptstyle {M'}}} \!
\left( \gerFrhtau(h), g \right) = \pi^\varphi_{\varepsilon,
H^\tau_{\scriptscriptstyle {M'}}} \! \left( h, \calFrgtau(g) \right) \, $,
hence  $ \gerFrhtau $  is adjoint of  $ \calFrgtau(g) $,  is described by the
previous formulas and has image  $ \gerUunophiMh{M} $,  q.e.d.   $ \square $
\enddemo

\vskip7pt

   Similar arguments prove next result, which is the analogue of (3.9); as
usual, we set
  $$  \calUepsilonphiMh{M} := \, {\calUphiMh M } \Big/ (q - \varepsilon) \,
{\calUphiMh M } \, \cong \, {\calUphiMh M } \otimes_{k \left[ q, \qm \right]}
k \, .  $$

\vskip7pt

\proclaim{Theorem 7.9}
                                          \hfill\break
   \indent   (a) There exists a unique continuous monomorphism of formal Hopf
algebras
  $$  \calFrhtau : F^\infty \left[ G^\tau_{\scriptscriptstyle M}
\right] \cong \calUunophiMh{M} \llonghookrightarrow \calUepsilonphiMh{M}
\eqno (7.5)  $$
defined (for all  $ \, \alpha \in R^+ \, $,  $ \mu \in M \, $)  by
  $$  \calFrhtau : \;  \quad  \fbar^\varphi_\alpha \Big\vert_{q=1} \mapsto
{\big( \fbar^\varphi_\alpha \Big)}^{\! \ell} \Big\vert_{q=\varepsilon}
\, , \; \quad L^\varphi_\mu \Big\vert_{q=1} \mapsto {\big( L^\varphi_\mu
\big)}^\ell \Big\vert_{q=\varepsilon} \, ,  \; \quad \ebar^\varphi_\alpha
\Big\vert_{q=1} \mapsto {\Big( \ebar^\varphi_\alpha \Big)}^{\! \ell}
\Big\vert_{q=\varepsilon}   \eqno (7.6)  $$
which is the continuous extension of  $ \gerFrGtau $  (cf.~(4.5))  and is
adjoint of  $ \gerFrgtau $  (cf.~(3.6)) with respect to quantum Poisson
pairings, that is
  $$  \pi^\varphi_{\varepsilon,G^\tau_{\scriptscriptstyle M}} \! \big(
\calFrhtau(h), g \big) = \pi^\varphi_{1,G^\tau_{\scriptscriptstyle M}} \!
\big( h, \gerFrgtau(g) \big)   \eqno \forall \;\, h \in {\calUunophiMh
M }, \, g \in {\gerUepsilonphiMg M' } \, .  \qquad  $$
                                            \hfill\break
   \indent   (b) The image  $ \, Z^\varphi_0 \; \left( \, \cong_{\calFrhtau}
\calUunophiMh{M} \, \right) \, $  of  $ \calFrhtau $  is a formal Hopf
subalgebra contained in the centre of  $ \calUepsilonphiMh{M} $.
                                            \hfill\break
   \indent   (c) The set  $ \Big\{ \calF_{\ell \phi} \cdot
B^\varphi_{\scriptscriptstyle \ell \phi, \ell \tau, \ell \eta} \cdot
\calE_{\ell \eta} \, \Big\vert \, \phi \! \in \! \N^N, \tau \! \in \! \N^n,
\eta \! \in \! \N^N \Big\} $  is a pseudobasis of  $ Z^\varphi_0 $
over  $ k $.
                                            \hfill\break
   \indent   (d) The set  $ \, \Big\{\, \calF_\phi \cdot
B^\varphi_{\scriptscriptstyle \phi, \tau, \eta} \cdot \calE_\eta \,
\Big\vert \, \tau \in {\big\{ 0, 1, \dots, \ell-1 \big\}}^n; \, \phi,
\eta \in {\big\{ 0, 1, \dots, \ell-1 \big\}}^N \,\Big\} \, $  is a basis
of  $ \calUepsilonphiMh{M} $  over  $ Z_0 \, $;  therefore also
                            the set of ordered PBW monomials\break
$ \; \Big\{\, \calF_\phi \cdot M^\varphi_\mu \cdot \calE_\eta \,
\Big\vert \, \mu \in {\big\{ 0, 1, \dots, \ell-1 \big\}}^n; \, \phi, \eta
\in {\big\{ 0, 1, \dots, \ell-1 \big\}}^N \,\Big\} \; $  is a basis of
$ \calUepsilonphiMh{M} $  over  $ Z^\varphi_0 \, $.  Thus
$ \calUepsilonphiMh{M} $  is a free module of rank  $ \ell^{dim \left(
H^\tau \right)} $  over  $ Z^\varphi_0 \, $.
\endproclaim

\demo{Proof}  {\it (a)\/}  Since  $ \, \fbar^\varphi_\alpha \Big\vert_{q=1} $,
$ L^\varphi_\mu \Big\vert_{q=1} $,  $ \ebar^\varphi_\alpha \Big\vert_{q=1} $
($ \alpha \in R^+ $,  $ \mu \in M \, $)  are topological generators of
$ \calUunophiMh{M} $,  the formulas above uniquely determine a continuous
monomorphism  $ \calFrhtau $,  if any.  Now consider  $ \, \gerFrgtau \colon
\, \gerUepsilonphiMg{M'} \longtwoheadrightarrow \gerUunophiMg{M'} \cong
\ugtau \, $  (cf.~(3.6)), a Hopf epimorphism, and its dual, a formal Hopf
monomorphism  $ \, {\gerUunophiMg{M'} }^* \longhookrightarrow
{\gerUepsilonphiMg{M'} }^* \, $;  composing the latter with the
isomorphisms  $ \, \calUunophiMh{M} @>\cong>> {\gerUunophiMg{M'} }^* \, $,
$ \, {\gerUepsilonphiMg{M'} }^* @>\cong>> \calUepsilonphiMh{M} \, $  (given
by specialized quantum Poisson pairings) provides a monomorphism
$ \, \calFrhtau \colon \, \calUunophiMh{M} \longhookrightarrow
\calUepsilonphiMh{M} \, $;  then
  $$  \big\langle \calFrhtau(h), g \big\rangle =
\pi^\varphi_{\varepsilon,G^\tau_{\scriptscriptstyle M}} \! \big(
\calFrhtau(h), g \big) = \pi^\varphi_{1, G^\tau_{\scriptscriptstyle M}} \!
\big( h, \gerFrgtau(g) \big)   \eqno \forall \;\, h \in {\calUunophiMh M }, \, x
\in {\gerUepsilonMg M' }  \quad  $$
hence  $ \calFrhtau $  is described by formulas above.  Moreover, with
notation of \S 6.1 and \S 6.3,
  $$  \displaylines{
   \left\langle \calFrhtau \big( Y^\varphi_{\scriptscriptstyle \phi, \zeta,
\eta} \big), X_{\scriptscriptstyle \epsilon, \theta, \psi} \right\rangle =
\pi^\varphi_{\varepsilon,G^\tau_{\scriptscriptstyle M}} \left( \calFrhtau \big(
Y^\varphi_{\scriptscriptstyle \phi, \zeta, \eta} \big), X_{\scriptscriptstyle
\epsilon, \theta, \psi} \right) = \pi^\varphi_{1,G^\tau_{\scriptscriptstyle M}}
\left( Y^\varphi_{\scriptscriptstyle \phi, \zeta,
\eta}, \, \gerFrgtau \left( X_{\scriptscriptstyle \epsilon, \theta, \psi}
\right) \right) =  \cr
   = \! \chi_{{}_{\scriptstyle \ell \N^{\,N}}}\!(\epsilon) \cdot
\chi_{{}_{\scriptstyle \ell \N^{\,n}}}\!(\theta) \cdot
\chi_{{}_{\scriptstyle \ell \N^{\,N}}}\!(\psi) \cdot \left\langle
Y^\varphi_{\scriptscriptstyle \phi, \zeta, \eta},
X_{\scriptscriptstyle {\,1\, \over \,\ell\,} \cdot \epsilon,
{\,1\, \over \,\ell\,} \cdot \theta, {\,1\, \over \,\ell\,}
\cdot \psi} \right\rangle = \chi_{{}_{\scriptstyle \ell \cdot \left( \N^{\,N}
\times \N^{\,n} \times \N^{\,N} \right)}} \! (\epsilon, \theta, \psi) \cdot
\delta_{\ell (\phi, \zeta, \eta), (\epsilon, \theta, \psi)}  \cr }  $$
(where  $ \chi_{{}_{\scriptstyle \ell S}} $  is the characteristic function of
the sublattice  $ \ell \, S \subseteq S \, $,  $ S $  any abelian semigroup),
hence  $ \, \calFrhtau \big( Y^\varphi_{\scriptscriptstyle
\phi, \zeta, \eta} \big) = Y^\varphi_{\scriptscriptstyle \ell \phi, \ell \zeta,
\ell \eta} \, $  for all  $ \phi $,  $ \zeta $,  $ \eta $,  thus  $ \calFrhtau $
 maps elements of the pseudobasis (6.2) of  $ \calUunophiMh{M} $  onto elements
of the analogous pseudobasis of  $ \calUepsilonphiMh{M} $:  therefore
$ \calFrhtau $  is continuous.
                                                    \par
   Finally, since  $ \, \gerFrGtau \colon \, F \left[
G^\tau_{\scriptscriptstyle M} \right] \cong \gerFunophiMg{M}
\hookrightarrow \gerFepsilonphiMg{M} \, $  too is defined as (Hopf)
dual of  $ \, \gerFrgtau \colon \, \gerUepsilonphiMg{M'}
\twoheadrightarrow \gerUunophiMg{M'} \cong \ugtau \, $
(cf.~[DL], Proposition 6.4), then  $ \, \calFrhtau \colon
F^\infty \left[ G^\tau_{\scriptscriptstyle M} \right] \cong
\calUunophiMh{M} \hookrightarrow \calUepsilonphiMh{M} \; $  is
extension of  $ \; \gerFrGtau \colon F \left[ G^\tau_{\scriptscriptstyle M}
\right] \cong \gerFunophiMg{M} \hookrightarrow \gerFepsilonphiMg{M} \; $;
since  $ {\gerFphiMg M } $  is dense in  $ \, {\gerUphiMg {M'} }^* \cong
{\calUphiMh M } \, $  it is clear that this extension is by continuity.
 \eject
   {\it (b)\/}  This easily follows from the analogous result for
$ \calUepsilonphiMg{M} $  (cf.~[DP], Theorem 19.1)  and comparison
among  $ \calUepsilonphiMg{M} $  and  $ \calUepsilonphiMh{M} $.
                                                   \par
  {\it (c)\/}  This follows from the previous analysis, namely from  $ \,
\calFrhtau \big( Y^\varphi_{\scriptscriptstyle \phi, \zeta, \eta} \big) =
Y^\varphi_{\scriptscriptstyle \ell \phi, \ell \zeta, \ell \eta} \, $.
                                                   \par
  {\it (d)\/}  The span of  $ \, \big\{\, B^\varphi_{\scriptscriptstyle \phi,
\zeta, \eta} \,\big\vert\, (\phi, \zeta, \eta) \in \ell \big( \N^N \times
\N^n \times \N^N \big) \,\big\} \, $  (inside  $ \calUepsilonphiMh{M} $)
coincides with the span of  $ \, \big\{\, L^\varphi_\mu \,\big\vert\, \lambda
\in \ell \N^n = \ell M_+ \,\big\} \, $;  from this and from the explicit form of
the pseudobasis of  $ {\calUphiMh M } $  we get the claim.   $ \square $
\enddemo

\vskip7pt

   At last we prove the dual counterpart of (4.5), regarding
  $$  \calFepsilonphiMg{M} := \, {\calFphiMg M } \Big/ (q - \varepsilon)
\, {\calFphiMg M } \, \cong \, {\calFphiMg M }
\otimes_{k \left[ q, \qm \right]} k \; ;  $$
note that we obtain a quantum Frobenius morphism which is
{\sl surjective\/}  instead of  {\sl injective\/}.

\vskip7pt

\proclaim{Theorem 7.10}  There exists a Hopf algebra epimorphism
  $$  \calFrHtau \colon \, \calFepsilonphiMg{M} \llongtwoheadrightarrow
\calFunophiMg{M} \cong \uhtau   \eqno (7.7)  $$
dual of  $ \; \calFrgtau \colon \, F \left[ H^\tau_{\scriptscriptstyle M'}
\right] \cong \calUunophiMg{M'} \longhookrightarrow \calUepsilonphiMg{M'}
\, $  and adjoint of it with respect to the quantum Poisson pairings.
\endproclaim

\demo{Proof}  Since  $ \, \calFepsilonphiMg{M} \longhookrightarrow
\gerUepsilonphiMh{M} \, $,  we can restrict  $ \gerFrhtau $  to
$ \calFepsilonphiMg{M} $,  thus obtaining a Hopf algebra morphism
$ \, \calFrHtau \colon \, \calFepsilonphiMg{M} \llongrightarrow
\gerUunophiMh{M} \cong \uhtau \, $.  But Theorem 9.3 gives
$ \, \calFunophiMg{M} = \gerUunophiMh{M} \cong \uhtau \, $,  whence
the thesis.   $ \square $
\enddemo

\vskip7pt

   We call also  $ \gerFrhtau $,  $ \calFrhtau $,  and
$ \calFrHtau $  {\bf quantum Frobenius morphisms},  because they can be
thought of as liftings of classical Frobenius morphisms to characteristic
zero.

\vskip7pt

   {\bf 7.11  Specializations of quantum Poisson pairings.}  \;  From \S\S
7.2--6 we get that the Hopf pairings
$ \, \pi^\varphi_{q,H^\tau_{\scriptscriptstyle M}} \colon \, {\gerUphiMh M' }
\otimes {\calUphiMg M } \longrightarrow \kqqm \, $, $ \, \pi^\varphi_{q,
G^\tau_{\scriptscriptstyle M}} \colon \, {\calUphiMh M } \otimes
{\gerUphiMg M' } \longrightarrow \kqqm \, $  (cf.~6.12) respectively
specialize to the natural Hopf pairings  $ \, \pi_{H^\tau_{\scriptscriptstyle
M}} \colon \, \uhtau \otimes F \left[ H^\tau_{\scriptscriptstyle M} \right]
\longrightarrow k \, $,  $ \, \pi_{G^\tau_{\scriptscriptstyle M}} \colon \,
F^\infty \left[ G^\tau_{\scriptscriptstyle M} \right] \otimes \ugtau
\longrightarrow k \, $;  in other words,
$ \, \pi^\varphi_{q,H^\tau_{\scriptscriptstyle M}} \big( {\hat h}, {\tilde g}
\big){\big\vert}_{q=1} =  \pi_{H^\tau_{\scriptscriptstyle M}} \big({\hat
h}{\big\vert}_{q=1}, {\tilde g}{\big\vert}_{q=1} \big) \, $,
$ \, \pi^\varphi_{q,G^\tau_{\scriptscriptstyle M}} \big( {\tilde h},
{\hat g} \big){\big\vert}_{q=1} =  \pi_{G^\tau_{\scriptscriptstyle M}}
\big( {\tilde h}{\big\vert}_{q=1}, {\hat g}{\big\vert}_{q=1} \big) \, $.
Thus the quantum Poisson pairing is a quantization of the classical Hopf
pairing on both our Poisson groups dual of each other.  In addition
we show that it can also be thought of as a quantization of the
classical Poisson pairing  $ \, \pi^\tau_{\Cal P} \colon \, \gerh^\tau
\otimes \gerg^\tau \rightarrow k \, $,  and of  {\sl new}  pairings between
function algebras.  We use notations  $ \, [\ \, ,\ ] := m - m^{op} \, $,
$ \, \nabla := \Delta - \Delta^{op} \, $  (superscript "$ op $" denoting
opposite operation).
                                                      \par
   First of all, we define a suitable grading on  $ {\gerUphiMg Q } $  (as a
$ \kqqm $--module)  by
  $$  deg \left(\prod_{r=N}^1 \! {\left( E^\varphi_{\alpha^r} \right)}^{(m_r)}
\! \cdot \prod_{i=1}^n \left( K^\varphi_i; 0 \atop t_i \right) \!
{\left( K^\varphi_i \right)}^{-Ent(t_i/2)} \cdot \prod_{r=1}^N \!
{\left( F^\varphi_{\alpha^r} \right)}^{(n_r)} \! \right) \! := \sum_{r=N}^1 m_r
+ \sum_{i=1}^n t_i + \sum_{r=1}^N n_r  $$
and linear extension.  Then let  $ \, \kqqm =: \gerU^\varphi_0 \subset
\gerU^\varphi_1 \subset \cdots \subset \gerU^\varphi_h \subset \cdots (\,
\subset {\gerUphiMg Q }) \, $  be the associated filtration, and set  $ \;
\partial(x) := h \; $  for all  $ \, x \in \gerU^\varphi_h \setminus
\gerU^\varphi_{h-1} \, $.  Notice that a similar notion of degree exists for
$ \ugtau $,  defined by means of the filtration  $ \, U_0 \subset U_1 \subset
\cdots \subset U_N \subset \cdots \subset \ugtau \, $  induced by the canonical
filtration of  $ T \left( \gerg^\tau \right) $  (the tensor algebra on
$ \gerg^\tau $),  and similarly for  $ \, \ugtau \otimes \ugtau \, $.  Finally
define
  $$  \pi_{q,{\Cal P}}^\varphi (h,g) := {(q-1)}^{\partial (g)} \cdot
\pi_q^\varphi (h,g)   \eqno  \forall \;\; h \in {\gerUphiMh Q }, \, g \in
\gerUphiMg{Q} \; ;  \qquad  $$
this yields a perfect pairing  $ \, \pi_{q,{\Cal P}} : {\gerUphiMh Q } \times
{\gerUphiMg Q } \longrightarrow \kqqm_{(q-1)} \, $  (the latter being the
localized ring).  In particular  $ \pi^\varphi_{q,{\Cal P}} $  can be
specialized at  $ q = 1 $.

\vskip7pt

\proclaim{Theorem 7.12}  \  $ \pi^\varphi_{q,{\Cal P}} : {\gerUphiMh Q } \times
{\gerUphiMg Q } \rightarrow \kqqm_{(q-1)} $  \  specializes to a pairing
  $$  \pi^\tau_{\Cal P} : \uhtau \times \ugtau \llongrightarrow k  $$
which extends the Lie bialgebra pairing  $ \; \pi^\tau_{\Cal P} \colon \,
\gerh^\tau \otimes \gerg^\tau \longrightarrow k \; $  (cf.~\S 1.2) and is such
that
  $$  \eqalign{
   \pi_{\Cal P} (\alpha \cdot x + \beta \cdot y, z)  &  = \alpha \cdot
\pi_{\Cal P} (x,z) + \beta \cdot \pi_{\Cal P} (y,z)  \cr
   \pi_{\Cal P} (x, \alpha \cdot u + \beta \cdot v)  &  = \alpha \cdot
\pi_{\Cal P} (x,u) + \beta \cdot \pi_{\Cal P} (x,v)  \cr
   \pi_{\Cal P} \big( x \cdot y, z \big) = \pi_{\Cal P} \big( x \otimes y,
\Delta(z) \big) \, ,  &  \qquad  \pi_{\Cal P} \big( x, z \cdot w \big) =
\pi_{\Cal P} \big( \Delta(x), z \otimes w \big)  \cr
   \pi_{\Cal P} \big( [x,y], z \big) = \pi_{\Cal P} \big( x \otimes y,
\delta(z) \big) \, ,  &  \qquad  \pi_{\Cal P} \big( x, [z,w] \big) =
\pi_{\Cal P} \big( \delta(x), z \otimes w \big)  \cr }
\eqno (7.8)  $$
for all  $ \, \alpha, \beta \in k $,  $ x, y \in \uhtau $,  $ z, w, u, v \in
\ugtau \, $  such that  $ \, \partial(\alpha \cdot u + \beta \cdot v) =
\partial(u) = \partial(v) \, $.
\endproclaim

\demo{Proof}  Let  $ \, x \in \uhtau $,  $ z \in \ugtau $,  and pick  $ \, x'
\in \gerUphiMh{Q} $,  $ z' \in \gerUphiMg{Q} $,  such that  $ \, x =
x'{\big\vert}_{q=1} $,  $ z = z'{\big\vert}_{q=1} $.  By definition,
$ \pi_{\Cal P} (x,z) $  is given by
  $$  \pi_{\Cal P} (x,z) := {\pi_{q,{\Cal P}} \big( x', z'
\big)}{\Big\vert}_{q=1} = {\Big( {(q-1)}^{\partial(z')} \cdot \pi_q \big( x', z'
\big) \Big)}{\Big\vert}_{q=1} \; ;  $$
in particular, we can select  $ x' $  and  $ z' $  such that  $ \, \partial
\big( x' \big) = \partial(x) \, $,  $ \, \partial \big( z' \big) = \partial(z)
\, $.  Now, the first two lines in (7.8) follows directly from similar
properties for  $ \, \pi_{q,{\Cal P}} \, $,  which are directly implied by
definitions.  As for the other relations in (7.8), using Leibnitz' and
co-Leibnitz' rules and identities  $ \, \partial(x \cdot y) = \partial(x) +
\partial(y) = \partial(x \otimes y) \, $  we are easily reduced to prove that
they holds for  $ \, x, y \in \gerh \, $  and  $ \, z, w \in \gerg \, $,  which
again follows from definition.  Finally to prove that  $ \pi_{\Cal P} $  is an
extension of the classical Poisson pairing a straightforward computation works.
$ \square $
\enddemo

\vskip7pt

   {\bf 7.13  The pairings  $ \, \fgtau \times \fhtau \longrightarrow k \, $,
$ \, \finfgtau \times \fhtau \longrightarrow k \, $.}  \  The construction in \S
7.11 can be reversed as follows.  Define a grading on  $ {\calUphiMg P } $  (as
a  $ \kqqm $--module)  by
  $$  deg \left(\prod_{r=N}^1 {\left( \ebar^\varphi_{\alpha^r}
\right)}^{m_r} \cdot \prod_{i=1}^n {\left( {\left( L^\varphi_i
\right)}^{\pm 1} - 1 \right)}^{l_i} \cdot \prod_{r=1}^N
{\left( \fbar^\varphi_{\alpha^r} \right)}^{n_r} \right) :=
\sum_{r=N}^1 m_r + \sum_{i=1}^n l_i + \sum_{r=1}^N n_r  $$
and linear extension; then let  $ \, \kqqm =: \calU^\varphi_0 \subset
\calU^\varphi_1 \subset \cdots \subset \calU^\varphi_h \subset \cdots (\,
\subset {\calUphiMg P }) \, $  be the associated filtration, and set  $ \;
\partial(x) := h \; $  for all  $ \, x \in \calU^\varphi_h \setminus
\calU^\varphi_{h-1} \, $  ($ h \in \N $).  Then extend  $ \, \pi^\varphi_q
\colon {\uqphiMh P } \otimes {\uqphiMg Q } \llongrightarrow \kq \, $  to a
perfect pairing of formal Hopf algebras  $ \, \pi^\varphi_q \colon {\uqphiMh P }
\otimes {\uqphiMg P } \llongrightarrow k \left( q^{1/D}, q^{-1/D} \right) \, $
(where  $ D $  is the determinant of the Cartan matrix) by the rule  $ \,
\pi^\varphi_q \left( L_\lambda, L_\mu \right) := q^{(\lambda \vert \mu)} \, $
(where  $ \, (\lambda \vert \mu) \, $  is defined in \S 1.1).  Finally define
  $$  \pi_{q,\varphi}^{\Cal P} (h,g) := {(q-1)}^{-\partial (g)} \cdot
\pi^\varphi_q(h,g)   \eqno  \forall \;\; h \in {\calUphiMh P }, \, g \in
{\calUphiMg P } \; ;  \qquad  $$
this yields a perfect pairing  $ \, \pi_{q,\varphi}^{\Cal P} : {\calUphiMh P }
\times {\calUphiMg P } \longrightarrow k \left[ q^{1/d}, q^{-{1/d}} \right]
\, $,  whose set of values is an ideal coprime with the principal ideal
$ \left( q^{1/D} - 1 \right) $;  furthermore, restriction gives also a similar
pairing  $ \, \pi_{q,\varphi}^{\Cal P} : {\gerFphiMg P } \times {\calUphiMg P }
\longrightarrow k \left[ q^{1/d}, q^{-{1/d}} \right] \, $.
                                            \par
   Now we can specialize these pairings at  $ \, q = q^{1/d} = 1 \, $,  which
gives the following:

\vskip7pt

\proclaim{Teorema 7.14}  The pairing  $ \; \pi_{q,\varphi}^{\Cal P} :
\calUphiMh{P} \times \calUphiMg{P} \llongrightarrow k \left[ q^{1/d}, q^{-{1/d}}
\right] \; $  and the pairing  $ \; \pi_{q,\varphi}^{\Cal P} : \gerFphiMg{P}
\times \calUMg{P} \llongrightarrow k \left[ q^{1/d}, q^{-{1/d}} \right] \; $
specialize to pairings
  $$  \pi_\tau^{\Cal P} : \, \finfgtau \otimes \fhtau \llongrightarrow k \, ,
\qquad  \pi_\tau^{\Cal P} : \, \fgtau \otimes \fhtau \llongrightarrow k  $$
such that
  $$  \eqalign{
   \pi_\tau^{\Cal P} (\alpha \cdot x + \beta \cdot y, z)  &  = \alpha \cdot
\pi_\tau^{\Cal P} (x,z) + \beta \cdot \pi_\tau^{\Cal P} (y,z)  \cr
   \pi_\tau^{\Cal P} (x, \alpha \cdot u + \beta \cdot v)  &  = \alpha \cdot
\pi_\tau^{\Cal P} (x,u) + \beta \cdot \pi_\tau^{\Cal P} (x,v)  \cr
   \pi_\tau^{\Cal P} \big( x \cdot y, z \big) = \pi_\tau^{\Cal P} \big( x
\otimes y, \Delta(z) \big) \, ,  &  \qquad  \pi_\tau^{\Cal P} \big( x, z \cdot w
\big) = \pi_\tau^{\Cal P} \big( \Delta(x), z \otimes w \big)  \cr
   \pi_\tau^{\Cal P} \big( \{x,y\}, z \big) = \pi_\tau^{\Cal P} \big( x \otimes
y, \nabla(z) \big) \, ,  &  \qquad  \pi_\tau^{\Cal P} \big( x, \{z,w\} \big) =
\pi_\tau^{\Cal P} \big( \nabla(x), z \otimes w \big)  \cr }   \eqno (7.9)  $$
for all  $ \, \alpha, \beta \in k $,  $ x, y \in \fgtau $  or  $ x, y \in
\finfgtau $,  $ z, u, v \in \fhtau \, $  such that  $ \, \partial(\alpha \cdot u
+ \beta \cdot v) = \partial(u) = \partial(v) \, $  (with  $ \, \partial(x) :=
\partial \left( x' \right) \, $  for any  $ \, x' \in \calUphiMg{P} \, $  such
that  $ \, x'{\big\vert}_{q=1} = x \, $).
\endproclaim

\demo{Proof}  Just mimick the proof of Theorem 7.12 above.   $ \square $
\enddemo

\vskip1,7truecm

    \centerline {\bf  \S \; 8 \,  Formal quantum groups }

\vskip10pt

   {\bf 8.1  Formal quantum groups versus quantum formal groups.}  \  The
title of this subsection  {\sl is not}  a play on words: in fact we wish to
discuss the possibility of develop  {\sl two}  different notions which are
to be quantum analog of the notion of formal group; the different position
of the word  {\sl quantum}  in the previous expressions just refer to two
different way of conceive the notion of formal group, which give rise to
two different "quantizations".
                                                     \par
   In \S 7.1 we start from the fact that a formal group is given by
a  {\sl commutative}  formal Hopf algebra, which can be realized as
$ {\ug}^* $   --- the dual of  $ \ug $ ---   thus we defined the quantum formal
groups as spectra of formal Hopf  algebras, and we looked at
$ {\uqphiMg M }^* $.
                                                     \par
   An alternative method stems from the fact that the topological Hopf algebra
of a formal group may be obtained as a suitable completion of a usual Hopf
algebra.  Namely, let  $ F^{\scriptscriptstyle \infty}[G] $  be the formal Hopf
algebra of a given formal group; let  $ G $  be an algebraic group with
associated formal group equal to the given one; let  $ {\frak m}_e $  be the
maximal ideal of  $ \fg $  associated to the identity  $ \, e \in G \, $;  then
$ \finfg $  is the  $ {\frak m}_e $--adic  completion of  $ \fg $.  Moreover we
remark that  $ \, {\frak m}_e = {\frak e} := Ker(\epsilon) \, $,  where
$ \epsilon $  is the counit of  $ \fg $.
                                                     \par
   The previous remarks motivate the following way of "quantizing"
$ \finfg $:  first, constructing a Hopf algebra  $ F_q[G] $  which
quantizes  $ \fg $;  second, constructing the  $ \gerE $--adic  completion
of  $ F_q[G] $,  with  $ \, \gerE := Ker \big( \epsilon \colon \, F_q[G]
\rightarrow \kq \big) \, $.  We shall call an object obtained in this way
{\bf formal quantum group}.  When considering formal  {\sl Poisson\/}  groups
we require also that such a quantization is one of the Poisson structure.
                                                    \par
   We have all the ingredients to perform this
construction.  The first steps are trivial.

\vskip7pt

\proclaim{Definition 8.2}  Let  $ M $  be a lattice as in \S 2.2,  and let
$ {\fqphiMg M } $,  $ {\gerFphiMg M } $,  and  $ {\calFphiMg M } $  be the
quantum function algebras defined in \S 4.
                                              \hfill\break
   \indent   Let  $ \, \hbox{\bf E}_\varphi := Ker \left( \epsilon \colon \,
{\fqphiMg M } \llongrightarrow \kq \right) \, $,  $ \, \gerE_\varphi :=
Ker \left( \epsilon \colon \, {\gerFphiMg M } \llongrightarrow \kqqm \right)
\, $,  and  $ \, \calE_\varphi := Ker \left( \epsilon \colon \, {\calFphiMg M }
\llongrightarrow \kqqm \right) \, $.  Then we define
  $$  \eqalign{
   {\fqphiMinfg M }  &  := \, \text{ $ \hbox{\bf
E}_\varphi $--adic  completion of} \; {\fqphiMg M }  \cr
   {\gerFphiMinfg M }  &  := \, \text{ $ \gerE_\varphi $--adic
completion of} \; {\gerFphiMg M }  \cr
   {\calFphiMinfg M }  &  := \, \text{ $ (q-1) \cdot \calE_\varphi $--adic
completion of} \; {\calFphiMg M } \, .  \cr }  $$
\endproclaim

\vskip7pt

\proclaim{Lemma 8.3}  Let  $ H $  be a Hopf algebra over a ring  $ R $,
let  $ {\Bbb E} $  be the kernel of the counit of  $ H $,  and let
$ \, u \in R \, $  be a non-invertible element of  $ R $.
                                             \hfill\break
   \indent   (a) \; Let  $ \widehat{H} $  be the  $ {\Bbb E} $--adic
completion of  $ H $.  There exists a unique structure of topological Hopf
algebra over  $ R $  on  $ \widehat{H} $  which extends by continuity that of  $
H $.
                                             \hfill\break
   \indent   (b) \; Let  $ \widehat{H}_u $  be the
$ u \cdot {\Bbb E} $--adic  completion of  $ H $.  There exists a unique
structure of topological Hopf algebra over  $ R $  on  $ \widehat{H}_u $
which extends by continuity that of  $ H $.   $ \square $
\endproclaim

\vskip7pt

\proclaim{Proposition 8.4}  $ {\gerFphiMinfg M } $  and
$ {\calFphiMinfg M } $  are  $ \kqqm $--integer  forms of
$ {\fqphiMinfg M } $  as topological Hopf algebras.   $ \square $
\endproclaim

\vskip7pt

\proclaim{Definition 8.5}  Let  $ \, \epsilon' \colon \, \hbox{\bf
H}_\varphi^{\scriptscriptstyle M} \longrightarrow \kq \, $  be the  $ \kq
$--algebra  morphism  defined by  $ \, \epsilon' \left( F^\varphi_i \right)
:= 0 $,  $ \, \epsilon' \left( L^\varphi_\mu \right) := 1  $,  $ \, \epsilon'
\left( E^\varphi_i \right) := 0  $,  ($ \, \forall \, i=\unon, \, \mu
\in M $)  and set  $ \, {\Bbb E}'_\varphi := Ker\left(\epsilon'\right) \, $,
$ \, {\widetilde{\Bbb E}}'_\varphi := {\Bbb E}'_\varphi \cap {\Cal
H}_\varphi^{\scriptscriptstyle M} \, $,  $ \, {\widehat{\Bbb E}}'_\varphi :=
{\Bbb E}'_\varphi \cap {\frak H}_\varphi^{\scriptscriptstyle M} \, $.  We
call  $ {\uqphiMinfh M } $  the  $ {\Bbb E}'_\varphi $--adic
completion of  $ \hbox{\bf H}_\varphi^{\scriptscriptstyle M} $,
$ {\calUphiMinfh M } $  the  $ {\widetilde{\Bbb E}}'_\varphi $--adic
completion of  $ {\Cal H}_\varphi^{\scriptscriptstyle M} $,  and
$ {\gerUphiMinfh M } $  the  $ (q-1) \cdot
{\widehat{\Bbb E}}'_\varphi $--adic  completion
of  $ {\frak H}_\varphi^{\scriptscriptstyle M} $,
with its natural structure of topological  $ \kqqm $--algebra.
\endproclaim

\vskip7pt

\proclaim{Proposition 8.6}  There exists a unique isomorphism of topological
              $ \kq $--algebras\footnote{ Of course by  {\it morphism
of topological algebras}  we mean a morphism of algebras which is  {\it
continuous}. }\break
$ \; \xi_{\scriptscriptstyle M}^{\varphi,\infty} \colon \, {\fqphiMinfg M }
\, {\buildrel \cong \over \llongrightarrow} \, {\uqphiMinfh M } \; $
which extends  $ \, \xi^\varphi_{\scriptscriptstyle M} \colon \,
{\fqphiMg M } \longhookrightarrow \hbox{\bf H}_\varphi^{\scriptscriptstyle M}
\, $  and  $ \, \xi^\varphi_{\scriptscriptstyle M} \colon \, {\fqphiMg M }
       \! \left[ \psi_{-\rho}^{-1} \right] $\break
$\, {\buildrel \cong \over \longrightarrow} \, \hbox{\bf
H}_\varphi^{\scriptscriptstyle M} \, $.  It restricts to  $ \; {\gerFphiMinfg
M } \, {\buildrel \cong \over \llongrightarrow} \, {\calUphiMinfh M } \; $
(which extends  $ \, \xi^\varphi_{\scriptscriptstyle M} \colon \,
{\gerFphiMg M } \longhookrightarrow {\Cal A}_\varphi^{\scriptscriptstyle M}
\, $  and  $ \, \xi^\varphi_{\scriptscriptstyle M} \colon \, {\gerFphiMg M }
\! \left[ \psi_{-\rho}^{-1} \right] \, {\buildrel \cong \over
\longrightarrow} \, {\Cal A}_\varphi^{\scriptscriptstyle M} \, $)  and to
$ \; {\calFphiMinfg M } \, {\buildrel \cong \over \llongrightarrow} \,
{\gerUphiMinfh M } \, $  (which extends  $ \, \xi^\varphi_{\scriptscriptstyle
M} \colon {\calFphiMg M } \longhookrightarrow
{\frak A}_\varphi^{\scriptscriptstyle M} \, $  and
$ \, \xi^\varphi_{\scriptscriptstyle M} \colon \, {\calFphiMg M }
\! \left[ \psi_{-\rho}^{-1} \right] \, {\buildrel \cong \over
\longrightarrow} \, {\frak A}_\varphi^{\scriptscriptstyle M} \, $.
Then by push-out the right-hand-side algebras get structures of topogical
Hopf algebras, so that  $ \xi_{\scriptscriptstyle M}^{\varphi,\infty} $  is
always an isomorphism of topological Hopf algebras.
\endproclaim

\demo{Proof}  Consider  {\it (a)}.  From definitions, Theorem 6.11, and
formulas for  $ \, \epsilon \colon \, {\uqphiMh M } \rightarrow \kq \, $  in
\S 6.9 it follows that  $ \, \xi^\varphi_{\scriptscriptstyle M} \left(
\hbox{\bf E}_\varphi \right) \subseteq {\Bbb E}'_\varphi \, $,  hence there
exists a unique continuous extension of
$ \xi^\varphi_{\scriptscriptstyle M} \, $,
$ \, \xi^{\varphi,\infty}_{\scriptscriptstyle M} \colon \, {\fqphiMinfg M }
\longhookrightarrow {\uqphiMinfh M } \, $,  which is a monomorphism of
topological  $ \kq $--algebras.  On the other hand
$ \, \xi^\varphi_{\scriptscriptstyle M} \colon \, {\fqphiMg M } \! \left[
\psi_{-\rho}^{-1} \right] @>\cong>> \hbox{\bf H}_\varphi^{\scriptscriptstyle
M} \, $,  with  $ \, \xi^\varphi_{\scriptscriptstyle M} \left( \psi_{-\rho}
\right) = L^\varphi_{-\rho} \, $  (cf.~the proof of Theorem 5.14);  then
$ \, \epsilon \left( 1 - \psi_{-\rho} \right) = \epsilon \big(
\xi^\varphi_{\scriptscriptstyle M} (1 - \psi_{-\rho}) \big) = \epsilon
\big( 1 - L^\varphi_{-\rho} \big) = 0 \, $,  hence
$ \, (1 - \psi_{-\rho}) \in Ker(\epsilon) =: \hbox{\bf E}_\varphi
\, $;  but then  $ \, \psi_{-\rho}^{-1} = \sum_{n=0}^{+\infty} {\left( 1 -
\psi_{-\rho} \right)}^n \in {\fqphiMinfg M } \, $,  whence  $ {\fqphiMg M }
\! \left[ \psi_{-\rho}^{-1} \right] $  canonically embeds into
$ {\fqphiMinfg M } $,  thus  $ \, \xi^\varphi_{\scriptscriptstyle M} \left(
{\fqphiMinfg M } \right) \supseteq \hbox{\bf H}_\varphi^{\scriptscriptstyle
M} \, $  and then by continuity  $ \, \xi^\varphi_{\scriptscriptstyle M}
\left( {\fqphiMinfg M } \right) = {\uqphiMinfh M } \, $,  so that  {\it (a)}
is proved.
                                                  \par
   For  {\it (b)}  and  {\it (c)}  we proceed like for  {\it (a)};  we have
only to notice, for case  {\it (c)},  that  $ \, L^\varphi_{-\rho} =
\prod_{i=1}^n L^\varphi_{-\mu_i} = \prod_{i=1}^n {\left( M^\varphi_i
\right)}^{-1} \, $,  hence  $ \, L^\varphi_\rho = \prod_{i=1}^n M^\varphi_i
\, $,  and  $ \, M^\varphi_i = \sum_{n=0}^{+\infty} {\left( 1 -
{\left( M^\varphi_i \right)}^{-1} \right)}^n = \sum_{n=0}^{+\infty}
{{(-d_i)}_q}^{\! n} \cdot {\left( q - 1 \right)}^n \cdot
{\Big( {{\left( M^\varphi_i \right)}^{-1} ; \, 0 \atop 1} \Big)}^{\! n} \, $
with  $ \, \Big( {{\left( M^\varphi_i \right)}^{-1} ; \, 0 \atop 1} \Big) \in
{\widehat{\Bbb E}}' \, $;  but  $ \, {\left( M^\varphi_i \right)}^{-1} =
\xi^\varphi_{\scriptscriptstyle M} \big( \psi_{-\mu_i} \big) \, $,  so
$ \, \Big( {{\left( M^\varphi_i \right)}^{-1} ; \, 0 \atop 1} \Big) =
\xi^\varphi_{\scriptscriptstyle M} \left( \left( \psi_{-\mu_i} ; \, 0 \atop 1
\right) \right) \, $,  with  $ \, \left( \psi_{-\mu_i} ; \, 0 \atop 1 \right)
\in \calE \, $;  then  $ \, \psi_{-\rho}^{-1} = \prod_{i=1}^n
\psi_{-\mu_i}^{-1} = \prod_{i=1}^n \sum_{n=0}^{+\infty} {{(-d_i)}_q}^{\! n}
\cdot {\left( q - 1 \right)}^n \cdot {\left( \psi_{-\mu_i} ; \, 0 \atop 1
\right)}^{\! n} \in {\calFphiMinfg M } \, $  and we conclude like
for  {\it (a)}.   $ \square $
\enddemo

\vskip7pt

\proclaim{Theorem 8.7}  The topological Hopf algebra
$ {\gerFphiMinfg M } $  specializes to  $ \, F^\infty \left[
G^\tau_{\scriptscriptstyle M} \right] = \finfgtau \, $  as topological
Poisson Hopf algebra for  $ \, q \rightarrow 1 \, $,  that is
  $$  \gerFunophiMinfg{M} := {\gerFphiMinfg M } \Big/ \! (q-1) \,
{\gerFphiMinfg M } \, \cong \, F^\infty \left[ G^\tau_{\scriptscriptstyle M}
\right] \, \cong \, {\calUphiMinfh M } \Big/ \! (q-1) \, {\calUphiMinfh M }
=: \calUunophiMinfh{M}  $$
\endproclaim

\demo{Proof}  Recall that  $ F^\infty \left[ G^\tau_{\scriptscriptstyle M}
\right] $  is the  $ {\frak e} $--adic  completion of  $ F \left[
G^\tau_{\scriptscriptstyle M} \right] $.  But  $ {\gerFphiMinfg M } $  is by
definition the  $ \gerE_\varphi $-adic  completion of  $ {\gerFphiMg M } $;
since  $ \, {\gerFphiMg M } @>{q \rightarrow 1}>> F \left[
G^\tau_{\scriptscriptstyle M} \right] \, $  as Poisson Hopf algebra
(cf.~(4.6)),  $ {\gerFphiMinfg M } $  does specialize   --- for
$ \, q \rightarrow 1 \, $  ---   to the  $ \gerE^\tau_1 $--adic
completion of  $ F \left[ G^\tau_{\scriptscriptstyle M} \right] $,
with  $ \, \gerE^\tau_1 := \gerE_\varphi{\Big\vert_{q=1}} \, $;
but  $ \gerE^\tau_1 = {\frak e} \, $,  whence the thesis.   $ \square $
\enddemo

\vskip7pt

   {\bf Remark 8.8.}  So far we found two topological Hopf algebras, that is
$ \, {\gerFphiMinfg M } = {\calUphiMinfh M } \, $  and  $ \, {\calUphiMh M }
= {\gerUphiMg M' }^* \, $,  which both contain  $ {\gerFphiMg M } $  and for
$ \, q \rightarrow 1 \, $  do specialize to the same object, namely
$ \; {\calUphiMinfh M }{\Big\vert}_{q=1} = {\gerFphiMinfg
M }{\Big\vert}_{q=1} \, \cong \, {\gerUphiMg M' }^*{\Big\vert}_{q=1} =
{\calUphiMh M }{\Big\vert}_{q=1} \; $.
                                               \par
   Now, next theorem shows that this is "singular fact", i.~e.~for
"general  $ q $"  we have
  $$  \calUphiMinfh{M} = \gerFphiMinfg{M} \, \not\cong \,
{\gerUphiMg{M'} }^* = \calUphiMh{M} \, .  $$

\vskip7pt

\proclaim{Theorem 8.9}  There does not exist any isomorphism of topological
Hopf  $ \kq $--algebras  among  $ \, {\uqphiMinfh M } = {\fqphiMinfg M }
\, $  and  $ \, {\uqphiMg {M'\!} }^* = {\uqphiMh {M\!} } \, $  whose
restriction to  $ {\fqphiMg M } $  is the identity.  Hence similar statements
hold for the integer forms too.
\endproclaim

\demo{Proof}  The second part of the claim follows from the first because
of Proposition 8.4.  Let now  $ \, \Theta \colon \, {\uqphiMinfh M } =
{\fqphiMinfg M } {\buildrel {\cong} \over \loongrightarrow}
{\uqphiMg {M'} }^* = {\uqphiMh M } \, $  be an isomorphism
of the above type; then  $ \, \Theta \left(
L^\varphi_{-\mu} \right) = L^\varphi_{-\mu} \, $  for all  $ \, \mu \in M_+
\, $.  Let  $ \, {\{a_n\}}_{n \in \N} \subseteq \kq \, $  be any sequence in
$ \kq $;  since  $ \, \big( L^\varphi_{-\mu_i} - 1 \big) = \big( {\left(
M^\varphi_i \right)}^{-1} - 1 \big) \in {\Bbb E}'_\varphi \, $,  we have
$ \, \sum_{n=0}^{+\infty} a_n {\big( {\left( M^\varphi_i \right)}^{-1} - 1
\big)}^n \in {\fqphiMinfg M } \, $;  therefore continuity implies
$ \; \Theta \left( \sum_{n=0}^{+\infty} a_n {\big( {\left( M^\varphi_i
\right)}^{-1} - 1 \big)}^n \right) = \sum_{n=0}^{+\infty} \Theta \left( a_n
{\big( {\left( M^\varphi_i \right)}^{-1}
                          - 1 \big)}^n \right) = $\break
$ \sum_{n=0}^{+\infty} a_n {\big( {\left( M^\varphi_i \right)}^{-1} - 1
\big)}^n \, ; \; $  the last should belong to  $ \,
{\uqphiMg {M'} }^* \, $,  i.~e.~it should be a linear functional on
$ \, {\uqphiMg {M'} } \, $:  but on  $ \, \Lambda_i^{-1} := L_{-\nu_i} \, $  its
value would be
  $$  {\left\langle \sum_{n=0}^{+\infty} a_n {\left( {\left( M^\varphi_i
\right)}^{-1} - 1 \right)}^n, \Lambda_i^{-1} \right\rangle}_{\!
\overline{\pi}} = \sum_{n=0}^{+\infty} a_n {\left\langle {\left( {\left(
M^\varphi_i \right)}^{-1} - 1 \right)}^n, \Lambda_i^{-1}
\right\rangle}_{\overline{\pi}} = \sum_{n=0}^{+\infty} a_n {(q_i - 1)}^n  $$
and for  {\sl general}  $ {\{a_n\}}_{n \in \N} $ the right-hand-side  {\sl is
not} an element of  $ \kq $,  contradiction.   $ \square $
\enddemo

\vskip7pt

   Thus  $ \, \calUphiMinfh{M} = \gerFphiMinfg{M} \, $  is a
quantization of  $ \finfgtau $  {\sl different}  from $ \, \calUphiMh{M}
= {\gerUphiMg M' }^* $;  so also  $ \, {\gerUphiMinfh M } =
{\calFphiMinfg M } \, $  is another quantization of  $ \uhtau $,
{\sl different}  from  $ \, {\gerUphiMh M } \, $:

\vskip7pt

\proclaim{Theorem 8.10}  For  $ \, q \rightarrow 1 \, $  the topological Hopf
algebra  $ \, {\gerUphiMinfh M } = {\calFphiMinfg M } \, $  does specialize to
$ \uhtau $  as a Poisson Hopf coalgebra,  that is
  $$  \calFunophiMinfg{M} := \, {\calFphiMinfg M } \Big/ (q-1) \,
{\calFMinfg M } \, \cong \, \uhtau \, \cong \, {\gerUphiMinfh M } \Big/ (q-1) \,
{\gerUphiMinfh M } \, =: \gerUunophiMinfh{M}  $$
\endproclaim

\demo{Proof}  From definitions follows  $ \, \calFunophiMinfg{M} =
\calFunophiMg{M} \, $  as Poisson Hopf coalgebras; but for Theorem 7.3 is
$ \, \calFunophiMg{M} \cong \uhtau \, $  (as Poisson Hopf coalgebras), whence
the claim.   $ \square $
\enddemo

\vskip7pt

   We finish with a quantum Frobenius morphism.  Let  $ \varepsilon $  and
$ \ell $  be as in \S 4.3,  and set
  $$  \gerFepsilonphiMg{M} := \, {\gerFphiMinfg M } \Big/ (q - \varepsilon)
\, {\gerFphiMinfg M } \, \cong \, {\calUphiMinfh M } \Big/ (q - \varepsilon)
\, {\calUphiMinfh M } \, =: \calUepsilonphiMinfh{M}  $$

\vskip7pt

\proclaim{Theorem 8.11}  There exists a unique monomorphism of topological
Hopf algebras
  $$  \gerFrGtau^\infty : \, F^\infty \left[ G^\tau_{\scriptscriptstyle
M} \right] \cong \calUunophiMinfh{M} = \gerFunophiMinfg{M}
\llonghookrightarrow \gerFepsilonphiMg{M} = \calUepsilonphiMinfh{M}  $$
which extends  $ \, \gerFrGtau \colon \, F \left[
G^\tau_{\scriptscriptstyle M} \right] \cong \gerFunophiMg{M}
\llonghookrightarrow \gerFepsilonphiMg{M} \, $  (cf.~(4.7)): it is defined by
  $$  \calFrgtau^\infty :  \; \quad \fbar^\varphi_\alpha \Big\vert_{q=1}
\mapsto {\left( \fbar^\varphi_\alpha \right)}^{\! \ell}
\Big\vert_{q=\varepsilon} \, , \; \quad L^\varphi_\mu \Big\vert_{q=1} \mapsto
{\left( L^\varphi_\mu \right)}^\ell \Big\vert_{q=\varepsilon} \, ,  \; \quad
\ebar^\varphi_\alpha  \Big\vert_{q=1} \mapsto {\left(
\ebar^\varphi_\alpha \right)}^{\! \ell} \Big\vert_{q=\varepsilon} \; ;  $$
its image  $ F_0^{\varphi,\infty} $  is the topological Hopf subalgebra of
$ \, \calUepsilonphiMinfh{M} \, $  topologically generated by
$ \, \Big\{\, {\big( \fbar^\varphi_\alpha \big)}^{\! \ell},
{\big( L^\varphi_\mu \big)}^\ell, {\big( \ebar^\varphi_\alpha \big)}^{\!
\ell} \,\Big\vert\, \alpha \in R^+, \mu \in M \,\Big\} \, $,  and it is
contained in the centre of  $ \gerFepsilonphiMg{M} $.
\endproclaim

\demo{Proof}  Since  $ \, \gerFrGtau \colon \, F \left[
G^\tau_{\scriptscriptstyle M} \right] \cong \gerFunophiMg{M}
\longhookrightarrow \gerFepsilonphiMg{M} \, $  is a Hopf algebra monomorphism
we have  $ \, \gerFrGtau \big( \gerE_\varphi{\big\vert}_{q=1} \big) =
\gerE_\varphi{\big\vert}_{q=\varepsilon} \, $;  but then  $ \gerFrGtau $
extends uniquely by continuity to a topological Hopf algebra monomorphism
$ \, \gerFunophiMinfg{M} \longhookrightarrow \gerFepsilonphiMg{M} \, $  that
we call  $ \gerFrGtau^\infty $.  Now both  $ \gerFrGtau^\infty $  and
$ \calFrhtau $  are continuous extensions of  $ \gerFrGtau $,  hence they
coincide on  $ \gerFunophiMg{M} $;  in particular  $ \, \calFrhtau \left(
\psi_{-\rho} \right) = \gerFrGtau \left( \psi_{-\rho} \right) \, $,  so that
$ \, \calFrhtau \left( \psi_{-\rho}^{-1} \right) = \gerFrGtau \left(
\psi_{-\rho}^{-1} \right) \, $;  therefore  $ \gerFrGtau^\infty $  and
$ \calFrhtau $  coincide on  $ \, \gerFunophiMg{M} \! \left[
\psi_{-\rho}^{-1} \right] = {\Cal H}_\varphi^{\scriptscriptstyle M} \, $,
thus from (7.6) we get the formulas above for  $ \gerFrGtau^\infty $:
these uniquely determine it because the elements
$ \fbar^\varphi_\alpha{\big\vert}_{q=1} $,
$ L^\varphi_\mu{\big\vert}_{q=1} $,
$ \ebar^\varphi_\alpha\big\vert_{q=1} $  are topological
generators of  $ \, \calUunophiMinfh{M} = \gerFunophiMinfg{M} \, $.
Then the description of  $ \, \gerFrGtau^\infty \Big( \gerFunophiMinfg{M}
\Big) = F_0^\infty \, $  is obvious, while the fact
that  $ F_0^\infty $  is contained in the centre of
$ \, \gerFepsilonphiMg{M} = \calUepsilonphiMinfh{M} \, $  easily follows
either from Theorem 7.9{\it (b)}  or from [CV-2], \S 3.3.   $ \square $
\enddemo
 \eject

\vskip1,7truecm

 \centerline{ \bf  Appendix: the case  $ \, G = SL(2,k)
\, $ }

\vskip10pt

   For  $ \, G = SL(2,k) \, $  the algebra  $ {\uqMh P } $,
resp.~$ {\uqMh Q } $,  is generated by  $ F $,  $ L^{\pm 1} $,
resp.~$  K^{\pm 1} = L^{\pm 2} $,  $ E $.  The formal Hopf
algebra structure is given by
  $$  \epsilon \left( F \right) = 0 \, ,  \quad  \epsilon \left(
L^{\pm 1} \right) = 1 \, ,  \quad  \epsilon \left( K^{\pm 1} \right) = 1 \, ,
\quad  \epsilon \left( E \right) = 0  $$
  $$  \displaylines{
   \Delta \left( F \right) = F \otimes 1 + \sum_{n=0}^\infty q^{-n}
{\left( q - \qm \right)}^{2n} \! \cdot K E^n \otimes F^{n+1} ,
\; \;  \Delta \left( L \right) = \sum_{n=0}^\infty {\left( q -
\qm \right)}^{2n} \! \cdot L E^n \otimes F^n L  \cr
   \Delta \left( L^{-1} \right) = L^{-1} \otimes L^{-1} -
{\left( q - \qm \right)}^2 \cdot L^{-1} E \otimes F L^{-1} ,
\; \;  \Delta \left( K \right) = \sum_{n=0}^\infty
{\left( q - \qm \right)}^{2n} \! \cdot K E^n \otimes F^n K  \cr
   \Delta \left( K^{-1} \right) = K^{-1} \! \otimes K^{-1} -
{\left( q - \qm \right)}^2 \! \cdot {[2]}_q \cdot K^{-1} E \otimes F K^{-1} +
{\left( q - \qm \right)}^4 \! \cdot K^{-1} E^2 \otimes F^2 K^{-1}  \cr
   \Delta \left( E \right) = 1 \otimes E + \sum_{n=0}^\infty q^{+n}
{\left( q - \qm \right)}^{2n} \cdot E^{n+1} \otimes F^n K  \cr
   S \left( F \right) = - q^{-2} \cdot \sum_{n=0}^\infty {\left( q -
\qm \right)}^{2n} \! \cdot F^{n+1} K^{-(n+1)} E^n ,
\; \;  S \left( L \right) = \sum_{n=0}^\infty {\left( q - \qm \right)}^{2n}
\! \cdot F^n K^{-(n+1)} E^n  \cr
   S \left( L^{-1} \right) = L - {\left( q - \qm \right)}^2 F L^{-1} E ,
\; \quad  S \left( K \right) = \sum_{n=0}^\infty {\left( q -
\qm \right)}^{2n} \! \cdot F^n K^{-(n+1)} E^n  \cr
   S \left( K^{-1} \right) = K - {[2]}_q \! \cdot \! {\left( q - \qm \right)}^2
\! \cdot F E + {\left( q - \qm \right)}^4 \! \cdot F^2 K^{-1} E^2  \cr
   S \left( E \right) = - q^{+2} \cdot \sum_{n=0}^\infty {\left( q -
\qm \right)}^{2n} \cdot F^n K^{-(n+1)} E^{n+1}  \cr }  $$
   \indent   In particular from this one can prove directly all the
specialization results of \S 7.
                                                  \par
   The quantum function algebra  $ {\fqMg P } = F_q^{\scriptscriptstyle P}
[{SL(2,k)}] \, $  is known (cf.~[APW], [SV]) to be generated by elements
$ \, a $,  $ b $,  $ c $,  $ d $  with relations
  $$  \displaylines{
   a b = q \, b a \; ,  \qquad  c d = q \, d c \; ,  \qquad a c = q \, c a \; ,
\qquad  b d = q \, d b  \cr
   b c = c b \; ,  \qquad  a d - d a = (q - \qm) \, b c \; ,
\qquad  a d - q \, b c = 1  \cr }  $$
with Hopf algebra structure defined by formulas
  $$  \displaylines{
   \Delta (a) = a \otimes a + b \otimes c \; ,   \qquad  \Delta (b) = a
\otimes b + b \otimes d  \cr
   \Delta (c) = c \otimes a + d \otimes c \; ,   \qquad  \Delta (d) = c
\otimes b + d \otimes d  \cr
   S(a) = d \; ,  \qquad S(b) = -q \, b \; ,  \qquad  S(c) = -\qm c \; ,
\qquad  S(d) = a  \cr
   \epsilon (a) = 1 \; ,  \qquad  \epsilon (b) = 0 \; ,  \qquad
\epsilon (c) = 0 \; ,  \qquad \epsilon (d) = 1  \cr }  $$
moreover  $ {\gerFMg P } $  is nothing but the  $ \kqqm $--subalgebra  of
$ {\fqMg P } $  generated by  $ a $,  $ b $,  $ c $,  $ d $.
                                                       \par
   The embedding  $ \, \xi_{\scriptscriptstyle P} \colon \, {\fqMg P }
\llonghookrightarrow {\uqMh P } \, $  is described by formulas
  $$  \xi_{\scriptscriptstyle P} \colon  \quad  a \mapsto L -
{\left( q - \qm \right)}^2 F L^{-1} E ,  \; b \mapsto -
\left( q - \qm \right) F L^{-1} ,  \; c \mapsto \left( q - \qm \right)
L^{-1} E ,  \; d \mapsto L^{-1} \, ;  $$
then one can check directly that this is a morphism of formal Hopf algebras.

\vskip1,7truecm

\Refs
\endRefs

\vskip7pt

\smallrm

[APW] \  H.~H.~Andersen, P.~Polo, Wen Kexin,  {\smallit Representations of
quantum algebras\/},  Invent.~Math.~{\smallbf 104} (1991), 1--59.

\vskip4pt

[CV-1] \  M.~Costantini, M.~Varagnolo,  {\smallit Quantum double and
multiparameter quantum group\/},  Comm.~in Alg.~{\smallbf 22} (1994),
6305--6321.

\vskip4pt

[CV-2] \  M.~Costantini, M.~Varagnolo, {\smallit Multiparameter Quantum
Function Algebra at Roots of 1\/},  Math. Ann.~{\smallbf 306} (1996), 759--780.

\vskip4pt

[DD] \  I.~Damiani, C.~De Concini,  {\smallit Quantum groups and Poisson
groups\/},  in W.~Baldoni, M.~Picardello (eds.),  {\smallit Representations
of Lie groups and quantum groups\/},  Longman Scientific  $ {\scriptstyle
\and} $  Technical.

\vskip4pt

[Di] \  J.~Dieudonn\'e,  {\smallit Introduction to the theory of formal
groups\/},  Pure and Applied Mathematics {\smallbf 20}, 1973.

\vskip4pt

[DKP] \  C.~De Concini, V.~G.~Kac, C.~Procesi,  {\smallit Quantum coadjoint
action\/},  Jour.~Am. Math.~Soc.~{\smallbf 5} (1992), 151--189.

\vskip4pt

[DL] \  C.~De Concini, V.~Lyubashenko,  {\smallit Quantum function algebra at
roots of 1\/},  Adv. Math.~{\smallbf 108} (1994), 205--262.

\vskip4pt

[DP] \  C.~De Concini, C.~Procesi,  {\smallit Quantum groups\/},  in
L.~Boutet de Monvel, C.~De Concini, C.~Procesi, P.~Schapira, M.~Vergne
(eds.),  {\smallit D-modules, Representation Theory, and Quantum Groups\/},
Lectures Notes in Mathematics  {\smallbf 1565},  Springer
$ {\scriptstyle \and} $  Verlag, Berlin--Heidelberg--New York, 1993.

\vskip4pt

[Dr] \  V.~G.~Drinfel'd,  {\smallit Quantum groups\/},  Proc.~ICM Berkeley 1
(1986), 789--820.

\vskip4pt

[EK-1] \  P.~Etingof, D.~Kazhdan,  {\smallit Quantization of Lie bialgebras,
I\/},  Selecta Mathematica, New Series  {\smallbf 2} (1996), 1--41.

\vskip4pt

[EK-2] \  P.~Etingof, D.~Kazhdan,  {\smallit Quantization of Poisson
algebraic groups and Poisson homogeneous spaces\/},  Preprint q-alg/9510020
(1995).

\vskip4pt

[Ga] \  F.~Gavarini,  {\smallit Quantum function agebras as quantum
enveloping algebras\/},  Preprint q-alg/9701010, to appear in Comm.~Algebra.

\vskip4pt

[LS] \  S.~Z.~Levendorskii, Ya.~S.~Soibelman,  {\smallit Algebras of
functions on compact quantum groups, Schubert cells and quantum tori\/},
Comm.~Math.~Phys.~{\smallbf 139} (1991), 141--170.

\vskip4pt

[Lu] \  G.~Lusztig,  {\smallit Quantum groups at roots of 1\/},
Geom.~Dedicata~{\smallbf 35} (1990), 89--113.

\vskip4pt

[Pa] \  P.~Papi,  {\smallit A characterization of a good ordering in a root
system\/},  Proc.~Am.~Math. Soc.~{\smallbf 120}, no.~3 (1994), 661--665.

\vskip4pt

[Re] \  N.~Reshetikin,  {\smallit Multiparameter Quantum Groups and Twisted
Quasitriangular Hopf Algebras\/},  Lett. Math.~Phys.~{\smallbf 20} (1990),
331--335.

\vskip4pt

[So] \  Ya.~S.~Soibelman,  {\smallit The algebra of functions on a compact
quantum group and its representations\/},  Leningrad Math.~J.~{\smallbf 2}
(1991), 161--178.

\vskip4pt

[SV] \  Ya.~S.~Soibelman, L.~L.~Vaksman,  {\smallit Algebra of functions on
the  quantum group SU(2)\/},  Functional Anal.~Appl.~{\smallbf 22} (1988),
170--181.

\vskip1truecm

{}

\enddocument